\documentclass[a4paper, 10pt, parskip=half, english]{article}

\usepackage{babel}
\usepackage[T1]{fontenc}
\usepackage[utf8]{inputenc}
\usepackage{lmodern}
\usepackage{amsmath,bm,mathtools,amssymb}
\usepackage{graphicx,color}
\usepackage[footnotesize]{caption}
\usepackage{booktabs}
\usepackage{multirow}
\usepackage{relsize}

\usepackage[a4paper, margin=3cm]{geometry}
\newcommand{\correctb}[1]{#1}

\usepackage[justification=justified,singlelinecheck=false]{subcaption}

\usepackage{siunitx}
\newcommand{\bigcdot}{\mathlarger{\cdot}}

\newcommand{\mR}{\mathbb{R}}

\newcommand{\GL}{G\!L}
\newcommand{\gl}{\mathfrak{gl}}

\newcommand{\SO}{S\!O}
\let\so\undefined
\newcommand{\so}{\mathfrak{so}}
\newcommand{\SE}{{S\!E}}
\newcommand{\se}{\mathfrak{se}}

\newcommand{\dexp}{\operatorname{dexp}}
\newcommand{\Exp}{\operatorname{Exp}}
\newcommand{\Log}{\operatorname{Log}}

\newcommand{\ad}[1]{\operatorname{ad}_{#1}}

\newcommand{\T}{^{\!\mathrm{T}}}

\newcommand{\diff}[1][]{\mathrm{d}#1}

\newcommand{\pd}[2]{\frac{\partial #1}{\partial #2}}
\newcommand{\td}[2]{\frac{\mathrm{d} #1}{\mathrm{d} #2}}

\newcommand{\vga}{\bm{\gamma}}

\newcommand{\vet}{\bm{\eta}}
\newcommand{\vth}{\bm{\theta}}

\newcommand{\vka}{\bm{\kappa}}

\newcommand{\vxi}{\bm{\xi}}

\newcommand{\vph}{\bm{\phi}} 

\newcommand{\vps}{\bm{\psi}}
\newcommand{\vom}{\bm{\omega}}

\newcommand{\vvep}{\bm{\varepsilon}}

\newcommand{\vTh}{\bm \Theta}

\newcommand{\vOm}{\bm \Omega}

\newcommand{\va}{\mathbf a}
\newcommand{\vb}{\mathbf b}
\newcommand{\vc}{\mathbf c}

\newcommand{\ve}{\mathbf e}
\newcommand{\vf}{\mathbf f}
\newcommand{\vg}{\mathbf g}

\newcommand{\vm}{\mathbf m}
\newcommand{\vn}{\mathbf n}

\newcommand{\vq}{\mathbf q}
\newcommand{\vr}{\mathbf r}
\newcommand{\vs}{\mathbf s}

\newcommand{\vu}{\mathbf u}
\newcommand{\vv}{\mathbf v}

\newcommand{\vx}{\mathbf x}
\newcommand{\vy}{\mathbf y}
\newcommand{\vz}{\mathbf z}

\newcommand{\vA}{\mathbf A}
\newcommand{\vB}{\mathbf B}
\newcommand{\vC}{\mathbf C}

\newcommand{\vH}{\mathbf H}
\newcommand{\vI}{\mathbf I}
\newcommand{\vJ}{\mathbf J}

\newcommand{\vL}{\mathbf L}
\newcommand{\vM}{\mathbf M}

\newcommand{\vS}{\mathbf S}
\newcommand{\vT}{\mathbf T}

\newcommand{\vX}{\mathbf X}
\newcommand{\vY}{\mathbf Y}
\newcommand{\vZ}{\mathbf Z}

\usepackage{hyperref}

\title{
	\bfseries A total Lagrangian, objective and intrinsically locking-free Petrov--Galerkin $\SE(3)$ Cosserat rod finite element formulation
}

\author{
	Jonas Harsch$^1$, Simon Sailer$^1$, Simon R. Eugster$^1$ \\[0.25cm]
	$^1$Institute for Nonlinear Mechanics \\
	University of Stuttgart \\
	Pfaffenwaldring 9, 70569 Stuttgart, Germany
}

\date{\today}

\begin{document}
\maketitle
\begin{abstract}
	Based on more than three decades of rod finite element theory, this publication combines the successful contributions found in literature and eradicates the arising drawbacks like loss of objectivity, locking, path-dependence and redundant coordinates. Specifically, the idea of interpolating the nodal orientations using relative rotation vectors, proposed by Crisfield and Jelenić in 1999, is extended to the interpolation of nodal Euclidean transformation matrices with the aid of relative twists; a strategy that arises from the $\SE(3)$-structure of the Cosserat rod kinematics. Applying a Petrov--Galerkin projection method, we propose a rod finite element formulation where the virtual displacements and rotations as well as the translational and angular velocities are interpolated instead of using the consistent variations and time-derivatives of the introduced interpolation formula. Properties such as the intrinsic absence of locking, preservation of objectivity after discretization and parametrization in terms of a minimal number of nodal unknowns are demonstrated by representative numerical examples in both statics and dynamics.
\end{abstract}
\section{Introduction}
The theory of shear-deformable (spatial) rod formulations dates back to the pioneer works of Cosserat\cite{Cosserat1909}, Timoshenko\cite{Timoshenko1921}, Reissner\cite{Reissner1981} and Simo\cite{Simo1986}. Thus, depending on the chosen literature, shear-deformable rods are called (special) Cosserat rods\cite{Antman2005}, Simo--Reissner beams, spatial Timoshenko beams, geometrically exact beams\cite{Betsch2002}, etc. Due to this ambiguity of naming conventions\cite{Reddy2020}, it is best to describe the used rod finite element formulation by the underlying kinematics only\cite{Eugster2020b, Harsch2021a}. Nonetheless, there is a wide agreement in English literature preferring the name special Cosserat rod. Thus, we adopt this notation but for simplicity suppress the prefix ``special'' in the subsequent treatment. In addition, there is a vast amount of rod finite element formulations found in literature. Introducing all of them concerning their historical development and their tiny differences or improvements is out of the scope of this publication and the interested reader is referred to the exhaustive literature survey of both shear-deformable and shear-rigid rods given by Meier et al.\cite{Meier2019}. However, we want to introduce the most important developments found in literature that can be seen as individual milestones leading to the present formulation. 

Cardona and Geradin\cite{Cardona1988} introduced the first total Lagrangian shear-deformable spatial rod finite element formulation where the nodal rotations are parametrized in terms of total rotation vectors. To prevent the well-known singularity of this parametrization, they introduced a strategy of also using the complement rotation vector. By application of a Bubnov--Galerkin projection, the virtual work contributions are discretized in space. Focusing on details of the rotation vector parametrization, the identical rod formulation is found in Ibrahimbegović et al.\cite{Ibrahimbegovic1995}. In 1999, Crisfield and Jelenić\cite{Crisfield1999} made a groundbreaking discovery in the theory of shear-deformable rod finite element formulations. All previously existing rod finite element formulations violate the objectivity requirement in the discrete approximation. Moreover, they presented a solution to the described problem by interpolating the nodal orientations using relative rotation vectors. The follow-up publication\cite{Jelenic1999} introduces an objective rod finite element formulation using incremental rotation vectors for the description of the nodal rotations and thus working in an updated Lagrangian setting. To simplify the linearization of the internal virtual work functional, the authors applied a Petrov--Galerkin\cite{Petrov1940} projection by introducing the virtual rotation expressed in the inertial frame. \correctb{This results in a non-symmetric and configuration-dependent mass matrix. A remedy for this is to use instead the virtual rotation vector expressed in the cross-section-fixed basis together with the corresponding angular velocity (again expressed in the cross-section-fixed basis). As will be discussed in the course of the subsequent treatment, this combination yields a symmetric and possibly constant mass matrix, typically met in rigid body dynamics\cite{Sailer2021}, however, at the cost of a non-symmetric stiffness matrix.}

The previously mentioned rod finite element formulations rely on the uncoupled composition approximations of the centerline points as elements of $\mR^3$ and the cross-sectional orientations as elements of the special orthogonal group $\SO(3)$ denoted by $\mR^3 \times \SO(3)$. In contrast, helicoidal approximations\cite{Borri1994} and strain-based approaches\cite{Cesarek2013, Renda2016} use a coupled interpolation of these two fields. Based on these developments, Sonneville et al.\cite{Sonneville2014} extended the idea of interpolating the nodal orientations using relative rotation vectors\cite{Crisfield1999} to the interpolation of nodal Euclidean transformation matrices using relative twists. That is, an interpolation strategy where positions and orientations are intrinsically coupled. \correctb{Exclusively working on the Lie group $\SE(3)$, the authors\cite{Sonneville2014} applied a Bubnov--Galerkin projection, i.e., the virtual displacements are given by the variation of the nodal values.}

Inspired by this interpolation, the present paper introduces a novel rod finite element formulation using only the coupled interpolation strategy of $\SE(3)$ while retaining all other favorable properties of the classical formulations. In this way, the present rod finite element combines all the important properties of the formulations in the above literature and tries to eliminate all their drawbacks. In particular, the main contributions of this paper are the following:
\begin{itemize}
	\item We present a novel total Lagrangian (thus path-independent),  intrinsically locking-free and objective rod finite element formulation parametrized by the nodal total rotation vectors and centerline points. Using the complement rotation vector in combination with a relative interpolation strategy circumvents possible occurring singularities \correctb{for element deformations in which the relative rotation angle remains below $\pi$}. Therefore, a minimal number of nodal unknowns is obtained in contrast to formulations relying on redundant coordinates\cite{Betsch2002, Eugster2014c, Harsch2021a}.
	\item Application of a Petrov--Galerkin projection results in a very simple inertial virtual work functional that ultimately leads to a symmetric mass matrix that is constant for most applications. \correctb{In contrast to Sonneville et al.\cite{Sonneville2014}, this approach neither requires the involved expressions of the $\SE(3)$-tangent operator and its inverse nor their derivatives.}
	\item Depending on the specific application, a system of ordinary differential equations or the nonlinear generalized force equilibrium is obtained. These can respectively be solved using a standard ODE solver (e.g. Runge--Kutta family, generalized-$\alpha$) or root finding algorithms (e.g. Newton--Raphson, Riks) instead of requiring highly specialized Lie group versions of them\cite{Sonneville2014}. Hence, the formulation can easily be integrated into existing flexible multibody dynamic software packages.
	\item Besides the presentation of the linearized internal forces, static and dynamic benchmark examples demonstrate the power of the new formulation. Using a two-node element -- that can be integrated by a two-point Gauss--Legendre quadrature rule -- second-order spatial convergence is achieved for both centerline and orientation.
\end{itemize}

The remainder of this paper is organized as follows. In Section~\ref{sec:cosserat_rod_theory}, the Cosserat rod theory is briefly recapitulated. The section closes with the formulation of the internal, external and inertial virtual work functionals of the shear-deformable spatial rod. Section~\ref{sec:finite_element_formulation} introduces a novel Petrov--Galerkin rod finite element formulation based on a very efficient two-node $\SE(3)$-interpolation strategy that preserves objectivity of the discretized strain measures. Introducing some bookkeeping enables the precise formulation of the discrete virtual work functionals that result in the tuples of internal, external and gyroscopic forces as well as the symmetric and possibly constant mass matrix. Further, the linearization of internal forces is given. To demonstrate the performance of the presented finite element formulation carefully selected benchmark examples are studied in Section~\ref{sec:numerical_examples}, meaning that each of this minimal set of examples demonstrates at least one affirmed property of the proposed formulation. To close the paper, Section~\ref{sec:conclusion} presents concluding remarks.

While most parts of the paper can be read with a rudimentary knowledge of matrix Lie groups, the interested reader is referred to Appendix~\ref{sec:matrix_lie_groups} for a deeper understanding of some manipulations. Therein, a concise overview of matrix Lie groups is given with all equations and references relevant to this work. In particular, the special orthogonal group $\SO(3)$ and the special Euclidean group $\SE(3)$ are examined as subgroups of $\GL(3)$ and $\GL(4)$, respectively. Similar introductions are found in literature\cite{Bullo1995, Park2005, Barfoot2014, Sonneville2014, Arnold2016}. For didactic reasons, the variation of the rod's strain measures is shown step by step in Appendix~\ref{sec:variations_strain_measures}. Since some of the derivatives of the used Lie group formulas cannot be found in literature, all required derivatives of the proposed $\SE(3)$-interpolation strategy are presented in Appendix~\ref{sec:linearization}. \correctb{Finally, Appendix~\ref{sec:discrete_preservation_properties} discusses the discrete conservation properties of the proposed rod finite element formulation. That is, the conservation of total energy, linear and angular momentum. For this purpose, another Petrov-Galerkin formulation is presented, where the virtual nodal rotations and the nodal angular velocities are expressed in the inertial frame.}
\section{Cosserat rod theory}\label{sec:cosserat_rod_theory}
\begin{figure}
	\centering
	\includegraphics{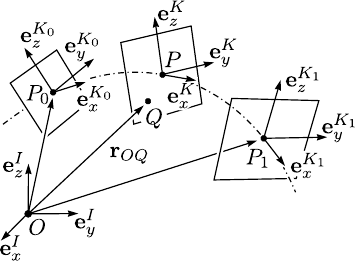}
	\caption{Kinematics of the centerline curve with their attached orthonormal basis vectors.}
	\label{fig:kinematics}
\end{figure}
\subsection{Centerline and cross-section orientation}
We introduce the Euclidean 3-space $\mathbb{E}^3$ as an abstract 3-dimensional real inner-product space\cite{Antman2005}. A basis for $\mathbb{E}^3$ is a linearly independent set of three vectors $\ve_x, \ve_y,  \ve_z \in \mathbb{E}^3$. The basis is said to be right-handed if $\ve_x \cdot (\ve_y \times \ve_z) > 0$ and orthonormal if their base vectors are mutually orthogonal and have unit length. In this paper, only right-handed orthonormal bases are considered. For a given basis $K = \{\ve_x^K, \ \ve_y^K, \ \ve_z^K\}$, the respective components $a_i^K$, $i \in \{x, y, z\}$ of a vector $\va = a_x^K \ve_x^K + a_y^K \ve_y^K + a_z^K \ve_z^K \in \mathbb{E}^3$ can be collected in the triple ${}_K \va = (a_x^K , \ a_y^K , \ a_z^K) \in \mR^3$. Thus, we carefully distinguish $\mR^3$ from the Euclidean 3-space $\mathbb{E}^3$. The rotation between two bases $K_0$ and $K_1$ is captured by the proper orthogonal transformation matrix $\vA_{K_0K_1} \in \SO(3)$, which relates the coordinate representations ${}_{K_0} \va$ and ${}_{K_1} \va$ in accordance with ${}_{K_0} \va = \vA_{K_0K_1} \, {}_{K_1} \va$. Let $\xi \in \mathcal{J} = [0, 1] \subset \mR$ denote the centerline parameter and $(\eta, \zeta) \in \mathcal{A}(\xi) \subset \mR^2$ the cross-section parameters of the cross-section area at $\xi$. Considering the rod as a three-dimensional continuum, a point $Q$ of the rod can be addressed by
\begin{equation}\label{eq:rod_kinematics}
	{}_I \vr_{OQ}(\vxi, t) = {}_I \vr_{OP}(\xi, t) + \vA_{IK}(\xi, t) \, {}_K \vr_{PQ}(\eta, \zeta) \, , \quad \vxi = (\xi, \ \eta , \ \zeta) \in \mathcal{B} = \mathcal{J} \times \mathcal{A}(\mathcal{J})\subset \mR^3 \, .
\end{equation}
Herein, ${}_I \vr_{OP}$ are the Cartesian components with respect to the inertial basis $I$ of the time-dependent centerline curve $\vr_{OP}$, where the subscripts $O$ and $P$ refer to the origin $O$ and the centerline points, see Figure~\ref{fig:kinematics}. To not overload the notation in \eqref{eq:rod_kinematics}, we abstained from explicitly denoting the dependence of the points $P$ and $Q$ on $\xi$ and $\vxi$, respectively. Moreover, ${}_K \vr_{PQ}$ denote the cross-section coordinates, which are transformed to the $I$-basis using the transformation matrices $\vA_{IK}$. Note that the base vectors of the cross-section-fixed $K$-bases are functions of the centerline parameter $\xi$ and time $t$, i.e., $\ve_i^K = \ve_i^K(\xi,t),\,i \in \{x,y,z\}$.
\subsection{Variations, velocities and curvature}
Let $\dot{(\bullet)}$ and $(\bullet)_{,\xi}$ respectively denote the derivative with respect to time $t$ and centerline parameter $\xi$. The variation is denoted by $\delta(\bullet)$. The virtual displacement ${}_I \delta \vr_{P}$ and the centerline velocity ${}_I \vv_{P}$ are given by the variation as well as the time derivative of the centerline
\begin{equation}\label{eq:virtual_displacement_velocity}
	{}_I \delta \vr_{P} = \delta \left({}_I \vr_{OP}\right) \, , \quad {}_I \vv_{P} = {}_I \dot{\vr}_{OP} \, .
\end{equation}
The angular velocity of the $K$-basis relative to the inertial $I$-basis, in components w.r.t. the $K$-basis, is defined by
\begin{equation}
	{}_K \vom_{IK} \coloneqq j^{-1}_{\SO(3)} \left( {}_K \widetilde{\vom}_{IK} \right) \, , \quad \text{with} \quad {}_K \widetilde{\vom}_{IK} \coloneqq \vA_{IK}\T \dot{\vA}_{IK} \, ,
\end{equation}
where $j_{\SO(3)} \colon \mR^3 \to \so(3) = \{\vB \in \mR^{3\times3} | \vB\T = -\vB\}$ is the linear and bijective map such that $\widetilde{\vom} \vr = j_{\SO(3)}(\vom) \vr = \vom \times \vr$ for all $\vom, \vr \in \mR^3$, see~\eqref{eq:hat_so(3)} from Appendix~\ref{sec:matrix_lie_groups}. Analogously, we define the scaled curvature as
\begin{equation}\label{eq:curvature}
	{}_K \bar{\vka}_{IK} \coloneqq j^{-1}_{\SO(3)} \left( {}_K \widetilde{\bar{\vka}}_{IK} \right) \, , \quad \text{with} \quad {}_K \widetilde{\bar{\vka}}_{IK} \coloneqq \vA_{IK}\T \vA_{IK,\xi}
\end{equation}
and the virtual rotation as
\begin{equation}\label{eq:virtual_rotation}
	{}_K \delta\vph_{IK} \coloneqq j^{-1}_{\SO(3)} \left( {}_K \delta \widetilde{\vph}_{IK} \right) \, , \quad \text{with} \quad {}_K \delta \widetilde{\vph}_{IK} \coloneqq \vA_{IK}\T \delta \vA_{IK} \, .
\end{equation}
\subsection{Reference arc length}
For the reference centerline curve ${}_I \vr_{OP}^0$, the length of the rod's tangent vector is $J = \|{}_I \vr_{OP, \xi}^0\|$. Thus, for a given centerline parameter $\xi$, the reference arc length $s$ is defined by
\begin{equation}
	s(\xi) \coloneqq \int_{0}^{\xi} J(\bar{\xi}) \, \diff[\bar{\xi}] \, .
\end{equation}
Following Harsch and Eugster\cite{Harsch2020a}, the derivative with respect to the reference arc length $s$ of a function $\vf \colon \mathcal{J} \times \mR \to \mR^3$ can be computed by
\begin{equation}\label{eq:reference_arc_length_derivative}
	\vf_{,s}(\xi,t) \coloneqq \vf_{,\xi}(\xi,t) \frac{1}{J(\xi)} \, .
\end{equation}
\subsection{Objective strain measures}
The objective strain measures of a Cosserat rod, see Antman\cite[Section 8.2]{Antman2005}, are given by
\begin{equation}\label{eq:continuous_strain_measures}
	{}_K \vka_{IK} = {}_K \bar{\vka}_{IK} / J \quad \text{and} \quad {}_K \vga = {}_K \bar{\vga} / J \, , \quad \text{with} \quad {}_K \bar{\vga} \coloneqq \vA_{IK}\T {}_I \vr_{OP, \xi} \, ,
\end{equation}
which can be gathered in the six-dimensional tuple $\vvep = ({}_K \vga , \ {}_K \vka_{IK}) \in \mR^6$. Therein, the dilatation and shear is captured by ${}_K \vga$, while ${}_K \vka_{IK}$ measures torsion and bending, see Antman\cite[Section 8.6]{Antman2005}.
\subsection{Internal virtual work}
Without loss of generality, we restrict ourselves to hyperelastic material models where the strain energy density with respect to the reference arc length $W = W({}_K \vga, {}_K \vka_{IK}; \xi)$ depends on the strain measures~\eqref{eq:continuous_strain_measures} and possibly explicitly on the centerline parameter $\xi$. By that, the internal virtual work functional is defined as
\begin{equation}\label{eq:internal_virtual_work1}
	\begin{aligned}
		\delta W^\mathrm{int} \coloneqq& -\int_{\mathcal{J}} \delta W J \diff[\xi] = -\int_{\mathcal{J}} \left\{ \delta ({}_K \bar{\vga})\T {}_K \vn + \delta ({}_K \bar{\vka}_{IK})\T {}_K \vm \right\} \diff[\xi] \, ,
	\end{aligned}
\end{equation}
where we have introduced the constitutive equations
\begin{equation}\label{eq:constitutive_equations}
	{}_K \vn \coloneqq \left(\pd{W}{{}_K \vga}\right)\T \, , \quad {}_K \vm \coloneqq \left(\pd{W}{{}_K \vka_{IK}}\right)\T \, .
\end{equation}
Note that even in the inelastic case, where no strain energy density $W$ is available, the internal virtual work~\eqref{eq:internal_virtual_work1} can be used, with internal forces and moments ${}_K \vn$ and ${}_K \vm$ given by different constitutive laws\cite{Eugster2020b}. Using $\delta ({}_K \vga)$ and $\delta ({}_K \vka_{IK})$ derived in Appendix~\ref{sec:variations_strain_measures}, the internal virtual work~\eqref{eq:internal_virtual_work1} takes the form
\begin{equation}\label{eq:internal_virtual_work2}
	\delta W^\mathrm{int} = -\int_{\mathcal{J}} \Big\{ 
	({}_I \delta \vr_{P})_{,\xi}\T \, \vA_{IK} \, {}_K \vn + ({}_K \delta \vph_{IK})_{,\xi}\T \, {}_K \vm 
	- {}_K \delta\vph_{IK}\T \left[ {}_K \bar{\vga} \times {}_K \vn + {}_K \bar{\vka}_{IK} \times {}_K \vm \right] 
	\Big\} \diff[\xi] \, .
\end{equation} 
As in Equation (2.10) of Simo and Vu-Quoc\cite{Simo1986}, we introduce the diagonal elasticity matrices $\vC_{\vga} = \operatorname{diag}(k_\mathrm{e}, k_\mathrm{s}, k_\mathrm{s})$ and $\vC_{\vka} = \operatorname{diag}(k_\mathrm{t}, k_{\mathrm{b}_y}, k_{\mathrm{b}_z})$ with constant coefficients. In the following, the simple quadratic strain energy density
\begin{equation}\label{eq:strain_energy_density}
	W({}_K \vga, {}_K \vka_{IK}; \xi) = \frac{1}{2} \left({}_K \vga - {}_K \vga^0\right)\T \vC_{\vga} \left({}_K \vga - {}_K \vga^0\right) + \frac{1}{2} \left({}_K \vka_{IK} - {}_K \vka_{IK}^0\right)\T \vC_{\vka} \left({}_K \vka_{IK} - {}_K \vka_{IK}^0\right)
\end{equation}
is used, where the superscript $0$ refers to the evaluation in the rod's reference configuration.
\subsection{External virtual work}
Assume the line distributed external forces ${}_I \vb \colon \mathcal{J} \times \mR \to \mR^3$ and moments ${}_K \vc \colon \mathcal{J} \times \mR \to \mR^3$ to be given as densities with respect to the reference arc length. Moreover, for $i\in\{0,1\}$, point forces ${}_I \vb_i \colon \mR \to \mR^3$ and point moments ${}_K \vc_i \colon \mR \to \mR^3$ can be applied to the rod's boundaries at $\xi_0=0$ and $\xi_1=1$. By that, the corresponding external virtual work functional is given by
\begin{equation}\label{eq:external_virtual_work}
	\delta W^\mathrm{ext} = \int_{\mathcal{J}} \left\{ {}_I\delta\vr_{P}\T {}_I \vb + {}_K \delta \vph_{IK}\T {}_K \vc \right\} J \diff[\xi]
	+ \sum_{i = 0}^1 \left[ {}_I\delta\vr_{P}\T {}_I \vb_i + {}_K \delta \vph_{IK}\T {}_K \vc_i \right]_{\xi_i} \, .
\end{equation}
\subsection{Inertial virtual work}
Let $\rho_0 \colon \mathcal{J} \to \mR$ denote the rod's mass density per unit reference volume and $\diff[A]$ the cross-section surface element. It is convenient to define the following abbreviations
\begin{equation}\label{eq:A_S_I_rho}
	A_{\rho_0}(\xi) \coloneqq \int_{\mathcal{A}(\xi)} \rho_0 \, \diff[A] \, , \quad 
	{}_K \vS_{\rho_0}(\xi) \coloneqq \int_{\mathcal{A}(\xi)} \rho_0 {}_K \widetilde{\vr}_{PQ} \, \diff[A] \, , \quad
	{}_K \vI_{\rho_0}(\xi) \coloneqq \int_{\mathcal{A}(\xi)} \rho_0 {}_K \widetilde{\vr}_{PQ} \, {}_K \widetilde{\vr}_{PQ}^{\,\mathrm{T}} \diff[A] \, .
\end{equation}
Further, using the mass differential $\diff[m] = \rho_0 J \diff[A] \diff[\xi]$ and expressing the variation and second time derivative of $\vr_{OQ}$ in terms of the rod's kinematics \eqref{eq:rod_kinematics}, the inertial virtual work functional of the Cosserat rod can be written as (cf. Equation (9.54) of Eugster and Harsch\cite{Eugster2020b} for a coordinate-free version)
\begin{equation}\label{eq:inertia_virtual_work}
	\begin{aligned}
		\delta W^\mathrm{dyn} &= -\int_{\mathcal{B}} {}_I\delta \vr_{OQ}\T {}_I\ddot{\vr}_{OQ} \, \diff[m]= -\int_{\mathcal{J}}
		\begin{pmatrix}
			{}_I \delta \vr_{P} \\
			{}_K \delta \vph_{IK}
		\end{pmatrix}\T
		\left\{
		\overline{\vM}
				\begin{pmatrix}
					{}_I \vv_{P} \\
					{}_K \vom_{IK}
				\end{pmatrix}^{\bigcdot}
		+ \vg
		\right\} J \diff[\xi] \, ,
	\end{aligned}
\end{equation}
where we have introduced the two quantities
\begin{equation}\label{eq:M_bar_and_g}
	\overline{\vM} = 
	\begin{pmatrix}
		A_{\rho_0} \mathbf{1}_{3 \times 3} & \vA_{IK} {}_K \vS_{\rho_0}\T \\
		{}_K \vS_{\rho_0} \vA_{IK}\T & {}_K \vI_{\rho_0}
	\end{pmatrix} \, , \quad
	\vg = \begin{pmatrix}
		\vA_{IK} {}_K \widetilde{\vom}_{IK} {}_K \vS_{\rho_0}\T {}_K \vom_{IK} \\
		{}_K \widetilde{\vom}_{IK} {}_K \vI_{\rho_0} {}_K \vom_{IK}
	\end{pmatrix} \, .
\end{equation}
Note that the quantity $\vS_{\rho_0}$ vanishes if the centerline points ${}_I \vr_{OP}(\xi, t)$ coincide with the cross-sections' center of mass. For this case $\overline{\vM}$ and $\vg$ get independent of the rod's configuration.
\section{Finite element formulation}\label{sec:finite_element_formulation}
\subsection{Rod kinematics in homogenous coordinates}
The frame $\mathcal{I} = \{O, \ve_x^I, \ve_y^I,  \ve_z^I\}$ is the set collecting the origin $O$ together with the inertial base vectors $\ve_i^I$, $i \in \{x, y, z\}$. A generic and possibly non-inertial frame $\mathcal{K} = \{P, \ve_x^K, \ve_y^K, \ve_z^K\}$ is given by a point $P$ together with a right-handed orthonormal basis $K$. The point $Q$ relative to the $\mathcal{K}$-frame is addressed by the triple ${}_{K} \vr_{P Q} \in \mR^3$, which are the components of the vector $\vr_{PQ} \in \mathbb{E}^3$ between $P$ and $Q$ with respect to the $K$-basis. Rigid body motions between two frames $\mathcal{K}_0$ and $\mathcal{K}_1$ are captured by the Euclidean transformation matrix
\begin{equation}\label{eq:homogenous_transformation}
	\vH_{\mathcal{K}_0\mathcal{K}_1} = \begin{pmatrix}
		\vA_{K_0K_1} & {}_{K_0} \vr_{P_0P_1} \\
		\mathbf{0}_{1 \times 3} & 1
	\end{pmatrix} \, ,
\end{equation}
which relates the triples ${}_{K_1} \vr_{P_1 Q} \in \mR^3$ and ${}_{K_0} \vr_{P_0 Q} \in \mR^3$ in accordance with
\begin{equation}\label{eq:Euclidean_transformation}
	\begin{pmatrix}
		{}_{K_0} \vr_{P_0Q} \\
		1
	\end{pmatrix} =
	\begin{pmatrix}
		{}_{K_0} \vr_{P_0P_1} + \vA_{K_0K_1} \, {}_{K_1} \vr_{P_1Q} \\
		1
	\end{pmatrix} =
	\vH_{\mathcal{K}_0\mathcal{K}_1}
	\begin{pmatrix}
		{}_{K_1} \vr_{P_1Q} \\
		1
	\end{pmatrix} \, .
\end{equation}
In fact, the Euclidean transformation matrix $\vH_{\mathcal{K}_0\mathcal{K}_1}$ is an element of the special Euclidean group $\SE(3)$ considered here as a Lie subgroup of the general linear group $\GL(4)$. A brief introduction to the Lie group setting and particularly an overview of the herein required mappings is given in Appendix~\ref{sec:matrix_lie_groups}. Direct computation readily verifies that the inverse of $\vH_{\mathcal{K}_0\mathcal{K}_1}$ is
\begin{equation}\label{eq:homogenous_transformation}
	\vH_{\mathcal{K}_0\mathcal{K}_1}^{-1} = \begin{pmatrix}
		\vA_{K_0K_1}\T & - \vA_{K_0K_1}\T {}_{K_0} \vr_{P_0P_1} \\
		\mathbf{0}_{1 \times 3} & 1
	\end{pmatrix} \, .
\end{equation}
Consequently, using~\eqref{eq:Euclidean_transformation}, the motion of the rod~\eqref{eq:rod_kinematics} can be written in homogenous coordinates as
\begin{equation}
	\begin{pmatrix}
		{}_{I} \vr_{OQ} \\
		1
	\end{pmatrix} =
	\vH_{\mathcal{I}\mathcal{K}}
	\begin{pmatrix}
		{}_{K} \vr_{PQ} \\
		1
	\end{pmatrix} \, , \quad \text{with} \quad \vH_{\mathcal{I}\mathcal{K}} = 
	\begin{pmatrix}
		\vA_{IK} & {}_I\vr_{OP} \\
		\mathbf{0}_{1 \times 3} & 1
	\end{pmatrix} \, .
\end{equation}
The group structure of $\SE(3)$ then allows to decompose the Euclidean transformation matrix $\vH_{\mathcal{I}\mathcal{K}}$ into
\begin{equation}\label{eq:decomposition_Euclidean_transformation}
	\vH_{\mathcal{I}\mathcal{K}} =
	\vH_{\mathcal{I}\mathcal{K}_0} \vH_{\mathcal{K}_0\mathcal{K}} \, .
\end{equation}
\subsection{$\SE(3)$-interpolation}\label{sec:SE(3)-interpolation}
In order to introduce the idea of the $\SE(3)$-interpolation strategy, we discretize the rod by a single two-node finite element. In contrast to Sonneville et al.\cite{Sonneville2014}, we choose a parametrization of $\vH_{\mathcal{I}\mathcal{K}_0}$ and $\vH_{\mathcal{I}\mathcal{K}_1}$ in terms of time-dependent generalized coordinates $\vq(t) = \big({}_I \vr_{OP_0}(t), \vps_0(t), {}_I \vr_{OP_1}(t), \vps_1(t)\big) \in \mR^{12}$ given by the nodal positions ${}_I \vr_{OP_0}(t) \in \mR^3$, ${}_I \vr_{OP_1}(t) \in \mR^3$ and the nodal rotation vectors $\vps_0(t) \in \mR^3$, $\vps_1(t) \in \mR^3$. Using the exponential map of $\SO(3)$, defined in~\eqref{eq:ExpLog_SO3}, the parametrization is
\begin{equation}\label{eq:parametrization_SE3}
	\vH_{\mathcal{I}\mathcal{K}_0}(\vq) = 
	\begin{pmatrix}
		\Exp_{ \SO(3)}(\vps_0) & {}_I \vr_{OP_0} \\
		\mathbf{0}_{3 \times 1} & 1
	\end{pmatrix} \quad \text{and} \quad
	\vH_{\mathcal{I}\mathcal{K}_1}(\vq) = 
	\begin{pmatrix}
		\Exp_{ \SO(3)}(\vps_1) & {}_I \vr_{OP_1} \\
		\mathbf{0}_{3 \times 1} & 1
	\end{pmatrix} \, .
\end{equation}
With the  $\SE(3)$-logarithm map from~\eqref{eq:ExpLog_SE3}, we can compute the relative twist $\vth_{\mathcal{K}_0\mathcal{K}_1}$ corresponding to the relative Euclidean transformation $\vH_{\mathcal{K}_0\mathcal{K}_1} = \vH_{\mathcal{I}\mathcal{K}_0}^{-1} \, \vH_{\mathcal{I}\mathcal{K}_1}$ as
\begin{equation}
	\vth_{\mathcal{K}_0\mathcal{K}_1} =
	\Log_{\SE(3)}\big(\vH_{\mathcal{K}_0\mathcal{K}_1}\big) =
	\begin{pmatrix}
		\vT_{\SO(3)}^{-\mathrm{T}}(\vps_{01}) \, {}_{K_0} \vr_{P_0P_1} \\
		\vps_{01}
	\end{pmatrix} \, .
\end{equation}
Herein $\vps_{01}$ corresponds to the relative rotation vector parameterizing the transformation between the $K_0$- and $K_1$-basis according to $\vA_{K_0K_1} = \Exp_{\SO(3)} (\vps_{01})$. \correctb{Due to the singularity of $\Log_{\SO(3)}$ for $\|\vps_{01}\|=\pi$, see Appendix~\ref{sec:SO3}, this interpolation strategy is restricted to applications in which the relative rotation angle satisfies $\omega=\|\vps_{01}\| < \pi$. A discretization with a higher number of elements always cures this problem.} By linearly scaling the relative twist $\vth_{\mathcal{K}_0\mathcal{K}_1}$ and using the $\SE(3)$-exponential map~\eqref{eq:ExpLog_SE3}, a relative Euclidean transformation from the $\mathcal{K}_0$-frame to the $\mathcal{K}$-frame can be constructed as
\begin{equation}\label{eq:linear_relative_Euclidean_transformation}
	\vH_{\mathcal{K}_0\mathcal{K}}(\xi, \vq) = \Exp_{\SE(3)} \big(\xi \, \vth_{\mathcal{K}_0\mathcal{K}_1}(\vq)\big) \, .
\end{equation}
The Euclidean transformation $\vH_{\mathcal{K}_0\mathcal{K}}$ satisfies $\vH_{\mathcal{K}_0\mathcal{K}}(0, \vq) = \mathbf{1}_{4 \times 4}$ and $\vH_{\mathcal{K}_0\mathcal{K}}(1, \vq) = \vH_{\mathcal{K}_0\mathcal{K}_1}$. To obtain the Euclidean transformation for a point within the finite element, the expressions \eqref{eq:decomposition_Euclidean_transformation} and~\eqref{eq:linear_relative_Euclidean_transformation} motivate the following $\SE(3)$-interpolation
\begin{equation}\label{eq:SE(3)_interpolation}
	\vH_{\mathcal{I}\mathcal{K}}(\xi, \vq) \coloneqq \vH_{\mathcal{I}\mathcal{K}_0}(\vq) \vH_{\mathcal{K}_0\mathcal{K}}(\xi, \vq) = \vH_{\mathcal{I}\mathcal{K}_0}(\vq) \Exp_{\SE(3)} \left(\xi \, \vth_{\mathcal{K}_0\mathcal{K}_1}(\vq)\right) \, ,
\end{equation}
originally proposed in Equation (55) of Sonneville et al.\cite{Sonneville2014}. In order to see the interpolation for the position and orientation, we explicitly compute the exponential map in~\eqref{eq:SE(3)_interpolation} resulting in
\begin{equation}
	\begin{aligned}
		\vH_{\mathcal{I}\mathcal{K}} &=
		\begin{pmatrix}
			\vA_{IK_0} & {}_I \vr_{OP_0} \\
			\mathbf{0}_{1 \times 3} & 1
		\end{pmatrix}
		\begin{pmatrix}
			\Exp_{\SO(3)}(\xi \vps_{01}) & \quad \vT_{\SO(3)}\T(\xi \vps_{01}) \, \xi \, \vT_{\SO(3)}^{-\mathrm{T}} (\vps_{01}) {}_{K_0} \vr_{P_0 P_1} \\
			\mathbf{0}_{1 \times 3} & 1
		\end{pmatrix} \\
		\begin{pmatrix}
			\vA_{IK} & {}_I \vr_{OP} \\
			\mathbf{0}_{1 \times 3} & 1
		\end{pmatrix}&=
		\begin{pmatrix}
			\vA_{IK_0} \Exp_{\SO(3)}(\xi \vps_{01}) & \quad {}_I \vr_{OP_0} + \xi \, \vA_{IK_0} \vT_{\SO(3)}\T(\xi \vps_{01}) \vT_{\SO(3)}^{-\mathrm{T}} (\vps_{01}) {}_{K_0} \vr_{P_0 P_1} \\
			\mathbf{0}_{1 \times 3} & 1
		\end{pmatrix} \, .
	\end{aligned}
\end{equation}
Consequently, the rod's orientation is discretized by
\begin{equation}\label{eq:discretization_relative_rotations}
	\vA_{IK}(\xi, \vq) = \vA_{IK_0}(\vq) \Exp_{\SO(3)}\big(\xi \vps_{01}(\vq) \big) \, ,
\end{equation}
which corresponds to Equation (4.7) of Crisfield and Jeleni{\'c}\cite{Crisfield1999}. Since $\vA_{IK}(0, \vq) = \vA_{IK_0}$ and $\vA_{IK}(1, \vq) = \vA_{IK_1}$, the discretization~\eqref{eq:discretization_relative_rotations} is a highly nonlinear interpolation of the nodal transformation matrices. The rod's centerline is discretized by
\begin{equation}
	{}_I \vr_{OP}(\xi, \vq) = {}_I \vr_{OP_0} + \xi \, \vA_{IK_0}(\vq) \vT_{\SO(3)}\T \big( \xi \vps_{01}(\vq) \big) \vT_{\SO(3)}^{-\mathrm{T}} \big( \vps_{01}(\vq) \big) {}_{K_0} \vr_{P_0 P_1}(\vq) \, ,
\end{equation}
which also interpolates the nodal positions ${}_I \vr_{OP}(0, \vq) = {}_I \vr_{OP_0}$ and ${}_I \vr_{OP}(1, \vq) = {}_I \vr_{OP_1}$ in a highly nonlinear manner.
\subsection{Objectivity, discretized strain measures and absence of locking of $\SE(3)$-interpolation}
Using the just introduced Euclidean transformation matrices, we can recognize the strain measures of the Cosserat rod theory \eqref{eq:continuous_strain_measures} in the following equation
\begin{equation}\label{eq:strain_measures}
	\vH_{\mathcal{I}\mathcal{K}}^{-1} \vH_{\mathcal{I}\mathcal{K},s}
	= 
	\begin{pmatrix}
		\vA_{IK}\T & -\vA_{IK}\T {}_I\vr_{OP} \\
		\mathbf{0}_{1 \times 3} & 1
	\end{pmatrix}
	\begin{pmatrix}
		\vA_{IK, s} & {}_I \vr_{OP, s} \\
		\mathbf{0}_{1 \times 3} & 0
	\end{pmatrix}
	= \begin{pmatrix}
		\vA_{IK}\T \vA_{IK, s} & \vA_{IK}\T {}_I \vr_{OP, s} \\
		\mathbf{0}_{1 \times 3} & 0			
	\end{pmatrix} 
	=
	\begin{pmatrix}
		{}_K \widetilde{\vka}_{IK} & {}_K \vga \\
		\mathbf{0}_{1 \times 3} & 0
	\end{pmatrix}
	\, .
\end{equation}
With the proposed interpolation~\eqref{eq:SE(3)_interpolation}, the strain measures in
\begin{equation}\label{eq:relative_strain_measures}
	\vH_{\mathcal{I} \mathcal{K}}^{-1} \vH_{\mathcal{I} \mathcal{K}, s} = \big( \vH_{\mathcal{I} \mathcal{K}_0} \vH_{\mathcal{K}_0 \mathcal{K}} \big)^{-1} \big( \vH_{\mathcal{I} \mathcal{K}_0} \vH_{\mathcal{K}_0 \mathcal{K}} \big)_{, s} = \vH_{\mathcal{K}_0 \mathcal{K}}^{-1} \vH_{\mathcal{K}_0 \mathcal{K}, s}
\end{equation}
depend only on the relative Euclidean transformation $\vH_{\mathcal{K}_0 \mathcal{K}}$ and its reference arc length derivative. Objectivity is the requirement that a quantity remains unaltered under a change of observer, i.e., a change of the $\mathcal{I}$-frame. Whether we use an $\mathcal{I}$-frame or an $\mathcal{I}^{+}$-frame, which are related by the time-dependent Euclidean transformation $\vH_{\mathcal{I}^{+} \mathcal{I}}$, the discretized strain measures in $\vH_{\mathcal{I}^{+} \mathcal{K}}^{-1} \vH_{\mathcal{I}^{+} \mathcal{K}, s}$ correspond to the strain measures obtained by~\eqref{eq:relative_strain_measures}. This proofs objectivity of the interpolation~\eqref{eq:SE(3)_interpolation}.

Using a little bit of Lie algebra introduced in the Appendix~\ref{sec:matrix_lie_groups}, the discretized strain measures can be drastically simplified. First, with the aid of~\eqref{eq:relative_strain_measures},~\eqref{eq:reference_arc_length_derivative} and~\eqref{eq:hat_se(3)}, they can be extracted from~\eqref{eq:strain_measures}, by
\begin{equation}
	\vvep = 
	\begin{pmatrix}
		{}_K \vga \\
		{}_K \vka_{IK}
	\end{pmatrix}
	= j_{\SE(3)}^{-1}\left(
	\vH_{\mathcal{K}_0 \mathcal{K}}^{-1} \vH_{\mathcal{K}_0 \mathcal{K}, \xi}
	\right) \frac{1}{J} \, .
\end{equation}
Inserting the relative interpolation~\eqref{eq:linear_relative_Euclidean_transformation} together with the definitions~\eqref{eq:ExpLog_SE3} and suppressing for a while the subscripts indicating $\SE(3)$, the discretized strain measures are
\begin{equation}
	\begin{aligned}
		\vvep(\xi, \cdot) &= 
		j^{-1} \circ \left( \Exp(\xi \vth_{\mathcal{K}_0\mathcal{K}_1})^{-1} \td{}{\xi} \Exp(\xi \vth_{\mathcal{K}_0\mathcal{K}_1}) \right) \frac{1}{J} 
		= j^{-1} \circ \left( \exp\big(j(\xi \vth_{\mathcal{K}_0\mathcal{K}_1})\big)^{-1} \td{}{\xi} \exp\big(j(\xi \vth_{\mathcal{K}_0\mathcal{K}_1})\big) \right) \frac{1}{J} \\
		&= j^{-1} \circ \left( \dexp_{-j(\xi \vth_{\mathcal{K}_0\mathcal{K}_1})} \Big(\td{}{\xi} j(\xi \vth_{\mathcal{K}_0\mathcal{K}_1})\Big) \right) \frac{1}{J} \, ,
	\end{aligned}
\end{equation}
where we have used the left-trivialized differential~\eqref{eq:left_and_right_differential}$_1$. Due to the linearity of $j$, we have $\diff / \diff[\xi]\big(j(\xi \vth_{\mathcal{K}_0\mathcal{K}_1})\big) = j(\vth_{\mathcal{K}_0\mathcal{K}_1})$. Since $j(\xi \vth_{\mathcal{K}_0\mathcal{K}_1})$ and $j(\vth_{\mathcal{K}_0\mathcal{K}_1})$ commute, according to the comment below~\eqref{eq:box}, the discretized strain measures simplify further to
\begin{equation}\label{eq:global_strain_measures_SE(3)}
	\vvep(\xi, \vq) = 
	\left( j^{-1} \circ \dexp_{-j(\xi \vth_{\mathcal{K}_0\mathcal{K}_1}(\vq))} \circ j \right) \vth_{\mathcal{K}_0\mathcal{K}_1}(\vq) \frac{1}{J}
	= \vth_{\mathcal{K}_0\mathcal{K}_1}(\vq) \frac{1}{J} \, .
\end{equation}
This means that the discretized strains are constant. 

\begin{figure}
	\centering
	\begin{subfigure}[b]{0.31\textwidth}
		\caption{}
		\includegraphics[width=\textwidth, trim=15cm 17cm 15cm 17cm, clip]{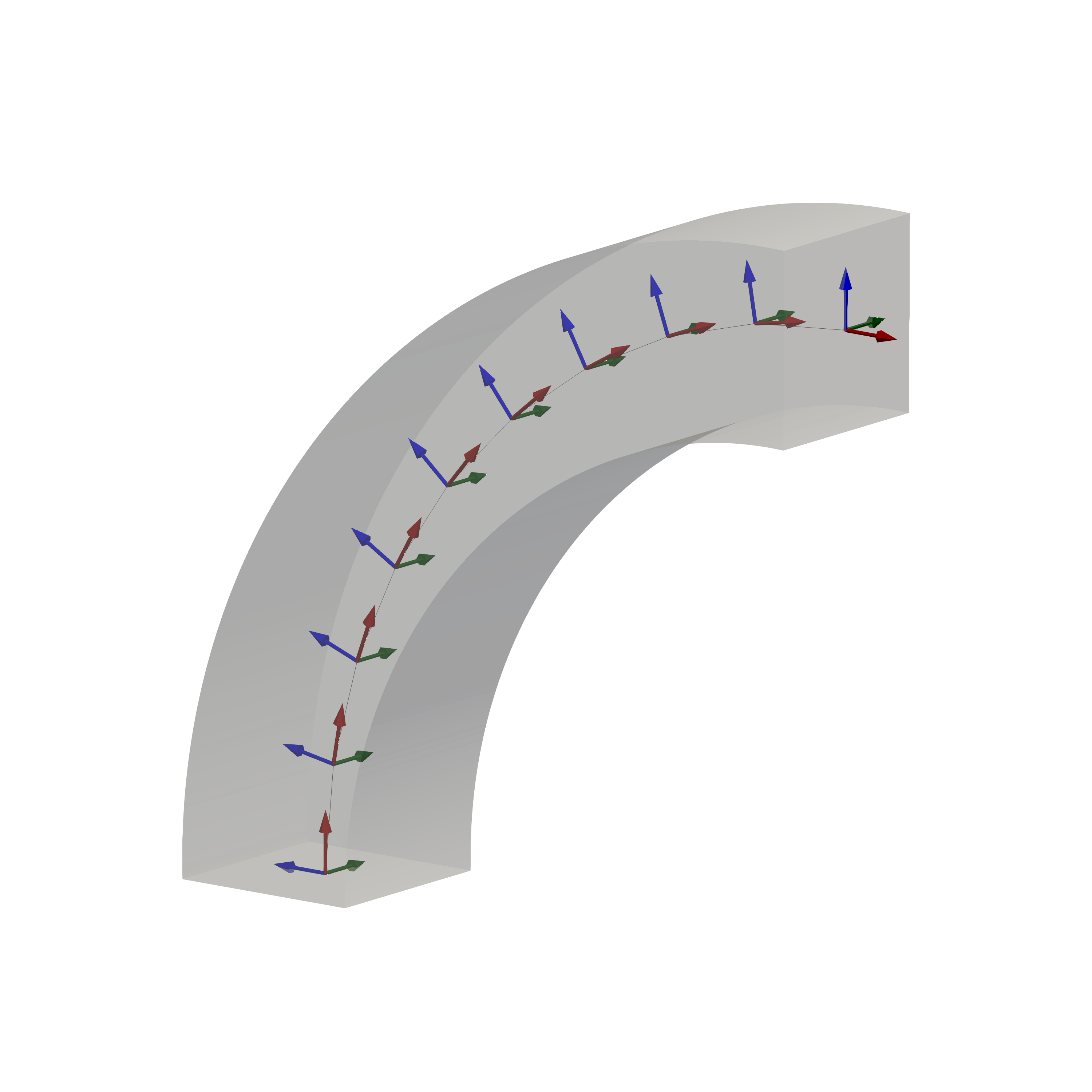}
	\end{subfigure}
	\hfill
	\begin{subfigure}[b]{0.31\textwidth}
		\caption{}
		\includegraphics{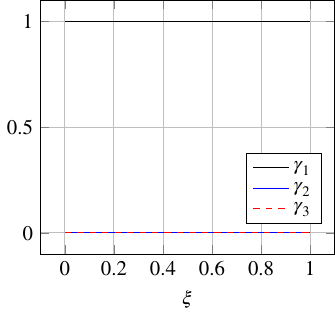}
	\end{subfigure}
	\hfill
	\begin{subfigure}[b]{0.31\textwidth}
		\caption{}
		\includegraphics{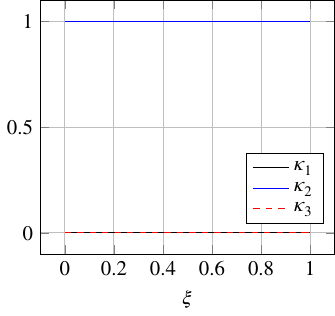}
	\end{subfigure}
	\begin{subfigure}[b]{0.31\textwidth}
		\caption{}
		\includegraphics[width=\textwidth, trim=15cm 17cm 15cm 17cm, clip]{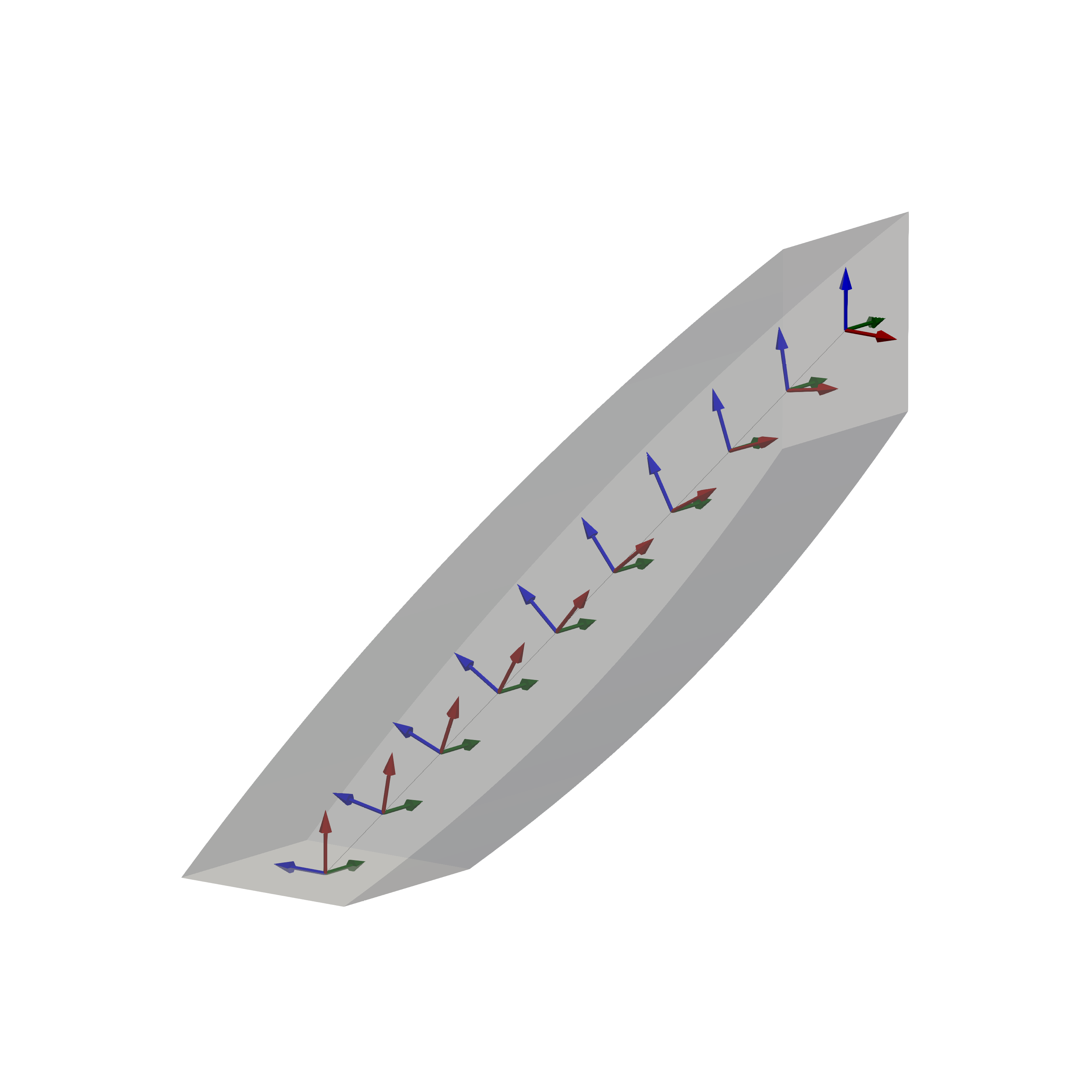}
	\end{subfigure}
	\hfill
	\begin{subfigure}[b]{0.31\textwidth}
		\caption{}
		\includegraphics{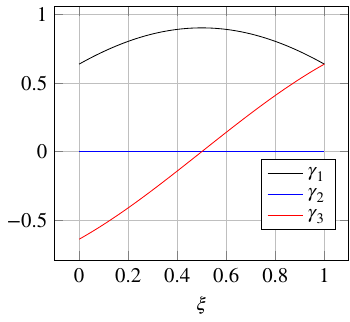}
	\end{subfigure}
	\hfill
	\begin{subfigure}[b]{0.31\textwidth}
		\caption{}
		\includegraphics{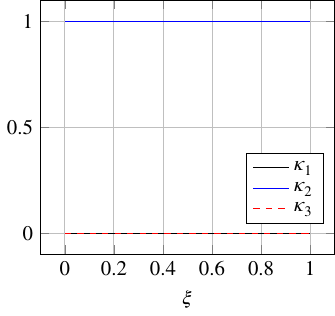}
	\end{subfigure}
	\caption{\correctb{Quarter circle in the $\ve_x^I$-$\ve_z^I$-plane defined by $\vps_0 = (0, \ -\pi / 2, \ 0)$, $\vps_1 = (0, \ 0, \ 0)$, ${}_I \vr_{OP_0} = (0, \ 0, \ 0)$ and ${}_I \vr_{OP_1} = (1, \ 0, \ 1)$. Comparison between $\SE(3)$-interpolation and $\mR^3 \times \SO(3)$-interpolation. (a) - (c) Visualization and strain measures of the two-node $\SE(3)$-element. (d) - (f) Visualization and strain measures of the two-node $\mR^3 \times \SO(3)$-element.}}
	\label{fig:Crisfield_vs_SE3}
\end{figure}

\correctb{As discussed in Balobanov and Niiranen\cite{Balobanov2018}, shear and membrane locking can occur if Kirchhoff and inextensibility constraints follow in the limit case of a parameter tending to zero.
	This appears for instance for the strain energy density \eqref{eq:strain_energy_density}, if the stiffnesses are computed in the sense of Saint--Venant by using the material's Young's and shear moduli $E$ and $G$, respectively, together with the cross-section geometry. For the particular choice of a quadratic cross-section (width $w$, area $A=w^2$, second moment of area $I =w^4/12$), the stiffnesses would be given as $k_\mathrm{e} = EA$, $k_\mathrm{s} = GA$, $k_{\mathrm{b}_y} = k_{\mathrm{b}_z} = EI$ and $k_\mathrm{t} = 2GI$. Dividing the strain energy density \eqref{eq:strain_energy_density} by $w^4$, in the limit of $w \to 0$, the scaled axial $k_\mathrm{e}/w^4$ and shear stiffnesses $k_\mathrm{s}/w^4$ tend to infinity. The strain energy remains only bounded if the dilatation remains one and the shear  deformations zero, i.e., if ${}_K \vga = (1, \ 0, \ 0)$; inextensibility $\|{}_K \vga\| - 1 = 0$ readily follows. Formulations that are prone to locking cannot fulfill these constraints exactly over the entire element and introduce parasitic dilatation and shear strains, see Figure~\ref{fig:Crisfield_vs_SE3} (e) for the two-node $\mR^3 \times \SO(3)$-interpolation\cite{Crisfield1999} in a configuration where the nodal coordinates are given by the pure bending situation of a quarter circle. In contrast, according to \eqref{eq:global_strain_measures_SE(3)}, the proposed $\SE(3)$-interpolation can exactly represent constant strain measures for all the individual contributions, can thus satisfy ${}_K \vga = (1, \ 0, \ 0)$, and is hence intrinsically relieved from both shear as well as membrane locking. Obviously, the $\SE(3)$-interpolation can represent the quarter circle exactly, see Figure~\ref{fig:Crisfield_vs_SE3} (a).}
\subsection{Kinematics of element nodes}
The concepts of the previous paragraphs can be synthesized in a rod finite element with a piecewise two-node $\SE(3)$-interpolation that results in constant strains within an element. Let the rod be discretized by $n_\mathrm{el}$ two-node elements. The nodal Euclidean transformation matrices $\vH_{\mathcal{I} \mathcal{K}_i}$, $i = 0, \dots, n_\mathrm{el}$ are parameterized in terms of the time-dependent nodal generalized coordinates $\vq_i(t) = \big({}_I \vr_{OP_i}(t), \vps_i(t)\big) \in \mR^6$ given by the nodal centerline points ${}_I \vr_{OP_i}(t) \in \mR^3$ and the nodal total rotation vectors $\vps_i(t) \in \mR^3$. Using the exponential map of $\SO(3)$, see comment before~\eqref{eq:parametrization_SE3}, the nodal parametrization is
\begin{equation}
	\vH_{\mathcal{I}\mathcal{K}_i}(\vq_i) = 
	\begin{pmatrix}
		\Exp_{ \SO(3)}(\vps_i) & {}_I \vr_{OP_i} \\
		\mathbf{0}_{3 \times 1} & 1
	\end{pmatrix} \, .
\end{equation}
The rod's parameter space $\mathcal{J}$ is divided into $n_\mathrm{el}$ linearly spaced element intervals $\mathcal{J}^e = [\xi^{e}, \xi^{e+1})$ via $\mathcal{J} = \bigcup_{e=0}^{n_\mathrm{el}-1} \mathcal{J}^e$. This means a totality of $n_\mathrm{el} + 1$ nodes and $n_{\vq} = 6(n_\mathrm{el} + 1)$ unknowns. The nodal quantities can be assembled in the global tuple of generalized coordinates $\vq(t) = \big(\vq_0(t), \vq_1(t), \dots, \vq_{N-1}(t)\big) \in \mR^{n_{\vq}}$. Introducing an appropriate Boolean connectivity matrix $\vC_e \in \mR^{12 \times n_{\vq}}$, the element generalized coordinates $\vq^e(t) = \big({}_I \vr_{OP_e}(t), \vps_e(t), {}_I \vr_{OP_{e+1}}(t), \vps_{e+1}(t)\big) \in \mR^{12}$ can be extracted from the global generalized coordinates $\vq$ via $\vq^e = \vC_{e} \vq$. By introducing an additional set of Boolean connectivity matrices $\vC_{\vr, 0}, \vC_{\vr, 1}, \vC_{\vps, 0}, \vC_{\vps, 1} \in \mR^{3 \times 12}$, the centerline points ${}_I \vr_{OP_e}$, ${}_I \vr_{OP_{e+1}}$ and total rotation vectors $\vps_{e}$, $\vps_{e+1}$ can be extracted from the element generalized coordinates $\vq^e$ via ${}_I \vr_{OP_e} = \vC_{\vr, 0} \vq^e$, ${}_I \vr_{OP_{e+1}} = \vC_{\vr, 1} \vq^e$, $\vps_{e} = \vC_{\vps, 0} \vq^e$ and $\vps_{e+1} = \vC_{\vps, 1} \vq^e$. Note, the Boolean connectivity matrices are only used for the mathematical description of these extraction procedures. During a numerical implementation a different and much more efficient approach is used.

\subsection{Element-wise $\SE(3)$-interpolation}
For a given element interval $\mathcal{J}^e = [\xi^e, \xi^{e+1})$ the corresponding linear Lagrange basis functions $N_0^e$, $N_1^e$ and their first derivatives $N^e_{0,\xi}$, $N^e_{1,\xi}$ are defined by
\begin{equation}
	N^e_0(\xi) = \frac{\xi^{e+1} - \xi}{\xi^{e+1} - \xi^e} \, , \quad N^e_1(\xi) = \frac{\xi - \xi^e}{\xi^{e+1} - \xi^e} \, , 
	\quad N^e_{0,\xi}(\xi) = \frac{-1}{\xi^{e+1} - \xi^e} \, , \quad N^e_{1,\xi}(\xi) = \frac{1}{\xi^{e+1} - \xi^e} \, .
\end{equation}
Further, we introduce the characteristic function $\chi_{\mathcal{J}^e} \colon \mathcal{J} \to \{0, 1\}$, being one for $\xi \in \mathcal{J}^e = [\xi^e, \xi^{e+1})$ and zero elsewhere. By that, we can extend the global two-node $\SE(3)$-interpolation~\eqref{eq:SE(3)_interpolation} to a piecewise two-node $\SE(3)$-interpolation given by
\begin{equation}\label{eq:relative_interpolation_H_two_node}
	\begin{aligned}
		\vH_{\mathcal{I}\mathcal{K}}(\xi, \vq) &= \sum_{e=0}^{n_\mathrm{el} - 1} \chi_{\mathcal{J}^e}(\xi) \vH_{\mathcal{I}\mathcal{K}_e}(\vq) \Exp_{\SE(3)} \left(N^e_1(\xi) \,  \vth_{\mathcal{K}_e\mathcal{K}_{e+1}}(\vq)\right) \, , \\
		\vth_{\mathcal{K}_e\mathcal{K}_{e+1}}(\vq) &= \Log_{\SE(3)} \left[\vH^{-1}_{\mathcal{I}\mathcal{K}_e}(\vq) \, \vH_{\mathcal{I}\mathcal{K}_{e+1}}(\vq)\right] \, .
	\end{aligned}
\end{equation}
For symmetry reasons, cf. Crisfield and Jeleni{\'c}\cite{Crisfield1999}, it might be advantageous to use a modified $\SE(3)$-interpolation given by
\begin{equation}\label{eq:relative_interpolation_H_two_node_sym}
	\begin{aligned}
		\vH_{\mathcal{I}\mathcal{K}}(\xi, \vq) &= \sum_{e=0}^{n_\mathrm{el} - 1} \chi_{\mathcal{J}^e}(\xi) \vH_{\mathcal{I}\mathcal{R}_e}(\vq) \Exp_{\SE(3)} \left(N^e_0(\xi) \,  \vth_{\mathcal{R}_e\mathcal{K}_{e}}(\vq) + N^e_1(\xi) \,  \vth_{\mathcal{R}_e\mathcal{K}_{e+1}}(\vq)\right) \, , \\
		\vH_{\mathcal{I}\mathcal{R}_e}(\vq) &= \vH_{\mathcal{I}\mathcal{K}_e}(\vq) \Exp_{\SE(3)} \left(\tfrac{1}{2} \Log_{\SE(3)} \left[\vH^{-1}_{\mathcal{I}\mathcal{K}_e}(\vq) \, \vH_{\mathcal{I}\mathcal{K}_{e+1}}(\vq)\right]\right) \, , \\
		\vth_{\mathcal{R}_e\mathcal{K}_{e}}(\vq) &= \Log_{\SE(3)} \left[\vH^{-1}_{\mathcal{I}\mathcal{R}_e}(\vq) \, \vH_{\mathcal{I}\mathcal{K}_{e}}(\vq)\right] \, , \quad \vth_{\mathcal{R}_e\mathcal{K}_{e+1}}(\vq) = \Log_{\SE(3)} \left[\vH^{-1}_{\mathcal{I}\mathcal{R}_e}(\vq) \, \vH_{\mathcal{I}\mathcal{K}_{e+1}}(\vq)\right]  \, .
	\end{aligned}
\end{equation}
Since this strategy requires twice as many $\SE(3)$-exponential evaluations and three times as many $\SE(3)$-logarithm evaluations, the interpolation~\eqref{eq:relative_interpolation_H_two_node} is preferred due to superior efficiency and simplicity.

Using similar arguments that resulted in~\eqref{eq:global_strain_measures_SE(3)}, the interpolation~\eqref{eq:relative_interpolation_H_two_node} leads to piecewise constant strains

\begin{equation}\label{eq:discretized_picewise_strain_measures}
	\vvep(\xi, \vq) = \sum_{e=0}^{n_\mathrm{el} - 1} \chi_{\mathcal{J}^e}(\xi) \frac{\vth_{\mathcal{K}_e\mathcal{K}_{e+1}}(\vq)}{\xi^{e+1} - \xi^e} \frac{1}{J} \, .
\end{equation}
\correctb{Using the same reasoning following \eqref{eq:global_strain_measures_SE(3)} and since the piecewise two-node $\SE(3)$-interpolation can exactly represent constant strains within each element, neither membrane nor shear locking will appear with this discretization. As we need no further numerical strategies to avoid locking as for instance re-interpolation of strain measures\cite{Meier2015, Greco2017} or mixed formulations\cite{Santos2010, Santos2011, Betsch2016}, we call the finite element formulation intrinsically locking-free.}
\subsection{Petrov--Galerkin projection}
Following a Petrov--Galerkin projection\cite{Petrov1940}, we introduce the nodal virtual displacements $\delta \vs_i(t) = \big({}_I \delta \vr_{P_i}(t), {}_{K_i} \delta \vph_{IK_i}(t)\big) \in \mR^6$ given by the nodal centerline variation ${}_I \delta \vr_{P_i}(t) \in \mR^3$ and the nodal virtual rotation ${}_{K_i} \delta \vph_{IK_i}(t) \in \mR^3$. They are related to the variations of the nodal generalized coordinates by ${}_I \delta \vr_{P_i} = \delta({}_I \vr_{OP_i})$ and ${}_{K_i} \delta \vph_{IK_i} = \vT_{\SO(3)}(\vps_i) \delta \vps_i$, see Equation (49) of Cardona and Geradin~\cite{Cardona1988}. Led by the observation that $\overline{\vM}$ in~\eqref{eq:M_bar_and_g} is symmetric, it is advisable to introduce the nodal minimal velocities $\vu_i(t) = \big({}_I \vv_{P_i}(t), {}_{K_i} \vom_{IK_i}(t)\big) \in \mR^6$ given by the nodal centerline velocity ${}_I \vv_{P_i}(t) \in \mR^3$ and the nodal angular velocity ${}_{K_i} \vom_{IK_i}(t) \in \mR^3$.  Similar to the virtual displacements, the nodal minimal velocities are related to the time derivative of the nodal generalized coordinates by the nodal kinematic differential equation
\begin{equation}\label{eq:nodal_kinematic_equation}
	\dot{\vq}_i =
	\begin{pmatrix}
		{}_I\dot{\vr}_{OP_i} \\
		\dot{\vps}_i
	\end{pmatrix} =
	\begin{pmatrix}
		\mathbf{1}_{3 \times 3} & \mathbf{0}_{3 \times 3} \\
		\mathbf{0}_{3 \times 3} & \vT^{-1}_{\SO(3)}(\vps_i)
	\end{pmatrix}
	\begin{pmatrix}
		{}_I \vv_{P_i} \\
		{}_{K_i} \vom_{IK_i}
	\end{pmatrix} =
	\vB_i(\vq_i) \vu_i \, .
\end{equation}
Consequently, during a subsequent Galerkin projection, the symmetry of $\overline{\vM}$ is preserved and results in a symmetric mass matrix, see~\eqref{eq:discretized_inertial_virtual_work}.
Since the inverse tangent map $\vT^{-1}_{\SO(3)}$ used in~\eqref{eq:nodal_kinematic_equation} exhibits singularities for $\|\vps\| = k 2 \pi$ with $k = 0, 1, 2, \dots$, we
\correctb{apply the following strategy to avoid them. For $k = 0$, we use the first order approximation $\vT^{-1}_{\SO(3)}(\vom) = \mathbf{1}_{3 \times 3} + \frac{1}{2} \widetilde{\vom}$.} For $k > 0$, the concept of complement rotation vectors\cite{Cardona1988, Ibrahimbegovic1995} is applied. Due to the Petrov--Galerkin projection, it is sufficient to introduce a nodal update that is performed after each successful time step (in statics after each successful Newton increment). \correctb{This update, which corresponds to a change of coordinates for the orientation parametrization, is given by}
\begin{equation}
	\vps = 
	\begin{cases}
		\vps \, , & \|\vps\| \leq \pi \, , \\
		\vps^C = \big(1 - 2 \pi / \|\vps\|\big) \vps \, , & \|\vps\| > \pi \, ,
	\end{cases}
\end{equation}
with $\Exp_{\SO(3)}(\vps) = \Exp_{\SO(3)}(\vps^C)$, i.e., there is no difference whether the nodal transformation matrix $\vA_{IK}$ is described by the rotation vector $\vps$ or by its complement $\vps^C = \big(1 - 2 \pi / \|\vps\| \big) \vps$. Thus, for reasonable time steps (Newton increments) a both minimal and singularity free parametrization of $\SO(3)$ is obtained without changing the rod's virtual work formulation since all relevant nodal quantities are interpolated relatively\footnote{It should be mentioned that the proposed complement rotation vector formalism is restricted to numerical single step methods, e.g. Runge--Kutta family, generalized-$\alpha$ methods, etc., but cannot be used straightforwardly for multi step algorithms like BDF. Alternatively, another nodal rotation parametrization can be chosen, e.g., unit quaternions. Since in numerics, we cannot deal with unit quaternions, a general quaternion is constrained to be of unit length. While in statics these constraint equations have to be taken into account, during a numerical time integration they can be satisfied by applying an appropriate projection after each successful time step. Likewise, a modified kinematic differential equation can be introduced, see Section 6.9.6 of Egeland and Gravdahl\cite{Egeland2002}.}.

In analogy to the generalized coordinates, the nodal virtual displacements and minimal velocities are assembled in the global tuple of virtual displacements $\delta \vs(t) = \big(\delta\vs_0(t), \delta\vs_1(t), \dots, \delta\vs_{N-1}(t)\big) \in \mR^{n_{\vq}}$ and global tuple of minimal velocities $\vu(t) = \big(\vu_0(t), \vu_1(t), \dots, \vu_{N-1}(t)\big) \in \mR^{n_{\vq}}$. Using again the boolean connectivity matrix $\vC_e$, we can extract the element virtual displacements $\delta \vs^e(t) = \big({}_I \delta \vr_{P_e}(t), {}_{K_e} \!\delta\vph_{IK_e}(t), {}_I \delta \vr_{P_{e+1}}(t), {}_{K_{e+1}} \!\delta\vph_{IK_{e+1}}(t)\big) \in \mR^{12}$ and the element minimal velocities $\vu^e(t) = \big({}_I \vv_{P_e}(t), {}_{K_e} \! \vom_{IK_e}(t), {}_I \vv_{P_{e+1}}(t), {}_{K_{e+1}} \! \vom_{IK_{e+1}}(t)\big) \in \mR^{12}$ from the global quantities in agreement with $\delta \vs^e = \vC_e \delta \vs$ and $\vu^e = \vC_e \vu$. Further, the centerline variations ${}_I \delta \vr_{P_e}$, ${}_I \delta \vr_{P_{e+1}}$ and centerline velocities ${}_I \vv_{P_e}$, ${}_I \vv_{P_{e+1}}$ can be extracted from the tuple of virtual displacements $\delta \vs^e$ and minimal velocities $\vu^e$ using ${}_I \delta \vr_{P_e} = \vC_{\vr, 0} \delta \vs^e$, ${}_I \delta \vr_{P_{e+1}} = \vC_{\vr, 1} \delta \vs^e$  and ${}_I \vv_{P_e} = \vC_{\vr, 0} \vu^e$, ${}_I \vv_{P_{e+1}} = \vC_{\vr, 1} \vu^e$, respectively. Identical extraction operations can be defined for the nodal virtual rotations and angular velocities via $\vC_{\vps,0}, \vC_{\vps,1}$.

In the sense of a Petrov--Galerkin projection\cite{Petrov1940}, a different interpolation for the virtual displacements is chosen. We chose a simple interpolation strategy based on the previously introduced linear Lagrange basis function. As outlined in the previous paragraph, in order to obtain a symmetric and possibly constant mass matrix, the angular velocities are interpolated similarly. Formalizing that, we introduce the interpolations
\begin{equation}\label{eq:var_vel_interpolation}
	\begin{aligned}
		{}_I \delta \vr_{P}(\xi, \delta \vs) &= \sum_{e=0}^{n_\mathrm{el}-1} \chi_{\mathcal{J}^e}(\xi) \left(
		N^e_0(\xi) {}_I \delta \vr_{P_e} + N^e_1(\xi) {}_I \delta \vr_{P_{e+1}}
		\right) \, , \\
		{}_I \vv_{P}(\xi, \vu) &= \sum_{e=0}^{n_\mathrm{el}-1} \chi_{\mathcal{J}^e}(\xi) \left(
		N^e_0(\xi) {}_I \vv_{P_e} + N^e_1(\xi) {}_I \vv_{P_{e+1}}
		\right) \, , \\
		{}_K \delta \vph_{IK}(\xi, \delta \vs) &= \sum_{e=0}^{n_\mathrm{el}-1} \chi_{\mathcal{J}^e}(\xi) \left(
		N^e_0(\xi) \, {}_{K_e} \! \delta \vph_{IK_e} + N^e_1(\xi) \, {}_{K_{e+1}} \! \delta \vph_{IK_{e+1}}
		\right) \, , \\
		{}_K \vom_{IK}(\xi, \vu) &= \sum_{e=0}^{n_\mathrm{el}-1} \chi_{\mathcal{J}^e}(\xi) \left(
		N^e_0(\xi) \, {}_{K_e} \! \vom_{IK_e} + N^e_1(\xi) \, {}_{K_{e+1}} \! \vom_{IK_{e+1}}
		\right) \, .
	\end{aligned}
\end{equation}
Choices like ${}_K \delta \vr_{P}$, i.e., centerline variations in the cross-section $K$-basis, introduce cumbersome changes of bases when for instance the generalized contact force directions should be computed, cf. Bosten et al.\cite{Bosten2022}. This is quite inconvenient especially when dealing with complex flexible multibody system connected by perfect bilateral constraints as treated in G{\'e}radin and Cardona~\cite{Geradin2001}.
\subsection{Discrete virtual work formulations}
Now we have all ingredients for discretizing the previously introduced continuous virtual work formulations in space. Inserting~\eqref{eq:var_vel_interpolation} into the internal virtual work~\eqref{eq:internal_virtual_work2}, their discretization is obtained
\begin{equation}\label{eq:discretized_internal_virtual_work}
	\begin{aligned}
		\delta W^\mathrm{int}(\delta \vs;\vq) &= \delta \vs\T \vf^{\mathrm{int}}(\vq) \, , \quad \vf^{\mathrm{int}}(\vq) = \sum_{e=0}^{n_\mathrm{el} - 1} \vC_{e}\T \vf^{\mathrm{int}}_e(\vC_{e}\vq) \, , \\
		\vf^{\mathrm{int}}_e(\vq^e) &= -\int_{\mathcal{J}^e} \sum_{i=0}^{1}\Big\{ N^e_{i,\xi} \vC_{\vr, i}\T \vA_{IK} \, {}_K \vn + N^e_{i,\xi} \vC_{\vps, i}\T {}_K \vm 
		-N^e_{i} \vC_{\vps, i}\T \left({}_K \bar{\vga}\times {}_K \vn + {}_K\bar{\vka}_{IK} \times {}_K \vm \right) \Big\} \diff[\xi] \, ,
	\end{aligned}
\end{equation}
where we have introduced the internal forces $\vf^{\mathrm{int}}$ together with their element contributions $\vf^{\mathrm{int}}_e$. The discrete ansatzfunction of the interpolation of~\eqref{eq:relative_interpolation_H_two_node} is used in the computation of the appearing contributions. For the sake of readability, above and subsequently, we partly suppress the function arguments, which should be clear from the context. Similarly, the discretization of the external virtual work~\eqref{eq:external_virtual_work} takes the form
\begin{equation}
	\begin{aligned}
		\delta W^\mathrm{ext}(\delta \vs; \vq) &= \delta \vs\T \vf^{\mathrm{ext}}(\vq) \, , \\  \vf^{\mathrm{ext}}(\vq) & = \sum_{e=0}^{n_\mathrm{el} - 1} \vC_e\T \vf^{\mathrm{ext}}_e(\vC_{e}\vq) 
		+ \vC_{n_\mathrm{el} - 1}\T \left[\vC_{\vr, 1}\T {}_I \vb_1 + \vC_{\vph, 1}\T {}_K \vc_1 \right]_{\xi=1}
		+ \vC_0\T \left[\vC_{\vr, 0}\T {}_I \vb_0 + \vC_{\vph, 0}\T {}_K \vc_0 \right]_{\xi=0} \, , \\
		\vf^{\mathrm{ext}}_e(\vq^e) &= \int_{\mathcal{J}^e} \sum_{i=0}^{1}\Big\{ N^e_{i}  \vC_{\vr, i}\T {}_I \vb + N^e_{i} \vC_{\vps, i}\T {}_K \vc \Big\} J \diff[\xi] \, ,
	\end{aligned}
\end{equation}
where we have introduced the external forces $\vf^{\mathrm{ext}}$ and their element contributions $\vf^{\mathrm{ext}}_e$. Finally, the discretization of the inertial virtual work contributions~\eqref{eq:inertia_virtual_work} is given by
\begin{equation}\label{eq:discretized_inertial_virtual_work}
	\begin{aligned}
		\delta W^\mathrm{dyn}(\delta \vs;\vq,\vu) &= -\delta \vs\T \left\{\vM(\vq) \dot{\vu} + \vf^{\mathrm{gyr}}(\vq, \vu) \right\} \, , \\
		\vM(\vq) &= \sum_{e=0}^{n_\mathrm{el} - 1} \vC_e\T \vM_e(\vC_{e}\vq) \vC_e \, , \quad
		\vf^{\mathrm{gyr}}(\vq, \vu) = \sum_{e=0}^{n_\mathrm{el}-1} \vC_e\T \vf_e^{\mathrm{gyr}}(\vC_{e} \vq, \vC_{e} \vu) \, , \\
\end{aligned}
\end{equation}
where
\begin{equation}
	\begin{aligned}
		\vM_e(\vq^e) &= 
		\begin{multlined}[t]
			\int_{\mathcal{J}^e} \sum_{i=0}^{1} \sum_{k=0}^{1} N^e_{i} N^e_{k} \Big\{
			\vC_{\vr, i}\T A_{\rho_0} \mathbf{1}_{3\times3} \vC_{\vr, k} 
			+ \vC_{\vr, i}\T \vA_{IK} {}_K \vS_{\rho_0}\T \vC_{\vps, k} \\
			~~~~~~~~~~~~~~~~~~~~~~~~~~~~+ \vC_{\vps, i}\T {}_K \vS_{\rho_0} \vA_{IK}\T \vC_{\vr, k} + \vC_{\vps, i}\T {}_K \vI_{\rho_0} \vC_{\vps, k}
			\Big\} J \diff[\xi] \, , 
		\end{multlined} \\
		\vf^{\mathrm{gyr}}_e(\vq^e, \vu^e) &= \int_{\mathcal{J}^e} \sum_{i=0}^{1} N^e_i \Big\{ \vC_{\vr, i}\T \vA_{IK} {}_K \widetilde{\vom}_{IK} {}_K \vS_{\rho_0}\T {}_K \vom_{IK}
		+ \vC_{\vps, i}\T {}_K \widetilde{\vom}_{IK} {}_K \vI_{\rho_0} {}_K \vom_{IK} \Big\} J \diff[\xi]  \, ,
	\end{aligned}
\end{equation}
where we have introduced the symmetric mass matrix $\vM$ and the gyroscopic forces $\vf^\mathrm{gyr}$ together with their element contributions $\vM_e$ and $\vf^\mathrm{gyr}_e$.

\correctb{Appendix~\ref{sec:discrete_preservation_properties} discusses the discrete conservation properties of the proposed rod finite element formulation. That is, the conservation of total energy, linear and angular momentum. Moreover, another $\SE(3)$ finite element formulation is presented, for which the nodal virtual rotations and angular velocities expressed in their cross-section-fixed bases are replaced by the nodal virtual rotations ${}_{I} \delta \vph_{IK_i}(t) \in \mR^3$ and nodal angular velocities ${}_{I} \vom_{IK_i}(t) \in \mR^3$ expressed in the inertial $I$-basis.}

The linearization of the internal forces $\vf^\mathrm{int}$ with respect to the generalized coordinates $\vq$ is given by
\begin{equation}
	\begin{aligned}
		\pd{\vf^{\mathrm{int}}}{\vq}(\vq) &= \sum_{e=0}^{n_\mathrm{el} - 1} \vC_{e}\T \pd{\vf^{\mathrm{int}}}{\vq^e}(\vC^e \vq) \, \vC^e \, , \\
		\pd{\vf^{\mathrm{int}}}{\vq^e}(\vq^e) &= -\int_{\mathcal{J}^e} \sum_{i=0}^{1}
		\begin{multlined}[t]
			\bigg\{
			N^e_{i,\xi} \vC_{\vr, i}\T \pd{\left(\vA_{IK} \, {}_K \vn\right)}{\vq^e} 
			+ N^e_{i,\xi} \vC_{\vps, i}\T \pd{{}_K \vm}{\vq^e} \\
			-N^e_{i} \vC_{\vps, i}\T \left({}_K \widetilde{\bar{\vga}} \pd{{}_K \vn}{\vq^e} - {}_K \widetilde{\vn} \pd{{}_K \bar{\vga}}{\vq^e} + {}_K \widetilde{\bar{\vka}}_{IK} \pd{{}_K \vm}{\vq^e} - {}_K \widetilde{\vm} \pd{{}_K \bar{\vka}_{IK}}{\vq^e} \right)
			\bigg\} \diff[\xi] \, ,
		\end{multlined} 
	\end{aligned}
\end{equation}
which requires the derivative of the rotation matrix $\vA_{IK}$ and the derivatives of the internal forces and moments
\begin{equation}
	\begin{aligned}
		\pd{{}_K \vn}{\vq^e} = \pd{{}_K \vn}{{}_K \vga} \pd{{}_K \vga}{\vq^e} + \pd{{}_K \vn}{{}_K \vka_{IK}} \pd{{}_K \vka_{IK}}{\vq^e} \, , \quad \pd{{}_K \vm}{\vq^e} = \pd{{}_K \vm}{{}_K \vga} \pd{{}_K \vga}{\vq^e} + \pd{{}_K \vm}{{}_K \vka_{IK}} \pd{{}_K \vka_{IK}}{\vq^e} \, .
	\end{aligned}
\end{equation}
Both are depending on the derivatives of the strain measures ${}_K \vga$ and ${}_K \vka_{IK}$ introduced in~\eqref{eq:discretized_picewise_strain_measures}. These derivatives together with the derivative of the interpolation formula~\eqref{eq:relative_interpolation_H_two_node}, i.e., the position vector ${}_I \vr_{OP}$ and the transformation matrix $\vA_{IK}$ are given in Appendix~\ref{sec:linearization}. As already noted, the quantity $\vS_{\rho_0}$ vanishes if the centerline points ${}_I \vr_{OP}$ coincide with the cross-sections' center of mass. For this case the gyroscopic forces $\vf^\mathrm{gyr}$ and the term $\vM \dot{\vu}$ get independent of the generalized coordinates $\vq$ which is why we do not consider their linearizations here. The linearization of the external forces is omitted as well since a possible $\vq$ dependence is defined by the specific form of the external forces and moments only.

Arising element integrals of the form $\int_{\mathcal{J}^e} f(\xi) \diff[\xi]$ encountered in the discretized internal, external and gyroscopic forces, as well as in the mass matrix, are subsequently computed using a two-point Gauss--Legendre quadrature rule.
\subsection{Equations of motion and static equilibrium}
The principle of virtual work states that the totality of virtual work contributions has to vanish for arbitrary virtual displacements at all time instants $t$, see Chapter 8 of dell'Isola and Steigmann\cite{dellIsola2020b}, i.e.,
\begin{equation}
	\delta W^\mathrm{tot} = \delta W^\mathrm{int} + \delta W^\mathrm{ext} + \delta W^\mathrm{dyn} \stackrel{!}{=} 0 \quad \forall \delta \vs(t), \forall t \, .
\end{equation}
Thus, the equations of motion
\begin{equation}
	\dot{\vu} = \vM^{-1}(\vq) \left(\vf^\mathrm{int}(\vq) + \vf^\mathrm{ext}(\vq) - \vf^\mathrm{gyr}(\vq, \vu)\right)
\end{equation}
have to be fulfilled for each instant of time $t$. Further, the nodal minimal velocities $\vu_i$ are related to the time derivatives of the nodal generalized coordinates $\dot{\vq}_i$ via the nodal kinematic differential equation~\eqref{eq:nodal_kinematic_equation}, which can be assembled to a global kinematic differential equation of the form $\dot{\vq} = \vB(\vq) \vu$. Depending on the specific application, either the system of ordinary differential equations
\begin{equation}
	\begin{aligned}
		\dot{\vq} &= \vB(\vq) \vu \, , \\
		\dot{\vu} &= \vM(\vq)^{-1} \left(\vf^\mathrm{int}(\vq) + \vf^\mathrm{ext}(\vq) - \vf^\mathrm{gyr}(\vq, \vu)\right) \, ,
	\end{aligned}
\end{equation}
or the nonlinear generalized force equilibrium
\begin{equation}\label{eq:nonlinear_generalized_force_equilibrium}
	\vf^\mathrm{int}(\vq) + \vf^\mathrm{ext}(\vq) = \mathbf{0}
\end{equation}
\correctb{is obtained. The system of ordinary differential equations can be solved using standard higher-order ODE solvers (e.g. family of explicit\cite{Hairer1993} and implicit\cite{Jay1995, Hairer2002} Runge--Kutta methods or structure-preserving algorithms\cite{Jay1996, Hairer2006}). This is a remarkable result since existing $\SE(3)$ rod formulations require highly customized Lie group solvers\cite{Bruels2012, Sonneville2014}. In order to apply well-established methods from structural dynamics like the Newmark-$\beta$\cite{Newmark1959} or the generalized-$\alpha$ method\cite{Chung1993}, a slightly modified update of the generalized coordinates has to be applied, see Equation (37a) and (37b) of Arnold et al.\cite{Arnold2016}. Alternatively, a generalized-$\alpha$ formulation for first-order differential equations\cite{Jansen2000} can be used without any modifications. The nonlinear generalized force equilibrium~\eqref{eq:nonlinear_generalized_force_equilibrium} is solved by any root finding algorithm (e.g. Newton--Raphson, Riks). Note that a system of linear equations with a non-symmetric matrix must be solved in each iteration.}
\section{Numerical experiments}\label{sec:numerical_examples}
\correctb{This section demonstrates the power of the developed $\SE(3)$ finite element formulation by a variety of selected benchmark examples. The \emph{cantilever experiment} slightly extends the investigations of Meier et al.\cite{Meier2015} and gives numerical evidence for the absence of locking and further demonstrates a second-order spatial convergence of the proposed formulation. The \emph{helix experiment}, which was used in Harsch et al. \cite{Harsch2021a} to study locking for a director beam formulation, is used as a second example to show that the $\SE(3)$-formulation is intrinsically locking-free. Objectivity after discretization is demonstrated in the example \emph{superimposed rotation of deformed cantilever}. Large inhomogeneous deformations are studied in the example \emph{rod bent to a helical form}\cite{Ibrahimbegovic1997}. Finally, dynamics is examined by the \emph{flexible heavy top} example.}

We used a Newton--Raphson method to solve all static experiments with an absolute tolerance $atol$ in terms of the max error of the total generalized forces. Prescribed boundary conditions were incorporated into the principle of virtual work using perfect bilateral constraints\cite{Geradin2001}. The dynamic problem was solved using a standard fourth order Runge--Kutta method \correctb{as well as the generalized-$\alpha$ method for first order differential equations\cite{Jansen2000}}. For simplicity, boundary conditions were enforced on acceleration level only. \correctb{Unless stated otherwise, arising element integrals were evaluated using a two-point Gauss--Legendre quadrature rule.}
\subsection{Cantilever experiment}\label{sec:quarter_circle}
\begin{figure}
	\centering
	\begin{subfigure}[b]{0.23\textwidth}
		\caption{}
		\includegraphics[width=\textwidth, trim=20cm 8cm 20cm 8cm, clip]{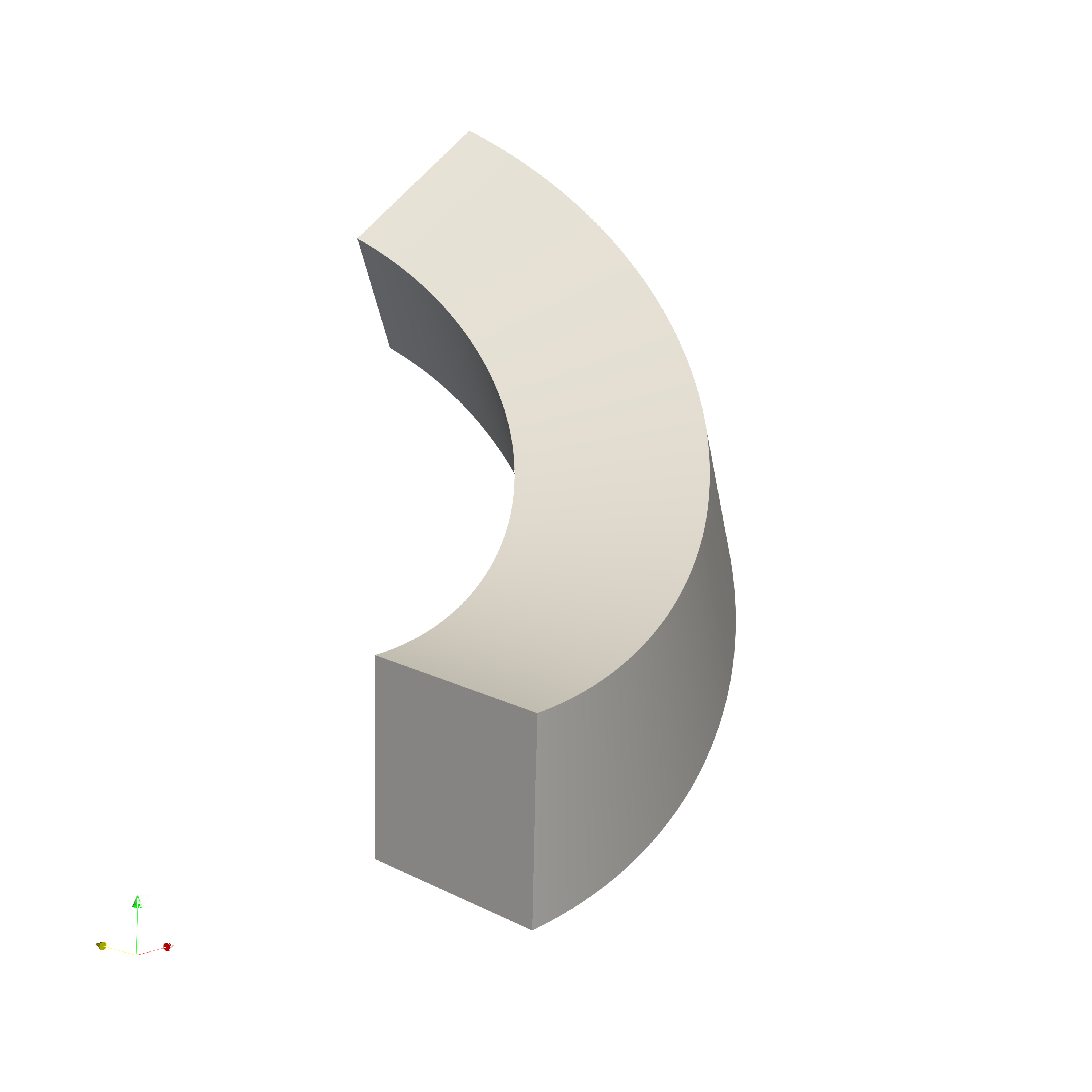}
	\end{subfigure}
	~
	\begin{subfigure}[b]{0.23\textwidth}
		\caption{}
		\includegraphics[width=\textwidth, trim=20cm 8cm 20cm 8cm, clip]{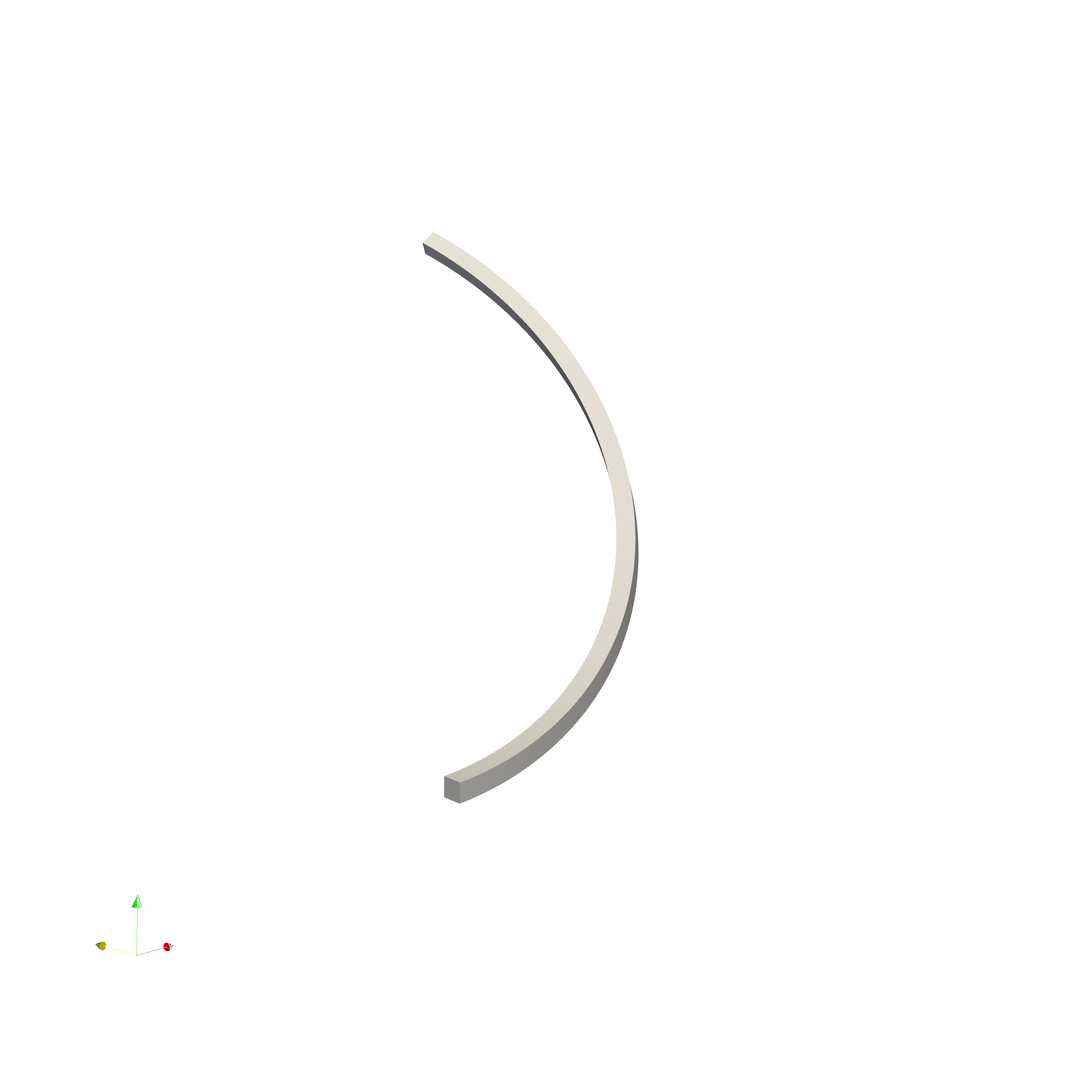}
	\end{subfigure}
	~
	\begin{subfigure}[b]{0.23\textwidth}
		\caption{}
		\includegraphics[width=\textwidth, trim=20cm 8cm 20cm 8cm, clip]{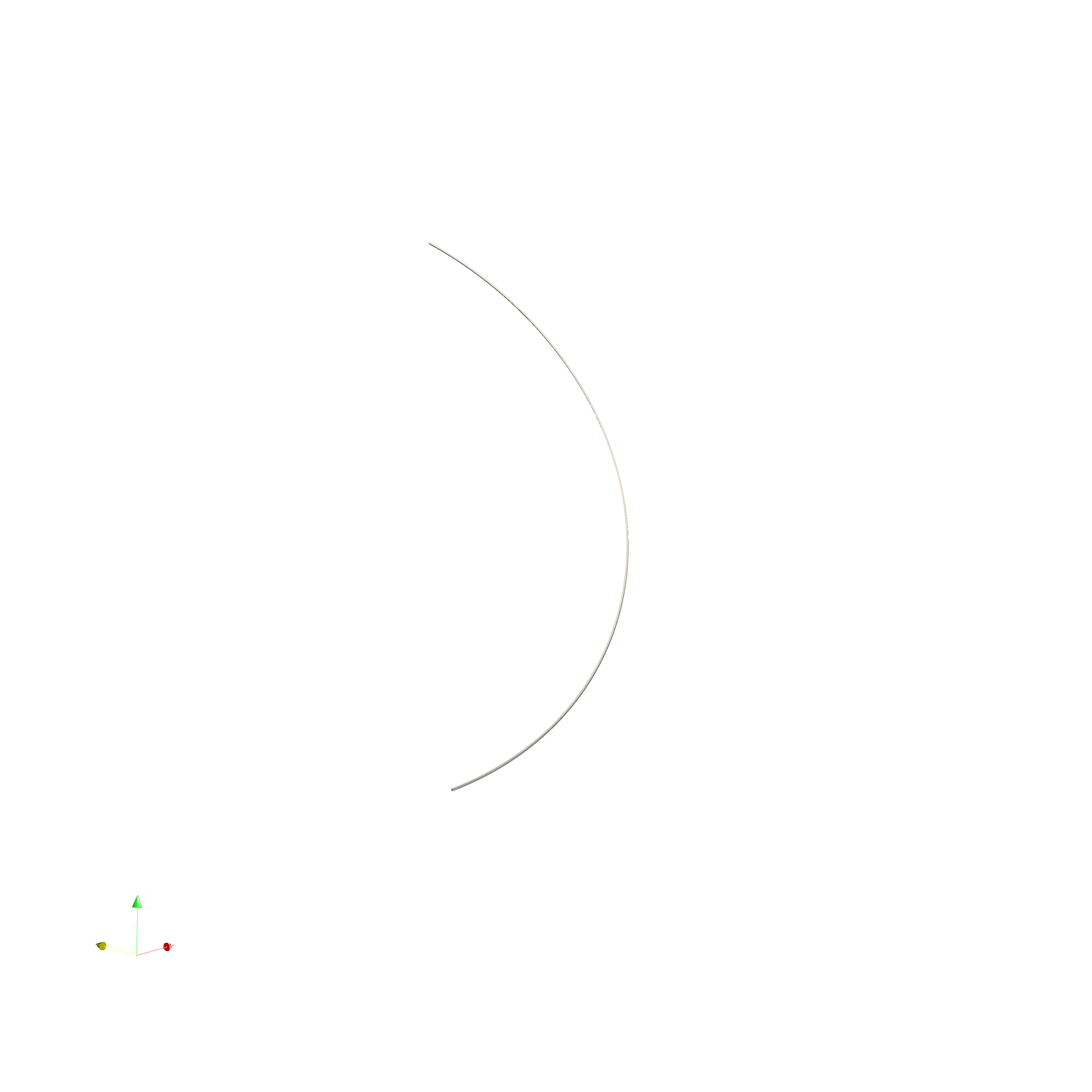}
	\end{subfigure}
	~
	\begin{subfigure}[b]{0.23\textwidth}
		\caption{}
		\includegraphics[width=\textwidth, trim=20cm 8cm 20cm 8cm, clip]{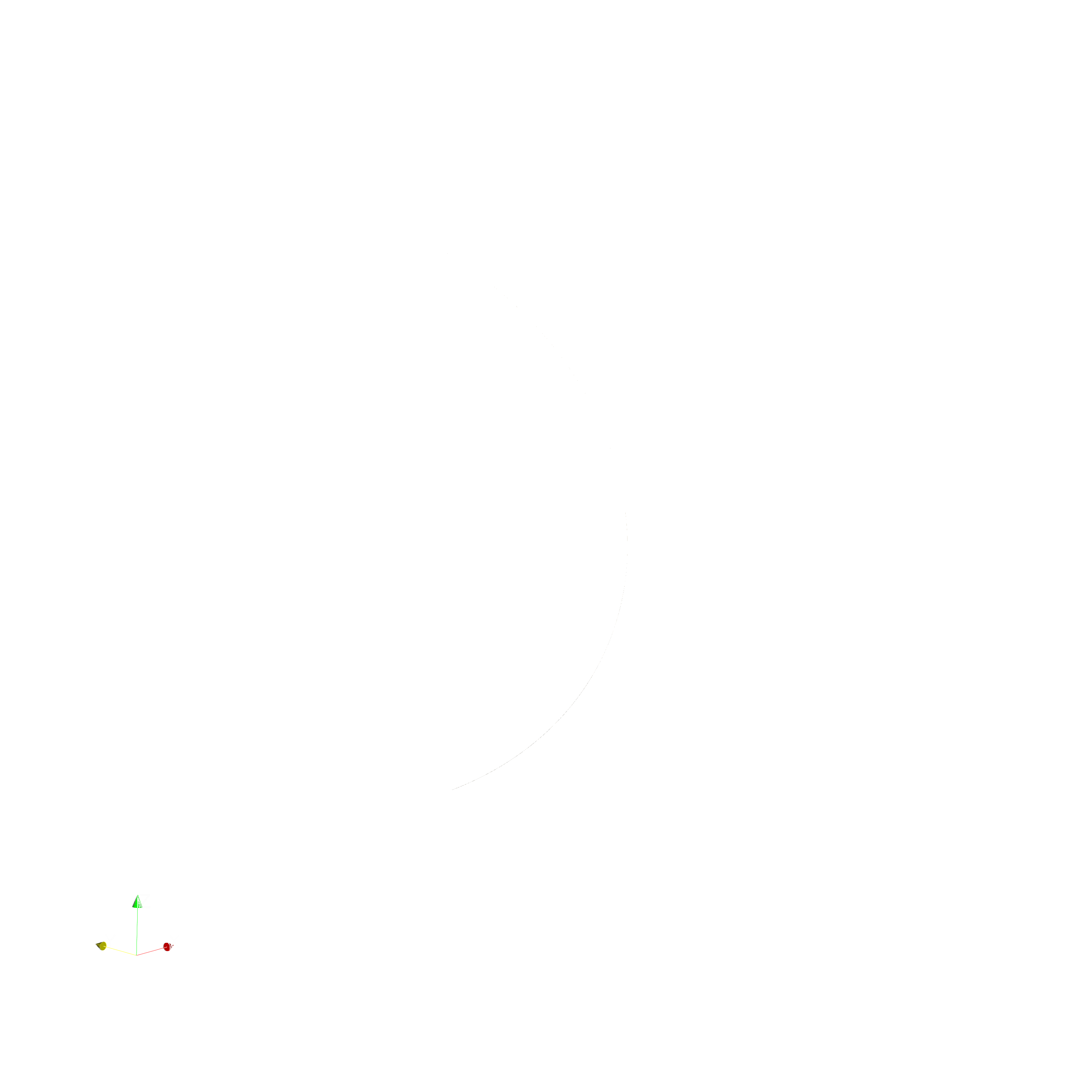}
	\end{subfigure}
	\caption{\correctb{Visualization of deformed configuration of the \emph{cantilever experiment} for (a) $\rho = \SI{e1}{}$, (b) $\rho = \SI{e2}{}$, (c) $\rho = \SI{e3}{}$ and (d) $\rho = \SI{e4}{}$ for five elements of the proposed $\SE(3)$-rod element.}}
	\label{fig:quarter_circle_final_configurations}
\end{figure}
\begin{figure}
	\centering
	\begin{subfigure}[b]{0.32\textwidth}
		\caption{}
		\vspace{-0.4cm}
		\includegraphics{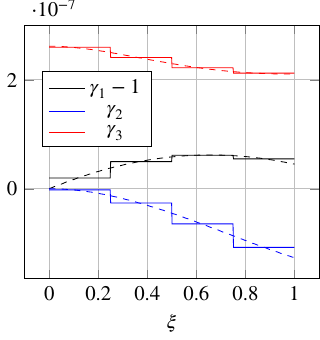}
	\end{subfigure}
	%
	\hspace{-0.25cm}
	\begin{subfigure}[b]{0.32\textwidth}
		\caption{}
		\vspace{-0.4cm}
		\includegraphics{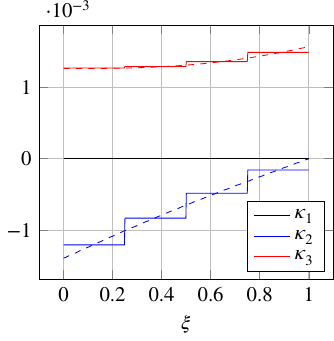}
	\end{subfigure}
	%
	%
	\begin{subfigure}[b]{0.33\textwidth}
		\caption{}
		\includegraphics{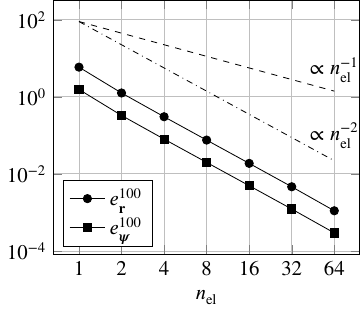}
	\end{subfigure}
	\caption{\correctb{(a), (b) Strain measures of the solution found by four $\SE(3)$-elements (solid) and those of the reference solution computed with $n_\mathrm{el} = 256$ elements (dashed). (c) Convergence behavior of the $\SE(3)$-formulation.}}
	\label{fig:convergence}
\end{figure}

We consider an initially straight cantilever rod of length $L = \SI{e3}{}$ with a quadratic cross-section of width $w$ subjected to a tip moment ${}_K\vc_1 = (0, \ 0, \ 0.5 \pi k_\mathrm{b} / L)$ and an out-of-plane tip load ${}_I\vb_1 = \vA_{IK} \, (0, \ 0, \ 0.5 \pi k_\mathrm{b} / L^2)$. In order to show the absence of locking, different slenderness ratios $\rho = L / w \in \{\SI{e1}{}, \SI{e2}{}, \SI{e3}{}, \SI{e4}{}\}$ are considered, i.e., widths $w \in \{\SI{e2}{}, \ \SI{e1}{}, \ \SI{e0}{}, \ \SI{e-1}{}\}$. ranging from $w=\SI{e-1}{}$ to $w=\SI{e2}{}$. Further, the elastic constants are given in terms of the Young's and shear moduli $E=1$ and $G=0.5$. That is, axial stiffness $k_\mathrm{e} = EA$, shear stiffness $k_\mathrm{s} = GA$, bending stiffnesses $k_{\mathrm{b}} = k_{\mathrm{b}_y} = k_{\mathrm{b}_z}= EI$ and torsional stiffness $k_\mathrm{t} = 2GI$, together with $A=w^2$ and $I =w^4/12$. As there is no analytical solution for this load case, the numerical solutions obtained for a single $\SE(3)$-element and for $\SI{64}{}$ $\SE(3)$-elements are compared to a reference solution found by a finite element implementation similar to Eugster et al.\cite{Eugster2014c} discretized by $\SI{256}{}$ linear elements. \correctb{To reduce shear locking in the latter formulation, reduced integration was applied. The strain measures of the reference solution are plotted in Figure~\ref{fig:convergence} (a) and (b).} For the centerline, we define the error
\begin{equation}
	e^k_{\vr} = \frac{1}{k} 
	\sqrt{
		\sum_{i=0}^{k-1} {}_I \Delta \vr_{P}\T(\xi_i) \, {}_I \Delta  \vr_{P}(\xi_i)
	} 
	\, , \quad {}_I \Delta \vr_{P}(\xi_i) = {}_I \vr_{OP}(\xi_i) - {}_I \vr_{OP}^{*}(\xi_i) \, , \quad \xi_i = \frac{i}{k - 1} \, ,
\end{equation}
where ${}_I \vr_{OP}^{*}$ denotes the centerline points of the reference solution and ${}_I \vr_{OP}$ the same quantity of the rod in comparison. In order to investigate the error in orientations, we use the measure\cite{Barfoot2014}
\begin{equation}
	e^k_{\vps} = \frac{1}{k} \sqrt{\sum_{i=0}^{k - 1} \Delta \vps\T(\xi_i) \Delta \vps(\xi_i) } \, , \quad \Delta \vps(\xi_i) = \Log_{\SO(3)}\big(\vA_{IK}\T(\xi_i) \vA_{IK}^{*}(\xi_i)\big) \, , \quad \xi_i = \frac{i}{k - 1} \, .
\end{equation}
Again, $\vA_{IK}^{*}$ denotes the orthogonal transformation matrix of the reference solution and $\vA_{IK}$ the same quantity of the rod in comparison. Since the relative rotation vector $\Delta \vps$ should be tending to zero during spatial convergence, this is a well-suited error measure \cite{Huynh2009}. The final loads were applied for all slenderness ratios within $\SI{20}{}$ increments and absolute tolerances as documented in Table~\ref{tab:locking}. \correctb{It can be observed that the spatial convergence behavior of the rod is unaffected by the chosen slenderness ratio, although no reduced integration was applied, which numerically proofs the absence of locking. For the different slenderness ratios, the deformed configurations are visualized in Figure~\ref{fig:quarter_circle_final_configurations}.}
\begin{table}
	\centering
	\caption{Experimental parameters and errors of the \emph{cantilever experiment}.}\label{tab:locking}
	\begin{tabular}{cccccccc}
		\toprule
		\multirow{2}{*}{slenderness} & \multirow{2}{*}{tolerance} & \multicolumn{2}{c}{$n_\mathrm{el} = 1$} & \multicolumn{2}{c}{$n_\mathrm{el} = 64$} \\
		& & $e^{100}_{\vr}$ & $e^{100}_{\vps}$ & $e^{100}_{\vr}$ & $e^{100}_{\vps}$ \\
		\midrule
		$\SI{e1}{}$ & $\SI{e-8}{}$ & $\SI{6.789e0}{}$ & $\SI{1.628e0}{}$ & $\SI{1.396e-03}{}$ & $\SI{3.355e-04}{}$ \\
		$\SI{e2}{}$ & $\SI{e-9}{}$ & $\SI{6.792e0}{}$ & $\SI{1.629e0}{}$ & $\SI{1.397e-03}{}$ & $\SI{3.357e-04}{}$ \\
		$\SI{e3}{}$ & $\SI{e-10}{}$ & $\SI{6.792e0}{}$ & $\SI{1.629e0}{}$ & $\SI{1.397e-03}{}$ & $\SI{3.357e-04}{}$ \\
		$\SI{e4}{}$ & $\SI{e-11}{}$ & $\SI{6.792e0}{}$ & $\SI{1.629e0}{}$ & $\SI{1.397e-03}{}$ & $\SI{3.356e-04}{}$ \\
		\bottomrule
	\end{tabular}
\end{table}

Due to the applied out-of-plane force, the resulting deformation is inhomogenous and cannot be described adequately by a small number of $\SE(3)$-elements with constant strain measures, see Figure~\ref{fig:convergence} (a) and (b). \correctb{This makes the example suitable for the investigation of the spatial convergence of the proposed formulation. For that, we used a moderate slenderness ratio $\rho = \SI{e3}{}$, an absolute tolerance of $atol = \SI{e-10}{}$ and $\SI{20}{}$ static increments. The numerical errors were computed with respect to the same reference solution as above.} Figure~\ref{fig:convergence} (c) visualizes the convergence behavior, which shows a second-order rate of spatial convergence in both centerline and rotation fields. This is in line with the observations made by Sonneville et al.\cite{Sonneville2014} within a different experiment.
\correctb{
	\subsection{Helix experiment}\label{sec:helix}
	\begin{figure}
		\centering
		\begin{subfigure}[b]{0.23\textwidth}
			\caption{}
			\includegraphics[width=\textwidth, trim=20cm 8cm 20cm 8cm, clip]{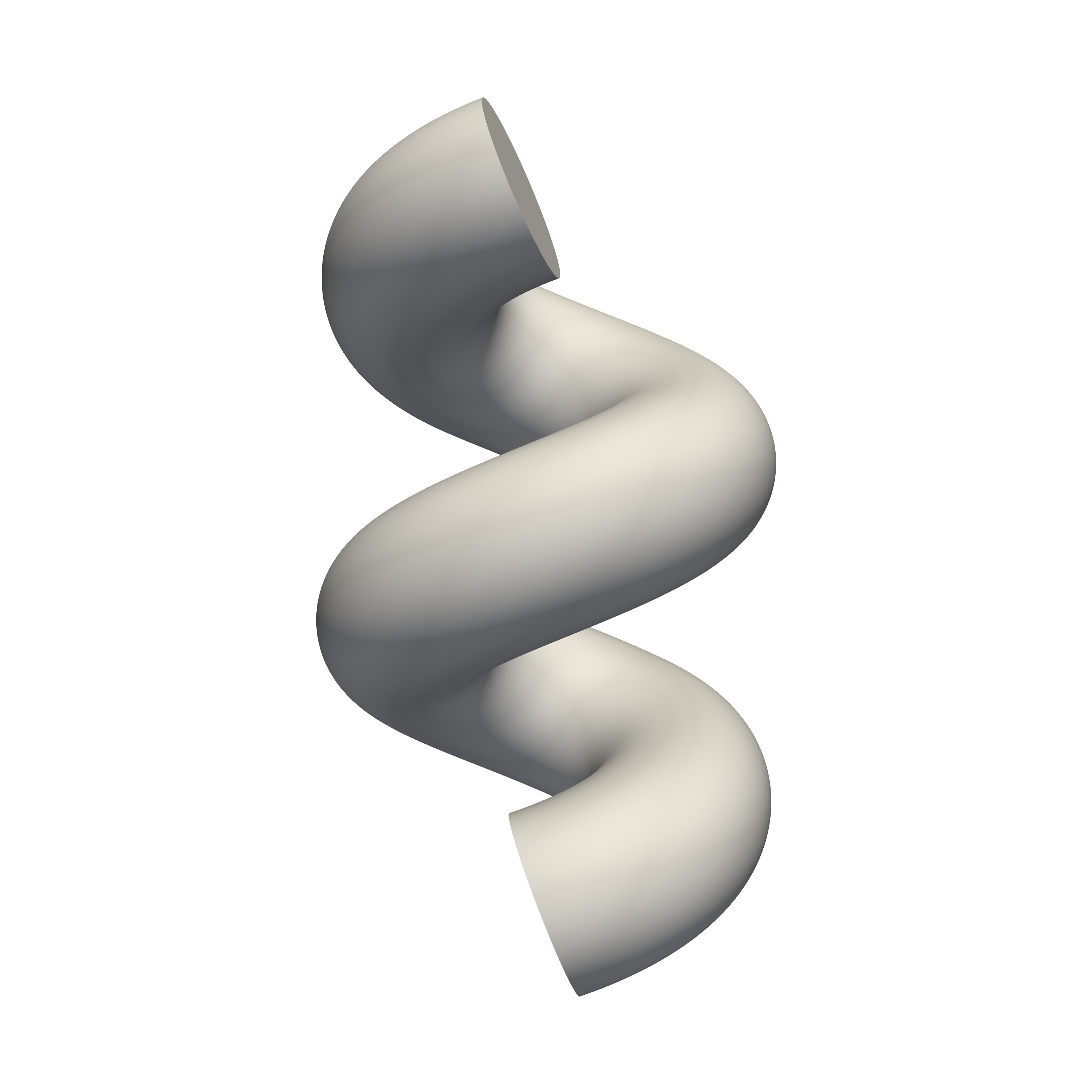}
		\end{subfigure}
		~
		\begin{subfigure}[b]{0.23\textwidth}
			\caption{}
			\includegraphics[width=\textwidth, trim=20cm 8cm 20cm 8cm, clip]{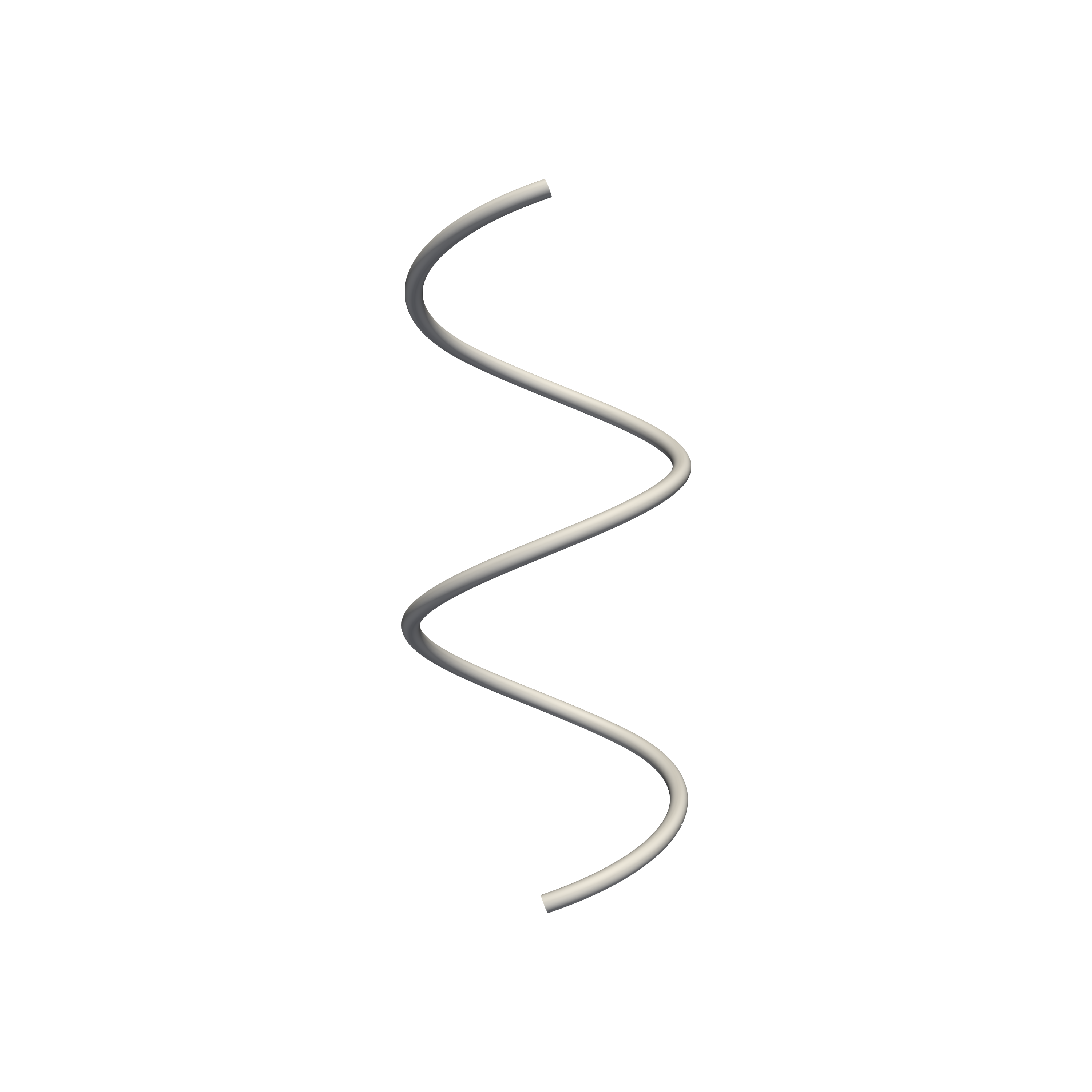}
		\end{subfigure}
		~
		\begin{subfigure}[b]{0.23\textwidth}
			\caption{}
			\includegraphics[width=\textwidth, trim=20cm 8cm 20cm 8cm, clip]{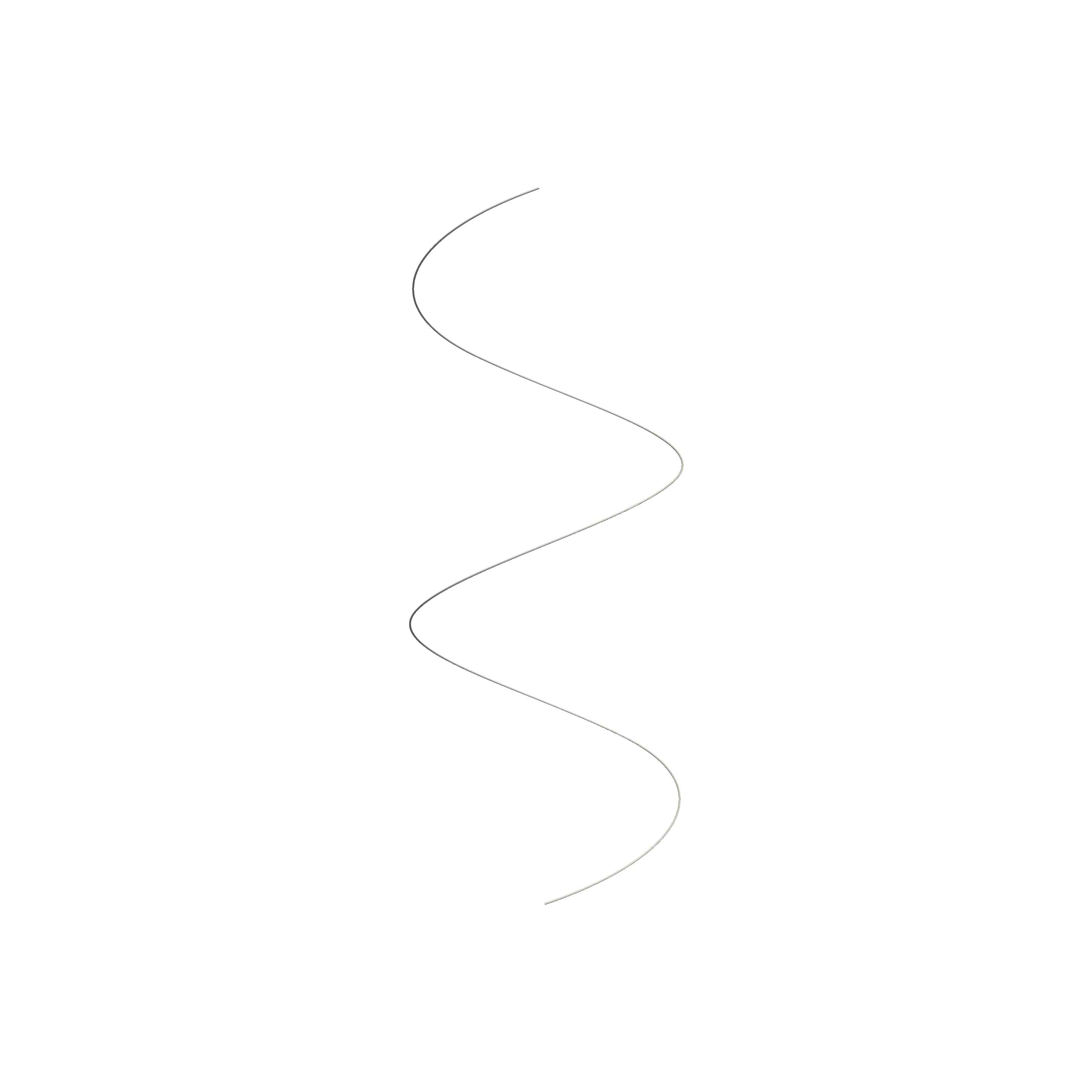}
		\end{subfigure}
		~
		\begin{subfigure}[b]{0.23\textwidth}
			\caption{}
			\includegraphics[width=\textwidth, trim=20cm 8cm 20cm 8cm, clip]{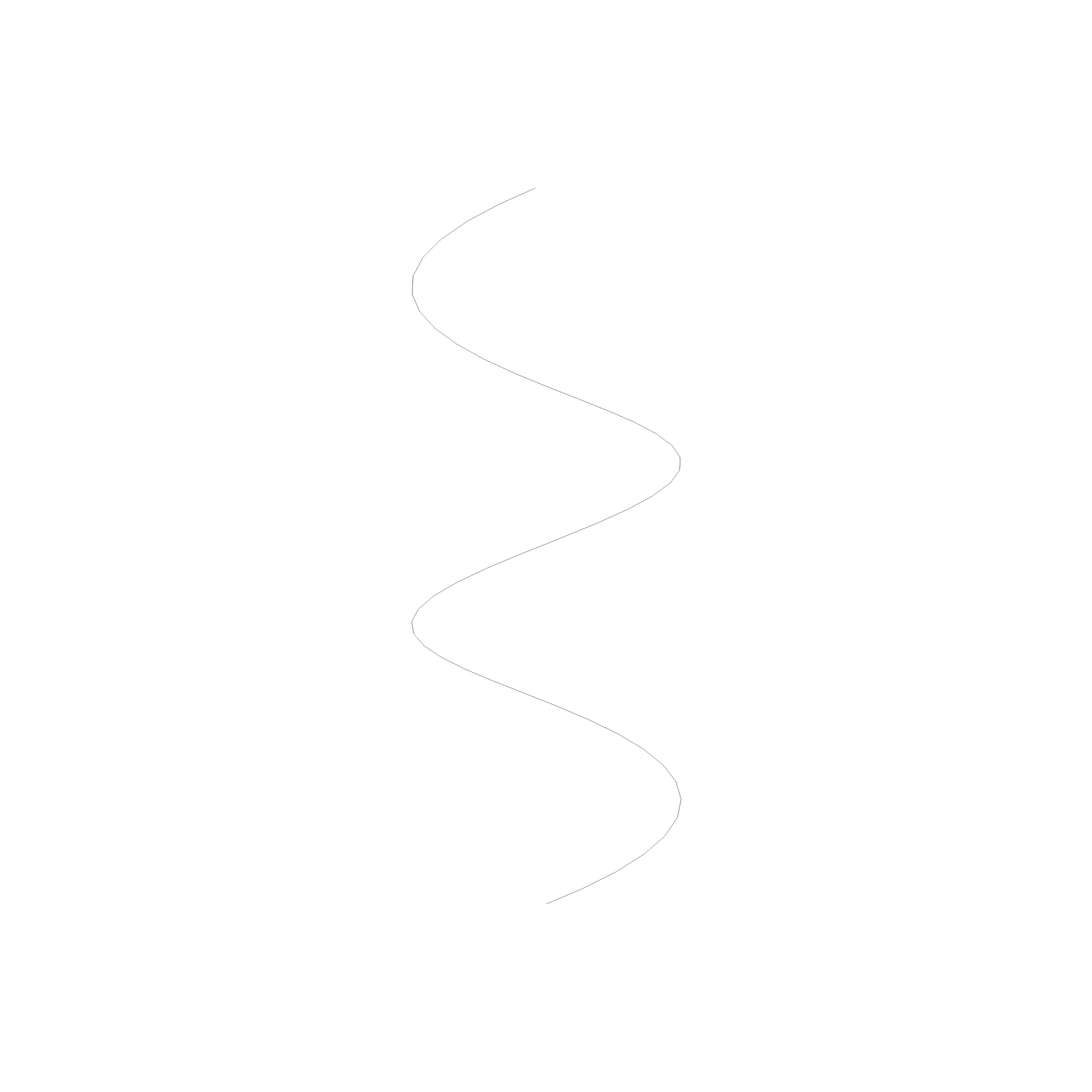}
		\end{subfigure}
		\caption{\correctb{Visualization of deformed configuration of the \emph{helix experiment} for (a) $\rho = \SI{e1}{}$, (b) $\rho = \SI{e2}{}$, (c) $\rho = \SI{e3}{}$ and (d) $\rho = \SI{e4}{}$ for five elements of the proposed $\SE(3)$-rod element.}}
		\label{fig:helix_final_configurations}
	\end{figure}
	In this example, an initially straight rod, is loaded such that it is deformed to a perfect helix\cite{Harsch2021a} with $n=2$ coils of radius $R_0=10$ and height $h=50$, see Figure~\ref{fig:helix_final_configurations}. Introducing the abbreviation $c = h / (R_0 2 \pi n) \geq 1$, the centerline of such a deformed rod is described by
	\begin{equation}\label{eq:analytical_solution_helix}
		{}_I \vr_{OP}^{*}(\xi) = 
		R_0 \begin{pmatrix}
			\sin\alpha(\xi) \\
			- \cos\alpha(\xi) \\
			c \alpha(\xi)
		\end{pmatrix} \, , \quad \alpha(\xi) = 2 \pi n \xi \, .
	\end{equation}
	Depending on the slenderness ratio $\rho$, the rod has a circular cross-section of diameter $d = L / \rho$, radius $r=d/2$, area $A=\pi r^2$ and second moment of area $I = \tfrac{\pi}{4}r^4$. Using the Young's and shear moduli $E=1$ and $G=0.5$, respectively, the elastic constants are given by $k_\mathrm{e} = EA$, $k_\mathrm{s} = GA$, $k_{\mathrm{b}} = k_{\mathrm{b}_y} = k_{\mathrm{b}_z}= EI$ and $k_\mathrm{t} = 2GI$. As shown in Harsch et al.\cite{Harsch2021a}, the force boundary conditions for this specific example are found by a so called inverse procedure, see Section 5.2 of Ogden\cite{Ogden1997}. Evaluating the derivative of~\eqref{eq:analytical_solution_helix}, the tangent vector of the curve is
	\begin{equation}
		{}_I \vr_{OP,\xi}^{*}(\xi) = 
		R_0 \alpha_{,\xi} \begin{pmatrix}
			\cos\alpha(\xi) \\
			\sin\alpha(\xi) \\
			c
		\end{pmatrix} \, , \quad \alpha_{,\xi}(\xi) = 2 \pi n \, .
	\end{equation}
	The rod should not be elongated during the deformation, i.e., $J(\xi) = \|{}_I \vr_{OP,\xi}^{*}(\xi)
	\| = L$, from which the rod's total length $L=\sqrt{1 + c^2} R_0 2\pi n$ follows. Further, the Serret--Frenet frame of the deformed curve is given by
	\begin{equation}\label{eq:Serret-Frenet_frame_helix}
		\begin{aligned}
			{}_I \ve_x^K(\xi) = \frac{1}{\sqrt{1 + c^2}}
			\begin{pmatrix}
				\cos \alpha(\xi) \\
				\sin \alpha(\xi) \\
				c
			\end{pmatrix} \, , \quad
			{}_I \ve_y^K(\xi) =
			\begin{pmatrix}
				-\sin \alpha(\xi) \\
				\cos \alpha(\xi) \\
				0
			\end{pmatrix} \, , \quad
			{}_I \ve_z^K(\xi) = \frac{1}{\sqrt{1 + c^2}}
			\begin{pmatrix}
				-c \cos \alpha(\xi) \\
				-c \sin \alpha(\xi) \\
				1
			\end{pmatrix} \, .
		\end{aligned}
	\end{equation}
	Using $(1 + c^2)^{-\frac{1}{2}} = 2 \pi n R_0/L = \alpha_{,\xi} R_0 / L$ together with the derivatives 
	of~\eqref{eq:Serret-Frenet_frame_helix}, the 
	rod's strain measures~\eqref{eq:continuous_strain_measures} compute as
	\begin{equation}\label{eq:helix_constant_strain}
		{}_K \vga = (1, \ 0, \ 0) \ , \quad {}_K \vka_{IK} = (c, \ 0, \ 1) R_0 \alpha^2_{,\xi} / L^2 \, .
	\end{equation}
	For the given strain energy density~\eqref{eq:strain_energy_density}, evaluating the constitutive equations~\eqref{eq:constitutive_equations} results in the force boundary conditions
	\begin{equation}\label{eq:constitutive_relations_helix}
		{}_I \vb_1 = (0, \ 0, \ 0) \, , \quad {}_K \vc_1 = (k_\mathrm{t} c, \ 0, \ k_\mathrm{b}) R_0 \alpha^2_{,\xi} / L^2
	\end{equation}
	that induce a pure torsion and bending deformation. Inserting~\eqref{eq:constitutive_relations_helix} into the differential equation of the nonlinear Cosserat rod\cite{Eugster2020b} results in the requirement $k_\mathrm{t} = k_\mathrm{b}$. This condition is satisfied for the given problem setup. For a numerical simulation, the rod has to be clamped at ${}_I \vr_{OP}|_{\xi=0} = {}_I \vr_{OP}^{*}|_{\xi=0}$ with an orientation given by $\vA_{IK} |_{\xi=0} = ({}_I \ve_x^K, \ {}_I \ve_y^K, \ {}_I \ve_z^K) |_{\xi=0}$.
	
	Depending on the used slenderness ratio, a different number of force increments and error tolerances were chosen, see Table~\ref{tab:helix}. As in the \emph{cantilever experiment}, no locking is observed.  In fact, the $\SE(3)$-interpolation can exactly represent a helix, which is a configuration with constant strains \eqref{eq:helix_constant_strain}. As mentioned in Section~\ref{sec:SE(3)-interpolation}, a single $\SE(3)$-element is restricted to a relative rotation vector of absolute value bounded by $\pi$. Thus four elements can exactly represent the helix~\eqref{eq:analytical_solution_helix} with two coils. Nevertheless, during the Newton iterations locally the relative rotation vector might exceed an absolute value of $\pi$. Hence, we use five elements for the numerical investigation.
	\begin{table}
		\centering
		\caption{Experimental parameters used for the \emph{helix example}.}\label{tab:helix}
		\begin{tabular}{ccc}
			\toprule
			slenderness & tolerance & force increments \\
			\midrule
			$\SI{e1}{}$ & $\SI{e-8}{}$ & 70 \\
			$\SI{e2}{}$ & $\SI{e-9}{}$ & 100 \\
			$\SI{e3}{}$ & $\SI{e-10}{}$ & 200 \\
			$\SI{e4}{}$ & $\SI{e-11}{}$ & 500 \\
			\bottomrule
		\end{tabular}
	\end{table}
}
\subsection{Superimposed rotation of deformed cantilever}
\begin{figure}
	\centering
	\begin{subfigure}[c]{0.49\textwidth}
		\caption{}
		\includegraphics[width=0.9\textwidth, trim=10cm 14cm 10cm 14cm, clip]{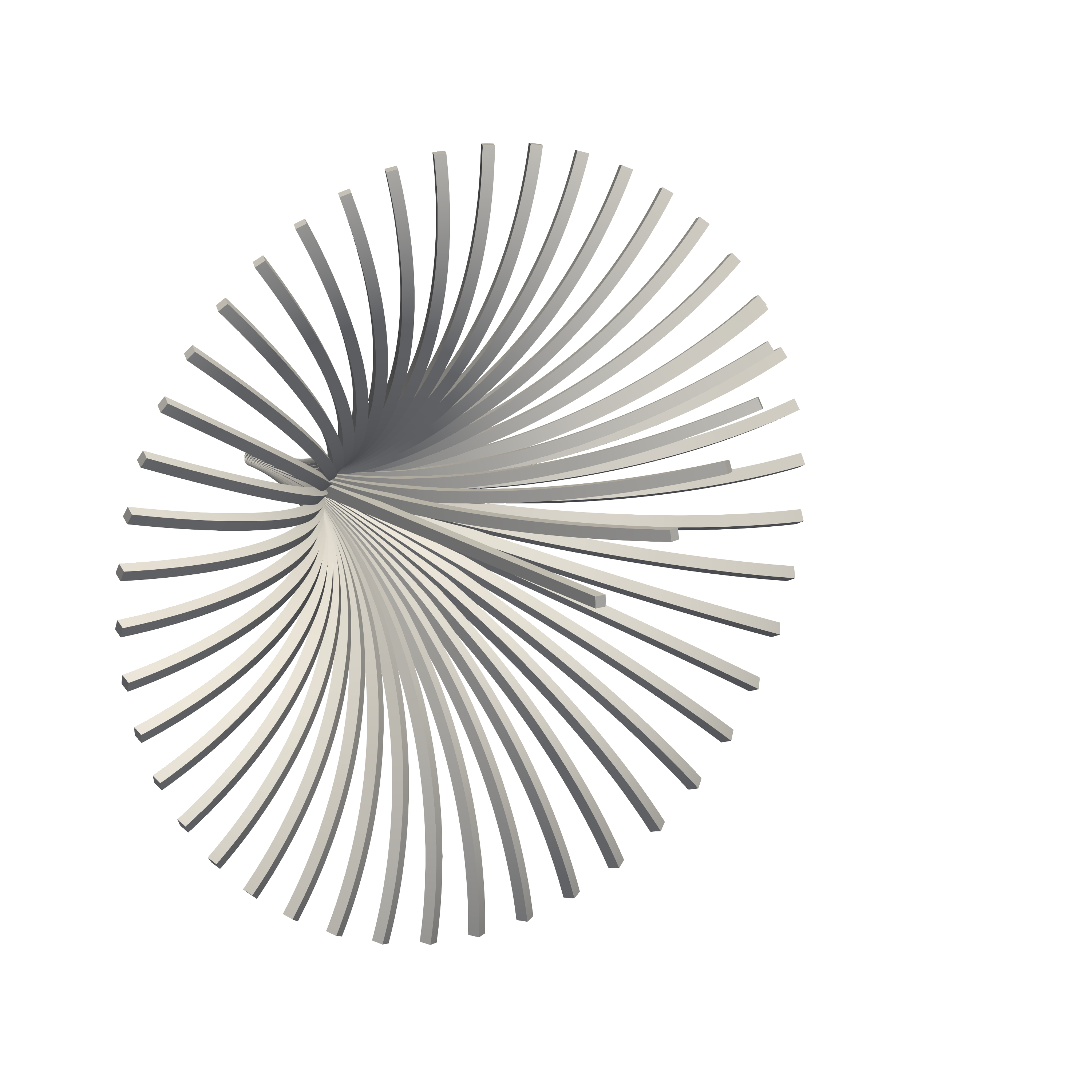}
	\end{subfigure}
	%
	\hspace{-0.75cm}
	\begin{subfigure}[c]{0.49\textwidth}
		\caption{}
		\includegraphics{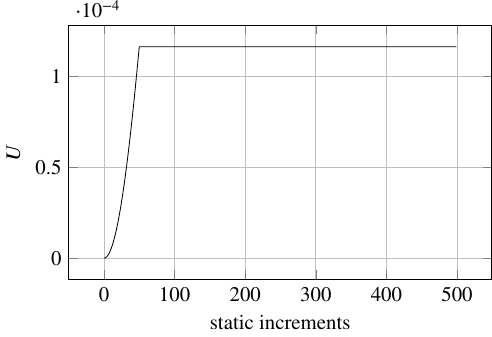}
	\end{subfigure}
	\caption{(a) Snapshots of deformed configurations of the rotated cantilever. (b) Potential energy $U$ vs. static increments.}
	\label{fig:objectivity_snapshot_tip_displacement}
\end{figure}

\begin{figure}
	\centering
	\begin{subfigure}[b]{0.49\textwidth}
		\caption{}
		\includegraphics{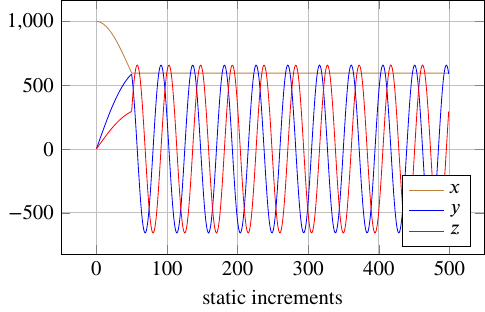}
	\end{subfigure}
	\hfill
	\begin{subfigure}[b]{0.49\textwidth}
		\caption{}
		\includegraphics{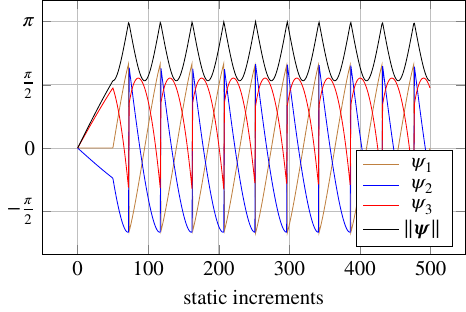}
	\end{subfigure}
	\caption{(a) Position vector of the second node vs. static increments. (b) Absolute rotation vector of the second node vs. static increments.}
	\label{fig:objectivity}
\end{figure}

The objectivity of the $\SE(3)$-formulation can  be demonstrated numerically by using the cantilever rod from Section~\ref{sec:quarter_circle} with slenderness ratio $\rho = \SI{e2}{}$, discretized by a single element. Both tip force and moment were successively applied during $50$ linearly spaced increments. During the subsequent $450$ increments the total rotation vector of the clamped node was linearly increased up to a value of $\vps_0 = (20 \pi, \ 0, \ 0)$, i.e., the deformed rod performed $10$ full rotations around the $\ve_x^I$-axis, see Figure~\ref{fig:objectivity_snapshot_tip_displacement} (a). Note, the very same problem can be solved using less than $100$ static increments. However, to get a better resolution of the derived quantities, we used a totality of $500$ static increments. It is well known that for non-objective finite element formulations, the total energy increases during a superimposed rigid body rotation\cite{Meier2014}. This is not the case for the present formulation, see Figure~\ref{fig:objectivity_snapshot_tip_displacement}~(b). After the final values of the tip force and moment are reached, the potential energy remains unaltered throughout the superimposed rotation of the clamped node. The tip displacement is shown in Figure~\ref{fig:objectivity}~(a). Moreover, Figure~\ref{fig:objectivity}~(b) visualizes that whenever the absolute value of the second nodal rotation vector would exceed $\pi$, its complement value is chosen. Thus, the  trajectory of the individual components of the nodal rotation vector exhibit discontinuities.
\subsection{Rod bent to a helical form}
\begin{figure}
	\centering
	\includegraphics{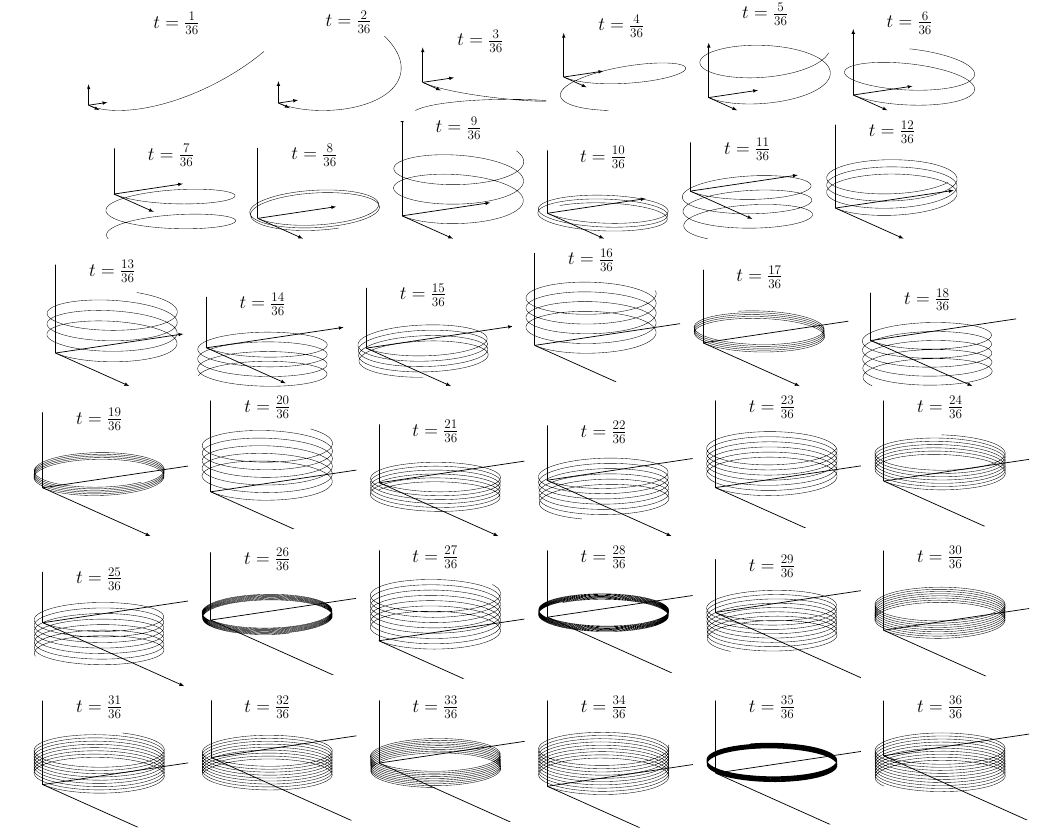}
	%
	%
	\caption{Linearly spaced snapshots of the helicoidal deformed rod centerline.}
	\label{fig:helix_snapshots}
\end{figure}

\begin{figure}
	\centering
	\begin{subfigure}[b]{0.32\textwidth}
		\caption{}
		\includegraphics{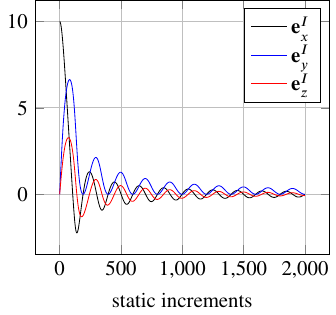}
	\end{subfigure}
	\hfill
	\begin{subfigure}[b]{0.32\textwidth}
		\caption{}
		\vspace{-0.4cm}
		\hspace{0.1cm}
		\includegraphics{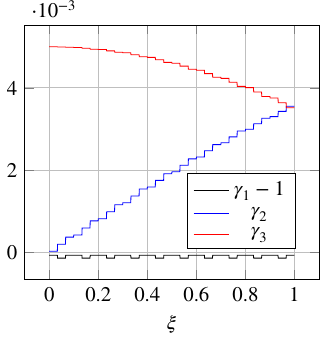}
	\end{subfigure}
	\hfill
	\begin{subfigure}[b]{0.32\textwidth}
		\caption{}
		\includegraphics{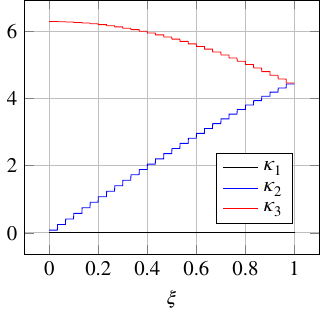}
	\end{subfigure}
	\caption{(a) Tip displacement vs. static increments. (b) and (c) Strain measures of the final configuration.}
	\label{fig:helix_tip_displacement}
\end{figure}

The next experiment investigates the capabilities of the presented formulation in describing large inhomogenous deformation. A well-suited example was introduced by Ibrahimbegovic\cite{Ibrahimbegovic1997} in which an initially straight cantilever rod of length $L=10$, axial stiffness $k_\mathrm{e} = \SI{e4}{}$, shear stiffness $k_\mathrm{s} = \SI{e4}{}$, bending stiffness $k_{\mathrm{b}} = k_{\mathrm{b}_y} = k_{\mathrm{b}_z} = \SI{e2}{}$ and torsional stiffness $k_\mathrm{t} = \SI{e2}{}$ is subjected to an incremental application of a tip moment ${}_K\vc_1 = \vA_{IK}\T (0, \ 0, \ 20 \pi k_\mathrm{b} / L)$ and an out-of-plane tip load ${}_I\vb_1= (0, \ 0, \ 50)$. By that, the rod is bend up to 10 full circles, while the $\ve_z^I$-component of the  tip displacement ``oscillates'' with decreasing magnitude, see Figure~\ref{fig:helix_snapshots}. This observation is in perfect agreement with the results found in literature and is further confirmed by comparing the tip displacements of the $\ve_z^I$-component, shown in Figure~\ref{fig:helix_tip_displacement} (a), with the results reported by Ibrahimbegovic\cite{Ibrahimbegovic1997}. Instead of the 100 or 200 linear finite elements used by Ibrahimbegovic\cite{Ibrahimbegovic1997} and Mäkinen\cite{Maekinen2007}, it was sufficient to discretize the rod with 30 elements of the present formulation.

Additionally, Figure~\ref{fig:helix_tip_displacement} (b) and (c) visualize the strain measures of the final configuration. Not surprisingly, the dilatation $\gamma_1$, both shear strains $\gamma_2$, $\gamma_3$ and both curvatures $\kappa_2$, $\kappa_3$ are step functions, since the used two-node element interpolation strategy~\eqref{eq:relative_interpolation_H_two_node} yields element-wise constant strain measures.	
\subsection{Flexible heavy top}
Inspired by the investigation of Mäkinen\cite{Maekinen2007}, the dynamics of an elastic heavy top is studied here. On the one hand, this example demonstrates the capability of the presented formulation to be solved using standard ODE solvers. On the other hand, in the limit case of an infinite stiff rod, the rod shows the well-known behavior of a heavy top. The motion of the heavy top is described by Euler's equations, see Equation (1.83) and (3.35) of Magnus~\cite{Magnus1971}, whose solution is taken as reference solution in the subsequent investigation. For high stiffnesses of the rod, we thus expect the solution to be close to the one of a rigid body.

Let the top be given by a cylinder of radius $r=0.1$ and length $L=0.5$ with cross-section area $A = \pi r^2$ and second moment of area $I = \pi r^4 / 4$. The cylinder is subjected to a constant distributed line force ${}_I \vb = \rho_0 A (0, \ 0, \ - 9.81)$. The stiff rod should be made out of steel with a uniform density $\rho_0 = 8000$, Young's modulus $E=\SI{210e6}{}$ and shear modulus $G = E / (2(1 + \nu))$, with a Poisson’s ratio $\nu = 1/3$. The stiff rod has consequently an axial stiffness $k_\mathrm{e} = EA$, shear stiffness $k_\mathrm{s} = GA$, bending stiffness $k_\mathrm{b} = k_{\mathrm{b}_y} = k_{\mathrm{b}_z} = EI$ and torsional stiffness $k_\mathrm{t} = 2GI$. We also considered a soft rod for which all stiffnesses were reduced by a factor $\SI{e3}{}$. The top was discretized using a single two-node $\SE(3)$-element. The initial position is such that the top points from the origin in positive $\ve_x^I$-direction, i.e., $\vq^0 = (\mathbf{0}_{3\times1}, \ \mathbf{0}_{3\times1}, \ \vr_1, \mathbf{0}_{3\times1})$ with $\vr_1 = (L, \ 0, \ 0)$. Its initial velocities are chosen such that in the case of a rigid rod a perfect precession motion\cite[Section 3.3.2 c)]{Magnus1971} is obtained, i.e., $\vu^0 = (\mathbf{0}_{3\times1}, \ \vOm, \ \vv_1, \vOm)$ with the tip velocity $\vv_1 = \vOm \times \vr_1$ and the angular velocity $\vOm = (\Omega, \ 0, \ \Omega_\mathrm{pr})$, where $\Omega = 50 \pi$ and $\Omega_\mathrm{pr} = g L / (r^2 \Omega)$. Finally, the motion of the top is constrained such that the first node coincides with the origin for all times. This can either be guaranteed by removing the corresponding degrees of freedom from the set of unknowns, or by using the concept of perfect bilateral constraints, cf. G\'eradin\cite{Geradin2001}.

\correctb{Using two different numerical time integration schemes (a standard fourth order Runge--Kutta method with the absolute and relative tolerances $atol = rtol = \SI{1e-8}{}$, see Hairer et al.\cite{Hairer1993} and the generalized-$\alpha$ method for first order differential equations proposed by Jansen et al.\cite{Jansen2000} with the spectral radius at infinity $\rho_\infty = 0.9$ and a step size $h=\SI{1e-5}{}$), the simulations were performed until a final time of $t_1 = 2 \pi / \Omega_\mathrm{pr}$ was reached, i.e., the rigid top performed a full rotation. Since the numerical solution of both methods only differ in their numerical digits, we only show the results obtained by the Runge--Kutta method.}

\begin{figure}
	\centering
	\begin{subfigure}[b]{0.47\textwidth}
		\caption{}
		\vspace{-0.15cm}
		\includegraphics{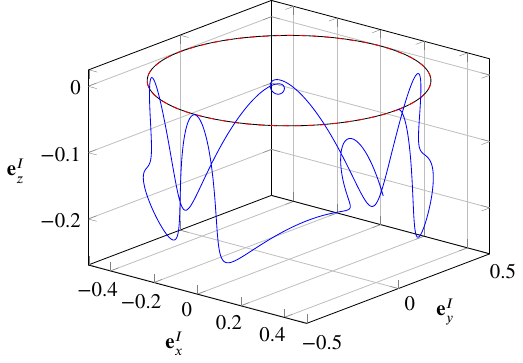}
	\end{subfigure}
	%
	\hspace{0.75cm}
	\begin{subfigure}[b]{0.47\textwidth}
		\caption{}
		\includegraphics{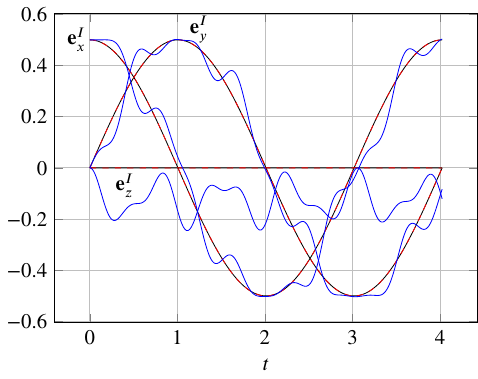}
	\end{subfigure}
	\caption{Tip displacement of the rigid top vs. stiff and soft rod solutions. The rigid top solution is drawn in black, the stiff rod in red and the soft rod in blue. (a) Spatial tip trajectory. (b) Tip displacements vs. time.}
	\label{fig:heavy_top}
\end{figure}

In Figure~\ref{fig:heavy_top} (a), the spatial trajectories of the different tops' free ends are shown. When comparing the projections of the rod's free tip, it can be seen that the assertion is true that the solution of a stiff rod cannot be distinguished from the rigid body solution, while the soft rod's tip performs a fascinating oscillatory motion superimposed to the rigid body solution.
\section{Conclusion}\label{sec:conclusion}
We have presented a rod finite element formulation based on a two-node $\SE(3)$-interpolation, that is, the interpolation of relative Euclidean transformation matrices using relative twists. In contrast to the typical uncoupled interpolation of centerline points and relative rotations using relative rotation vectors, the proposed $\SE(3)$-interpolation leads to element-wise constant strain measures. Thus, by construction, the resulting finite element formulation will show an absence of both membrane and shear locking and can be applied in scenarios where very high slenderness ratios are present. Objectivity of the discretized strain measures followed from the interpolation strategy with relative Euclidean transformation matrices.

With a Petrov--Galerkin projection method it is possible to discretize the arising virtual work functionals in space resulting in a set of ordinary differential equations (in dynamics) or a set of nonlinear equations (in statics) both of which can be solved using standard numerical schemes. Thus, the drawback of existing $\SE(3)$ rod formulations, which require highly specialized Lie group solvers is circumvented in an elegant way. Further, for typical applications, the arising mass matrix is symmetric and constant. Since the virtual displacements and rotations are interpolated instead of using the consistent variations of the proposed $\SE(3)$-interpolation formula, a non-symmetric stiffness matrix is obtained, which can be disadvantageous in terms of performance and storage.
The nodal kinematics of the proposed formulation, i.e., the nodal Euclidean transformation matrices, are parameterized using total centerline points and total rotation vectors together with the $\SO(3)$-exponential map. Arising singularities are circumvented by introducing the concept of complement rotation vectors together with the proposed relative interpolation strategy. Thus, a parametrization with a minimal number of six nodal unknowns has been achieved.
The paper is closed by demonstrating all the individual properties of the proposed formulation by numerical benchmark examples in statics and dynamics.
\appendix
\section{Matrix Lie groups}\label{sec:matrix_lie_groups}
\subsection{The general linear group $\GL(n)$}
In the course of this treatment, we restrict ourselves to matrix Lie groups over the real numbers, that is, Lie groups $G$ that are closed subgroups of the general linear group $\GL(n) = \{\vA \in \mR^{n \times n} | \det(\vA) \neq 0\}$, which is the group of all $n \times n$ invertible matrices with real entries\cite{Hall2015}. For $\vA, \vB \in \GL(n)$, the group operation is given by the usual matrix product $(\vA, \vB) \mapsto \vA \vB$. The inverse and transpose of $\vA$ are respectively denoted by $\vA^{-1}$ and $\vA\T$. The identity is given by the $n \times n$ identity matrix $\mathbf{1}_{n \times n}$. Analogously, we denote the $n \times n$ matrix containing only zeros by $\mathbf{0}_{n \times n}$.

Let $\gl(n) = T_{\mathbf{1}_{n \times n}} \GL(n)$ be the tangent space of $\GL(n)$ at the identity $\mathbf{1}_{n \times n}$. In~Section 3.3 of Baker\cite{Baker2002}, it is shown that $\gl(n)$ corresponds to the set of all real $n \times n$ matrices, which constitute an $n^2$-dimensional vector space. The Lie algebra of $\GL(n)$ consists of the vector space $\gl(n)$ together with the bilinear map
\begin{equation}
	[\bullet,\bullet] \colon \gl(n) \times \gl(n) \to \gl(n) \, , \quad (\vX, \vY) \mapsto [\vX, \vY] = \vX \vY - \vY \vX \, ,
\end{equation}
which is a Lie bracket as it is skew-symmetric and satisfies the Jacobi identity $[\vX, [\vY, \vZ]] + [\vY, [\vZ, \vX]] + [\vZ, [\vX, \vY]] = \mathbf{0}_{n \times n}$ $\forall \vX, \vY, \vZ \in \gl(n)$. For a fixed argument $\vX \in \gl(n)$ the Lie bracket defines the linear operator
\begin{equation}\label{eq:ad}
	\ad{\vX} \colon \gl(n) \to \gl(n) \, , \quad \vY \mapsto \ad{\vX}(\vY) = [\vX, \vY] \, ,
\end{equation}
called the adjoint map. Thus, instead of $[\vX, [\vX, [\vX, [\vX, \vY]]]]$, we can write $\ad{\vX}^4(\vY)$, see Hall\cite{Hall2015}~Section 3.1. Note that the adjoint map to the power of zero is defined as $\ad{\vX}^0(\vY) \coloneqq \vY$. Further, the identity $\ad{-\vX}^i(\vY) = (-1)^i\ad{\vX}^i(\vY)$ readily follows from iteratively applying the bilinearity of the Lie bracket. According to Baker\cite{Baker2002}~Section 2.1, the matrix exponential is defined by the power series
\begin{equation}\label{eq:exp}
	\exp \colon \gl(n) \to \GL(n) \, , \quad \vX \mapsto \vA = \exp(\vX) = \sum_{i=0}^{\infty} \frac{\vX^i}{i!} \, .
\end{equation}

As in~Section 2.6 of Iserles et al.\cite{Iserles2000}, for a smooth curve $\vZ(t) \in \gl(n)$, the derivative of the exponential map can be written as
\begin{equation}\label{eq:derivative_exponential_map}
	\td{}{t} \exp\left(\vZ(t)\right) = \dexp_{\vZ(t)}\left(\dot{\vZ}(t)\right) \exp\left(\vZ(t)\right) = \exp\left(\vZ(t)\right) \dexp_{-\vZ(t)}\left(\dot{\vZ}(t)\right) \, ,
\end{equation}
where, for $\vX \in \gl(n)$, the linear map $\dexp_{\vX}$ is defined in terms of a power series as
\begin{equation}\label{eq:dexp}
	\dexp_{\vX} \colon \gl(n) \to \gl(n) \, , \quad
	\vY \mapsto \dexp_{\vX}(\vY) \coloneqq \sum_{i=0}^{\infty} \frac{1}{(i + 1)!} \ad{\vX}^i(\vY) \, .
\end{equation}
Rearranging the equations in~\eqref{eq:derivative_exponential_map} to
\begin{equation}\label{eq:left_and_right_differential}
	\dexp_{-\vZ(t)}\left(\dot{\vZ}(t)\right) = \exp\left(\vZ(t)\right)^{-1} \td{}{t} \exp\left(\vZ(t)\right) \quad \text{and} \quad \dexp_{\vZ(t)}\left(\dot{\vZ}(t)\right) = \td{}{t} \exp\left(\vZ(t)\right) \exp\left(\vZ(t)\right)^{-1}  \, ,
\end{equation}
it is reasonable that $\dexp_{-\vX(t)}$ and $\dexp_{\vX(t)}$ are called left- and right-trivialized differentials, respectively. According to Equation~(2.46) of Iserles et al.\cite{Iserles2000}, the inverse of $\dexp_{\vX}$ can also be written by the power series
\begin{equation}\label{eq:dexp_inv}
	\begin{aligned}
		\dexp^{-1}_{\vX} \colon \gl(n) \to \gl(n) \, , \quad
		\vY \mapsto \dexp^{-1}_{\vX}(\vY) \coloneqq \sum_{i=0}^{\infty} \frac{B_i}{i!} \ad{\vX}^i(\vY) \, ,
	\end{aligned}
\end{equation}
where $B_i$ denotes the $i$-th Bernoulli number.

Let the tuple $(\mathfrak{g}, [\cdot, \cdot])$ be a subalgebra of $(\gl(n),  [\cdot, \cdot])$, i.e., $\mathfrak{g}$ is a subspace of $\gl(n)$ with dimension $k \leq n^2$ and $[\vX, \vY] \in \mathfrak{g}$ for $\vX, \vY \in \mathfrak{g}$. Then, we can always find a linear and bijective map
\begin{equation}\label{eq:j_map}
	j \colon \mR^k \to \mathfrak{g} \, , \quad \vx \mapsto \vX = j(\vx) \, .
\end{equation}
Since $j$ and $\dexp_{-j(\vx)}$ are linear maps, it is convenient to construct the tangent operator $\vT(\vx) \in \mR^{k \times k}$, which is the linear map defined as
\begin{equation}\label{eq:tangent_operator}
	\vT(\vx) \coloneqq j^{-1} \circ \dexp_{-j(\vx)} \circ j 
	= j^{-1} \circ \sum_{i=0}^{\infty} \frac{(-1)^i}{(i + 1)!} \ad{j(\vx)}^i \circ j
	= \sum_{i=0}^{\infty} \frac{(-1)^i}{(i + 1)!} \lceil \vx \rfloor^i \, ,
\end{equation}
where in the last equality we have introduced the mapping
\begin{equation}\label{eq:box}
	\lceil\bullet\rfloor
	\colon \mR^k \to \mR^{k \times k} \, , \quad \vx \mapsto \lceil \vx \rfloor = j^{-1} \circ \ad{j(\vx)} \circ j \, .
\end{equation}
Often this mapping is applied in the form $j\big(\lceil \vx \rfloor \vy\big) = \ad{j(\vx)}\big(j(\vy)\big)$, for $\vx, \vy \in \mR^k$. If $j(\vx)$ and $j(\vy)$ commute, then powers of the adjoint map $\ad{j(\vx)}^i(j(\vy))$ vanish for $i > 0$ resulting in $\vT(\vx) \vy = \left(j^{-1} \circ \dexp_{-j(\vx)} \circ j\right) \vy = \vy$. The inverse of the tangent operator is given by
\begin{equation}\label{eq:inverse_tangent_operator}
	\vT(\vx)^{-1} = \left(j^{-1} \circ \dexp_{-j(\vx)} \circ j\right)^{-1} = j^{-1} \circ \dexp^{-1}_{-j(\vx)} \circ j 
	= j^{-1} \circ \sum_{i=0}^{\infty} \frac{(-1)^i B_i}{i!} \ad{j(\vx)}^i \circ j
	= \sum_{i=0}^{\infty} \frac{(-1)^i B_i}{i!} \lceil \vx \rfloor^i \, .
\end{equation}
Since the matrix rank of $\lceil \vx \rfloor$ is at most $k$, the Cayley--Hamilton theorem guarantees the existence of functions $a_0(\lceil \vx \rfloor)$, $\dots$, $a_{k-1}(\lceil \vx \rfloor)$ such that
\begin{equation}
	\lceil \vx \rfloor^{k} = \sum_{i=0}^{k - 1} a_i(\lceil \vx \rfloor) \lceil \vx \rfloor^i \, .
\end{equation}
Hence, we can express the tangent operator~\eqref{eq:tangent_operator} and its inverse~\eqref{eq:inverse_tangent_operator} by a finite number of powers of $\lceil \vx \rfloor$. Using similar arguments, also the exponential map~\eqref{eq:exp} can be expressed by a finite number of matrix powers.
\subsection{The special orthogonal group $\SO(3)$}
\label{sec:SO3}
As a subgroup of $\GL(3)$, the special orthogonal group 
\begin{equation}
	\SO(3) = \{ \vA \in \GL(3) | \vA\T \vA = \mathbf{1}_{3\times3} , \, \operatorname{det}(\vA) = +1 \} 
\end{equation}
is the group of all proper orthogonal transformations. In Baker\cite{Baker2002}~Section 3.3, it is shown that the tangent space at the identity 
\begin{equation}
	\so(3) = T_{\mathbf{1}_{3 \times 3}}\SO(3) = \{\vOm \in \mR^{3 \times 3} | \vOm + \vOm\T = \mathbf{0}_{3 \times 3}\}
\end{equation} 
corresponds to the three-dimensional space of skew-symmetric $3 \times 3$ matrices, which in dynamics often corresponds to the space of angular velocities. Since $\so(3)$ is a subspace of $\gl(3)$ and the Lie bracket $[\vOm_1, \vOm_2] \in \so(3)$ for $\vOm_1, \vOm_2 \in \so(3)$, the tuple $(\so(3), [\cdot, \cdot])$ is a Lie subalgebra of $(\gl(3), [\cdot,\cdot])$.

For $\so(3)$, the map~\eqref{eq:j_map} can be identified explicitly as
\begin{equation}\label{eq:hat_so(3)}
	j_{\SO(3)}(\bullet) = \widetilde{(\bullet)}
	\colon \mR^3 \to \so(3) \, , \quad \vom =
	\begin{pmatrix}
		\omega_1 \\ \omega_2 \\ \omega_3
	\end{pmatrix} 
	\mapsto \vOm = j_{\SO(3)}(\vom) = \widetilde{\vom} 
	= \begin{pmatrix}
		0 & -\omega_3 & \omega_2 \\
		\omega_3 & 0 & -\omega_1 \\
		-\omega_2 & \omega_1 & 0
	\end{pmatrix} \, ,
\end{equation}
where we have introduced the tilde symbol for compact notation. This map is related to the cross product by $\widetilde{\vom} \vr = \vom \times \vr$ for $\vom, \vr \in \mR^3$.

In order to minimize the number of transcendental function evaluations\cite{Park2005, Sonneville2014}, we introduce the mappings
\begin{equation}\label{eq:auxiliary_alpha_beta}
	\alpha(\widetilde{\vom}) = \frac{\sin(\|\vom\|)}{\|\vom\|} \, , \quad \beta(\widetilde{\vom}) = 2\frac{1 - \cos(\|\vom\|)}{\|\vom\|^2} \, , \quad \gamma(\widetilde{\vom}) = \frac{\alpha(\widetilde{\vom})}{\beta(\widetilde{\vom})} = \frac{\|\vom\|}{2} \cot(\|\vom\| / 2) \, ,
\end{equation}
with $\lim_{\|\vom\| \to 0} \alpha(\widetilde{\vom}) = \lim_{\|\vom\| \to 0} \beta(\widetilde{\vom}) = \lim_{\|\vom\| \to 0} \gamma(\widetilde{\vom}) = 1$. Noting that $\widetilde{\vom}^3 = - \|\vom\|^2 \widetilde{\vom}$ and carefully separating even and odd parts of the power series of the exponential map~\eqref{eq:exp}, see Murray et al.\cite{Murray1994}~Section 2.2, the exponential map from $\so(3)$ to $\SO(3)$ has an analytical form known as Rodrigues' formula
\begin{equation}\label{eq:exp_SO3}
	\exp_{\SO(3)} \colon \so(3) \to \SO(3) \, , \quad
	\vOm \mapsto \vA = \exp_{\SO(3)}(\vOm) = \mathbf{1}_{3  \times 3} + \alpha(\vOm) \vOm + \frac{\beta(\vOm)}{2} \vOm^2 \, .
\end{equation}
For the case $\|\vOm\| \to 0$, following the typical approach in computational methods, the first order approximation $\exp_{\SO(3)}(\vOm) = \mathbf{1}_{3  \times 3} + \vOm$ is applied.

From explicitly computing~\eqref{eq:exp_SO3} and denoting $\omega = \|\vom\|$, the diagonal terms reveal the identity $\cos\omega =  \tfrac{1}{2}(\operatorname{tr}(\vA) - 1)$, which can be solved for $\omega = \omega(\vA)$. The off diagonals directly lead to $\widetilde{\vom} = \tfrac{\omega}{2 \sin(\omega)}(\vA - \vA\T)$. Using $\sin^2(\omega(\vA)) = 1 - (\operatorname{tr}(\vA) - 1)^2 / 4$, an analytic form of the $\SO(3)$-logarithm map is given by
\begin{equation}\label{eq:log_SO3}
	\log_{\SO(3)} \colon \SO(3) \to \so(3) \, , \quad
	\vA \mapsto \vOm = \log_{\SO(3)}(\vA) = \frac{\omega(\vA)}{2 \sin\big(\omega(\vA)\big)} (\vA - \vA\T) \, .
\end{equation}
Again, for $\omega \to 0$, we use the first order approximation $\vOm = (\vA - \vA\T) / 2$. For notational convenience, we further introduce the capitalized exponential and logarithm maps
\begin{equation}\label{eq:ExpLog_SO3}
	\begin{aligned}
		&\Exp_{ \SO(3)} \colon \mR^3 \to \SO(3) \, , \quad &\vom \mapsto \vA &= \Exp_{\SO(3)}(\vom) = \exp_{\SO(3)} \circ j_{\SO(3)} (\vom) \, , \\
		&\Log_{ \SO(3)} \colon \SO(3) \to \mR^3 \, , \quad 
		&\vA \mapsto \vom &= \Log_{\SO(3)}(\vA) = j_{\SO(3)}^{-1} \circ \log_{\SO(3)} (\vA) \\
		& &&= \frac{\omega(\vA)}{2 \sin\big(\omega(\vA)\big)} 
		\big(A_{32} - A_{23}, \, A_{13} - A_{31}, \, A_{21} - A_{12}\big) \, ,
	\end{aligned}
\end{equation}
which directly relate $\mR^3$ with $\SO(3)$ and vice versa.

For a lighter notation, we shortly suppress the subscripts of $j$ and $\lceil \bullet \rfloor$ indicating $\SO(3)$. Let $\vx, \vy, \vz \in \mR^3$, using the skew-symmetry and Jacobi's identity of the cross product, direct computation verifies $j\big(\lceil \vx \rfloor \vy\big) \vz = \ad{j(\vx)}\big(j(\vy)\big) \vz = \big(j(\vx) j(\vy) - j(\vy) j(\vx)\big) \vz = \big(\widetilde{\vx} \, \widetilde{\vy} - \widetilde{\vy} \, \widetilde{\vx}\big) \vz = j(\widetilde{\vx}\vy) \vz$. Comparison of both sides of the equality, readily implies $\lceil \vx \rfloor_{\SO(3)} = j_{\SO(3)}(\vx) = \widetilde{\vx}$. Again using $\widetilde{\vom}^3 = - \|\vom\|^2 \widetilde{\vom}$ and carefully separating the terms in the power series \eqref{eq:tangent_operator} (see Iserles et al.\cite{Iserles2000}~Equation~(B.10) and Equation~(18) of Park and Chung\cite{Park2005})\footnote{Note, both references define the tangent operator and its inverse in terms of the right-trivialized differential. Thus, their formulas differ in the sign of the skew-symmetric part.}, the $\SO(3)$-tangent operator can be written in the form
\begin{equation}\label{eq:tangent_map_so(3)}
	\vT_{\SO(3)}(\vom) = \sum_{i=0}^{\infty} \frac{(-1)^i}{(i + 1)!} \lceil \vom \rfloor^i_{\SO(3)} = \mathbf{1}_{3 \times 3} - \frac{\beta(\widetilde{\vom})}{2} \widetilde{\vom} + \frac{1 - \alpha(\widetilde{\vom})}{\|\vom\|^2} \widetilde{\vom}^2 \, .
\end{equation}
Similarly from \eqref{eq:inverse_tangent_operator} (see Iserles et al.\cite{Iserles2000}~Equation~(B.11), Park and Chung\cite{Park2005}~Equation~(19) and Equation~(2.20) of Bullo\cite{Bullo1995}), the inverse $\SO(3)$-tangent operator can efficiently be computed via
\begin{equation}\label{eq:inverse_tangent_map_so(3)}
	\vT^{-1}_{\SO(3)}(\vom) = \sum_{i=0}^{\infty} \frac{B_i}{i!} \lceil \vom \rfloor^i_{\SO(3)} = \mathbf{1}_{3 \times 3} + \frac{1}{2} \widetilde{\vom} + \frac{1 - \gamma(\widetilde{\vom})}{\|\vom\|^2}\widetilde{\vom}^2 \, .
\end{equation}
For $\|\vom\| \to 0$, we use the first order approximations $\vT_{\SO(3)}(\vom) = \mathbf{1}_{3 \times 3} - \tfrac{1}{2} \widetilde{\vom}$ and $\vT^{-1}_{\SO(3)}(\vom) = \mathbf{1}_{3 \times 3} + \frac{1}{2} \widetilde{\vom}$. Due to the skew-symmetry of  $\widetilde{\vom}$, the tangent map and its inverse fulfill the additional properties $\vT_{\SO(3)}(-\vom) = \vT_{\SO(3)}\T(\vom)$ and $\vT^{-1}_{\SO(3)}(-\vom) = \vT_{\SO(3)}^{-\mathrm{T}}(\vom)$. 
\subsection{The special Euclidean group $\SE(3)$}
As a subgroup of $\GL(4)$, the special Euclidean group 
\begin{equation}
	\SE(3) = \SO(3) \ltimes \mR^3  = \left\{
	\vH =
	\begin{pmatrix}
		\vA & \vr \\
		\mathbf{0}_{1 \times 3} & 1
	\end{pmatrix}  \in \GL(4) \, \bigg | \, \vA \in \SO(3) , \, \vr \in \mR^3
	\right\}
\end{equation}
is the the group of all Euclidean transformations $\vH$, which can be expressed in terms of a translation $\vr \in \mR^3$ and a proper orthogonal transformation $\vA \in \SO(3)$. Its linearization at the identity 
\begin{equation}
	\se(3) = T_{\mathbf{1}_{4 \times 4}} \SE(3) = \so(3) \ltimes \mR^3 = \left\{
	\vTh =
	\begin{pmatrix}
		\vOm & \vv \\
		\mathbf{0}_{1 \times 1} & 0
	\end{pmatrix} \in \mR^{4 \times 4} \bigg | \, \vOm \in \so(3) , \, \vv \in \mR^3
	\right\}
\end{equation}
is the space of all twists $\vTh$, which is a $4 \times 4$ matrix composed of a translational and angular velocity $\vv \in \mR^3$ and $\vOm \in \so(3)$, respectively. Since $\se(3)$ is a subgroup of $\gl(4)$ and the Lie bracket $[\vTh_1, \vTh_2] \in \se(3)$ for $\vTh_1, \vTh_2 \in \se(3)$, the tuple $(\se(3), [\cdot,\cdot])$ is a Lie subalgebra of $(\gl(4), [\cdot,\cdot])$.

For $\se(3)$, the map~\eqref{eq:j_map} can be identified explicitly as
\begin{equation}\label{eq:hat_se(3)}
	j_{\SE(3)}	\colon \mR^6 \to \se(3) \, , \quad \vth = 
	\begin{pmatrix}
		\vv \\ 
		\vom
	\end{pmatrix} \mapsto \vTh = j_{\SE(3)}(\vth) = 
	\begin{pmatrix}
		j_{\SO(3)}(\vom) & \vv \\
		\mathbf{0}_{1 \times 3} & 0
	\end{pmatrix} 
	= \begin{pmatrix}
		\widetilde{\vom} & \vv \\
		\mathbf{0}_{1 \times 3} & 0
	\end{pmatrix} \, .
\end{equation}

The exponential map from $\se(3)$ to $\SE(3)$ is given by~(see Example A.12 of Murray et al.\cite{Murray1994})
\begin{equation}\label{eq:exp_SE3}
	\exp_{\SE(3)} \colon \se(3) \to \SE(3) \, , \quad
	\vTh \mapsto \exp_{\SE(3)} (\vTh)
	= 
	\begin{pmatrix}
		\exp_{\SO(3)}(\vOm) & \quad \vT_{\SO(3)}\T(\vom) \vv \\
		\mathbf{0}_{1 \times 3} & 1
	\end{pmatrix} \, ,
\end{equation}
where $\vom = j^{-1}_{\SO(3)} (\vOm)$. By computing $\vTh = \log_{\SE(3)} \circ \exp_{\SE(3)}(\vTh)$ it can be verified that the $\SE(3)$-logarithm map is given in terms of $\SO(3)$ formulas only and reads as
\begin{equation}\label{eq:log_SE3}
	\log_{\SE(3)} \colon \SE(3) \to \se(3) \, , \quad
	\vH =
	\begin{pmatrix}
		\vA & \vr \\
		\mathbf{0}_{1 \times 3} & 1
	\end{pmatrix} \mapsto \vTh = 
	\begin{pmatrix}
		\vOm & \quad\vT_{\SO(3)}^{-\mathrm{T}}(\vom) \vr \\
		\mathbf{0}_{1 \times 3} & 0
	\end{pmatrix} \, ,
\end{equation}
where $\vOm = \log_{\SO(3)}(\vA)$ and $\vom = j_{\SO(3)}^{-1}(\vOm)$. Again, the capitalized exponential and logarithm maps relate $\mR^6$ with $\SE(3)$ and vice versa via
\begin{equation}\label{eq:ExpLog_SE3}
	\begin{aligned}
		&\Exp_{ \SE(3)} \colon \mR^6 \to \SE(3) \, , \quad &\vth \mapsto \vH &= \Exp_{ \SE(3)}(\vth) = \exp_{\SE(3)} \circ j_{\SE(3)} (\vth) \, , \\
		&\Log_{ \SE(3)} \colon \SE(3) \to \mR^6 \, , \quad &\vH \mapsto
		\vth &= \Log_{ \SE(3)}(\vH) = j_{\SE(3)}^{-1} \circ \log_{\SE(3)} (\vH) \\
		& &&= \begin{pmatrix}
			\vT_{\SO(3)}^{-\mathrm{T}}\big(\Log_{\SO(3)}(\vA)\big) \vr \\
			\Log_{\SO(3)}(\vA)
		\end{pmatrix} \, .
	\end{aligned}
\end{equation}
\section{Variation of strain measures}\label{sec:variations_strain_measures}
Using the skew-symmetry of ${}_K \delta \widetilde{\vph}_{IK}$, the expression \eqref{eq:virtual_rotation}$_2$ can be written as $\delta \vA_{IK}\T = -{}_K \delta \widetilde{\vph}_{IK} \vA_{IK}\T$. Consequently, the variation of ${}_K \bar{\vga}$ from \eqref{eq:continuous_strain_measures} computes to
\begin{equation}\label{eq:delta_K_gamma}
	\delta ({}_K \bar{\vga}) = \delta (\vA_{IK}\T {}_I \vr_{OP, \xi}) = \vA_{IK}\T \, \delta \left({}_I \vr_{OP, \xi}\right) + \delta \vA_{IK}\T {}_I \vr_{OP, \xi} = \vA_{IK}\T ({}_I \delta \vr_{P})_{,\xi} - {}_K \delta \vph_{IK} \times {}_K \bar{\vga} \, .
\end{equation}
In the last equality, we have recognized the virtual displacement \eqref{eq:virtual_displacement_velocity}$_1$ since the variation and $\xi$-derivative commute for vector components with respect to a constant basis. The same property gives rise to the equality
\begin{equation}\label{eq:delta_xi_interchangable}
	\mathbf{0} = \vA_{IK}\T \left[ 
	\left(\delta \left({}_I \ve_i^K\right)\right)_{,\xi} - \delta\left(\left({}_I \ve_i^K\right)_{,\xi}\right)
	\right] \, ,
\end{equation}
where ${}_I \ve_i^K = \vA_{IK} {}_K \ve_i^K$, $i \in \{x, y, z\}$. Straightforward computation together with~\eqref{eq:curvature} and~\eqref{eq:virtual_rotation} yields
\begin{equation}
	\begin{aligned}
		\vA_{IK}\T \left(\delta \left({}_I \ve_i^K\right)\right)_{,\xi} &= \vA_{IK}\T \left(\delta \vA_{IK} {}_K \ve_i^K\right)_{,\xi} \\
		&= \vA_{IK}\T \left[\vA_{IK} \left( {}_K \delta \vph_{IK} \times {}_K \ve_i^K\right)\right]_{,\xi} \\
		&= {}_K \bar{\vka}_{IK} \times \left( {}_K \delta \vph_{IK} \times {}_K \ve_i^K\right) + \left( {}_K \delta \vph_{IK}\right)_{,\xi} \times {}_K \ve_i^K
	\end{aligned}
\end{equation}
as well as
\begin{equation}
	\begin{aligned}
		-\vA_{IK}\T \delta \left(\left({}_I \ve_i^K\right)_{,\xi}\right) &= -\vA_{IK}\T \delta \left(\vA_{IK, \xi} {}_K \ve_i^K \right) \\
		&= -\vA_{IK}\T \delta \left[\vA_{IK} \left({}_K \bar{\vka}_{IK} \times {}_K \ve_i^K \right)\right] \\
		&= -{}_K \delta \vph_{IK} \times \left({}_K \bar{\vka}_{IK} \times {}_K \ve_i^K \right) - \delta \left({}_K \bar{\vka}_{IK}\right) \times {}_K \ve_i^K \\
		&= {}_K \delta \vph_{IK} \times \left({}_K \ve_i^K \times {}_K \bar{\vka}_{IK}  \right) - \delta \left({}_K \bar{\vka}_{IK}\right) \times {}_K \ve_i^K \, .
	\end{aligned}
\end{equation}
Inserting the latter two expressions in~\eqref{eq:delta_xi_interchangable} and making use of the Jacobi identity and the skew-symmetry of the cross product, one readily sees that
\begin{equation}\label{eq:delta_xi_interchangable2}
	\left( \left({}_K \delta\vph_{IK}\right)_{,\xi} - {}_K \delta\vph_{IK} \times {}_K \bar{\vka}_{IK} - \delta \left({}_K \bar{\vka}_{IK}\right)\right) \times {}_K \ve_{i}^K = \mathbf{0} \, .
\end{equation}
As this holds for all $i \in \{x, y, z\}$, the term in the bracket of~\eqref{eq:delta_xi_interchangable2} must vanish, resulting in the identity
\begin{equation}\label{eq:delta_K_kappa}
	\delta \left({}_K \bar{\vka}_{IK}\right) = \left({}_K \delta\vph_{IK}\right)_{,\xi} - {}_K \delta\vph_{IK} \times {}_K \bar{\vka}_{IK} \, .
\end{equation}
\section{Linearization of internal forces}
\label{sec:linearization}
In order to evaluate the linearization of the internal force vector $\vf^\mathrm{int}$, the derivatives of the exponential and logarithm maps of the underlying parametrization are required. Since most of these derivatives are not well established in literature and in order to be self-consistent, we briefly introduce them in the subsequent treatment. Most of the presented formulas have a removable singularity for $\omega = \|\vom\| \to 0$. Thus, during a numerical implementation a critical angle $\omega_\mathrm{crit} = \SI{1e-6}{}$ is defined. Whenever the value of $\omega$ falls below this critical angle, the first order approximations introduced in the first section are used. Thus, we have to introduce both, the derivative of the required formulas and the corresponding first order approximation.

In the following, we identify the $i$-th component of a tuple $\va$ by $a_i$, the $i,j$-th component of the matrix $\vA$ by $A_{ij}$ and the components of third order objects are indicated by three subscripts. Using the identities $\widetilde{\omega}_{ij} = \omega_k \varepsilon_{kji}$, $\partial\widetilde{\omega}_{ij} / \partial \omega_k = \varepsilon_{kji}$ and $\partial(\widetilde{\omega}_{il} \widetilde{\omega}_{lj}) / \partial \omega_k = \varepsilon_{kli} \widetilde{\omega}_{lj} + \widetilde{\omega}_{il} \varepsilon_{kjl}$, where $\varepsilon_{ijk}$ denotes the Levi--Civita permutation symbol, we can compute the derivatives of the auxiliary functions introduced in~\eqref{eq:auxiliary_alpha_beta}, with respect to $\omega_k$ by
\begin{equation}
	\pd{\alpha(\widetilde{\vom})}{\vom} = \frac{\vom}{\|\vom\|^2} \left(\cos(\omega) - \alpha(\widetilde{\vom})\right) \, , \quad
	\pd{\beta(\widetilde{\vom})}{\vom} = 2 \frac{\vom}{\|\vom\|^2} \left(\alpha(\widetilde{\vom}) - \beta(\widetilde{\vom})\right) \, , \quad 
	\pd{\gamma(\widetilde{\vom})}{\vom} = \vom \left(\frac{\gamma(\widetilde{\vom})}{\|\vom\|^2} \big(1 - \gamma(\widetilde{\vom})\big) - 1\right) \, .
\end{equation}

Thereby, the derivative of the rotation matrix components $A_{ij}$, obtained from the capitalized $\SO(3)$ exponential map~\eqref{eq:ExpLog_SO3}, with respect to $\omega_k$ is
\begin{equation}\label{eq:derivative_Exp_SO3}
	\left[ \pd{\Exp_{\SO(3)}}{\vom}(\vom) \right]_{ijk} =
	\begin{cases}
		-\alpha \varepsilon_{ijk} + \frac{\cos\omega - \alpha}{\omega^2} \widetilde{\omega}_{ij} \omega_k + \frac{\alpha - \beta}{\omega^2} \widetilde{\omega}_{il} \widetilde{\omega}_{lj} \omega_k + \frac{\beta}{2} \left(\varepsilon_{kli} \widetilde{\omega}_{lj} + \widetilde{\omega}_{il} \varepsilon_{kjl}\right) \, , & \omega > \omega_\mathrm{crit} \, , \\
		-\varepsilon_{ijk} \, , & \omega \leq \omega_\mathrm{crit} \, . \\
	\end{cases}
\end{equation}
Similar relations are found in Gallego and Yezzi\cite{Gallego2015}. 

Multiplying $\widetilde{\omega}_{nm} = \varepsilon_{kmn} \omega_{k}$ with $\varepsilon_{imn}$ and using $\varepsilon_{kmn} \varepsilon_{imn} = 2 \delta_{ki}$ gives $\omega_i = \frac{1}{2} \varepsilon_{imn} \widetilde{\omega}_{nm}$. Thus, we have to compute
\begin{equation}
	\left[ \pd{\Log_{\SO(3)}}{\vA}(\vA) \right]_{ijk} = \frac{1}{2} \varepsilon_{imn} \left[ \pd{\log_{\SO(3)}}{\vA}(\vA) \right]_{nmjk} \, .
\end{equation}
Further, from $\cos\omega =  \tfrac{1}{2}(\operatorname{tr}(\vA) - 1)$ with $\partial\arccos(x) /\partial x = - 1 / \sqrt{1 - x^2}$ we find
\begin{equation}
	\left[\pd{\omega(\vA)}{\vA}\right]_{jk} = - \frac{\delta_{jk}}{2 \sin\omega(\vA)} \quad \text{and} \quad \left[\pd{}{\vA}\left(\frac{\omega(\vA)}{2 \sin\omega(\vA)}\right)\right]_{jk} = \frac{\omega(\vA) \cos\omega(\vA) - \sin\omega(\vA)}{2 \sin^3\omega(\vA)} \delta_{jk} \, .
\end{equation}
Finally, using the identities $\partial A_{ab} / \partial A_{cd} = \delta_{ac} \delta_{bd}$ and $\varepsilon_{imn} (\delta_{mj}\delta_{nk} - \delta_{nj}\delta_{mk}) = 2 \varepsilon_{ijk}$, the derivative of $\omega_i$, obtained from the capitalized $\SO(3)$ logarithm map~\eqref{eq:ExpLog_SO3}, with respect to the rotation matrix components $A_{jk}$ is given by
\begin{equation}\label{eq:derivative_Log_SO3}
	\left[ \pd{\Log_{\SO(3)}}{\vA}(\vA) \right]_{ijk} =
	\begin{cases}
		\frac{\omega \cos\omega - \sin\omega}{4\sin^3\omega} \varepsilon_{imn} (A_{nm} - A_{mn}) \delta_{jk} + \frac{\omega}{2\sin\omega} \varepsilon_{ijk} \, , & \omega > \omega_\mathrm{crit} \, , \\
		\frac{1}{2} \varepsilon_{ijk} \, , & \omega \leq \omega_\mathrm{crit} \, .
	\end{cases}
\end{equation}

Further, the derivatives of the $\SO(3)$ tangent operators are required. Using the identities introduced above they can be computed as
\begin{equation}
	\left[ \pd{\vT_{\SO(3)}}{\vom}(\vom) \right]_{ijk} = 
	\begin{cases}
		\frac{\beta}{2} \varepsilon_{ijk} + \widetilde{\omega}_{ij} \omega_k \frac{\beta - \alpha}{\omega^2} + \frac{1 - \alpha}{\omega^2} \left(\varepsilon_{kli} \widetilde{\omega}_{lj} + \widetilde{\omega}_{il} \varepsilon_{kjl}\right) + \frac{3 \alpha - 2 - \cos\omega}{\omega^4} \widetilde{\omega}_{il} \widetilde{\omega}_{lj} \omega_k \, , & \omega > \omega_\mathrm{crit} \, , \\
		\frac{1}{2} \varepsilon_{ijk} \, , & \omega \leq \omega_\mathrm{crit}
	\end{cases}
\end{equation}
and
\begin{equation}
	\left[ \pd{\vT^{-1}_{\SO(3)}}{\vom}(\vom) \right]_{ijk} = 
	\begin{cases}
		-\frac{1}{2} \varepsilon_{ijk} + \frac{1 - \gamma}{\omega^2} \left(\varepsilon_{kli} \widetilde{\omega}_{lj} + \widetilde{\omega}_{il} \varepsilon_{kjl}\right) + \frac{1}{\omega^2} \left(
		2 \frac{1 - \gamma}{\omega^2} + \frac{\gamma}{\omega^2} (1 - \gamma) - \frac{1}{4}
		\right) \widetilde{\omega}_{il} \widetilde{\omega}_{lj} \omega_k \, , & \omega > \omega_\mathrm{crit} \, , \\
		-\frac{1}{2} \varepsilon_{ijk} \, , & \omega \leq \omega_\mathrm{crit} \, .
	\end{cases}
\end{equation}

The derivatives of the capitalized $\SE(3)$ exponential and logarithm introduced in~\eqref{eq:ExpLog_SE3} boils down to correctly assemble the just computed derivatives. Recalling that $\vth = (\vv , \ \vom)$, we get the derivative of the exponential map
\begin{equation}
	\left[ \pd{\Exp_{ \SE(3)}}{\vth}(\vth)\right]_{ijk} = 
	\begin{cases}
		\left[\pd{\Exp_{ \SO(3)}}{\vom}(\vom)\right]_{ij(k-3)} \, , & \mathrm{for} \quad 1 \leq i, j \leq 3 \, , \quad 4 \leq k \leq 6 \, , \\
		\left[\vT\T_{ \SO(3)}(\vom) \right]_{ik} \, , & \mathrm{for} \quad 1 \leq i \leq 3 \, , \quad j = 4 \, , \quad 1 \leq k \leq 3 \, , \\
		\left[\vv\right]_l \left[\pd{\vT_{ \SO(3)}}{\vom}(\vom)\right]_{li(k-3)} \, , & \mathrm{for} \quad 1 \leq i \leq 3 \, , \quad j = 4 \, , \quad 4 \leq k \leq 6 \, , \\
		0 \, , & \mathrm{else} \, .
	\end{cases}
\end{equation}
With 
\begin{equation}
	\vH = 
	\begin{pmatrix}
		\vA & \vr \\
		\mathbf{0}_{1 \times 3} & 1
	\end{pmatrix}
\end{equation}
the derivative of the $\SE(3)$ logarithm map with respect to its argument is given by
\begin{equation}
	\left[ \pd{\Log_{ \SE(3)}}{\vH}(\vH)\right]_{ijk} = 
	\begin{cases}
		\left[\vr\right]_l \left[\pd{\vT^{-1}_{ \SO(3)}}{\vom}(\vom)\right]_{lim} \left[\pd{\Log_{\SO(3)}}{\vA}(\vA)\right]_{mjk} \, , & \mathrm{for} \quad 1 \leq i, j, k \leq 3 \, , \\
		\left[\vT_{\SO(3)}^{-\mathrm{T}}(\vom)\right]_{ij} \, , & \mathrm{for} \quad \, 1 \leq i, j \leq 3 \, , \quad k = 4, \\
		\left[\pd{\Log_{\SO(3)}}{\vA}(\vA)\right]_{(i-3)jk} \, , & \mathrm{for} \quad \, 4 \leq i \leq 6 \, , \quad 1 \leq j, k \leq 3, \\
	\end{cases}
\end{equation}
where internally $\vom = \Log_{\SO(3)}(\vA)$ is used. Note that again the first order approximation of the required $\SO(3)$ quantities have to be used whenever $\omega \leq \omega_\mathrm{crit}$.

Finally, the derivative of the element-wise $\SE(3)$-interpolation~\eqref{eq:relative_interpolation_H_two_node} with respect to the generalized coordinates $\vq$ can be computed as
\begin{equation}
	\left[ \pd{\vH_{\mathcal{I}\mathcal{K}}}{\vq}(\xi, \vq) \right]_{ijk} = 
	\begin{multlined}[t]
		\sum_{e=0}^{n_\mathrm{el} - 1} \chi_{\mathcal{J}^e}(\xi)
		\Bigg\{
		\Bigg[\pd{\vH_{\mathcal{I}\mathcal{K}_e}}{\vq}\Bigg]_{ilk} \left[ \Exp_{\SE(3)} \left(N^e_1(\xi) \, \vth_{\mathcal{K}_e\mathcal{K}_{e+1}}\right) \right]_{lj} \\
		+ \left[\vH_{\mathcal{I}\mathcal{K}_e}\right]_{il} \left[\pd{\Exp_{ \SE(3)}}{\vth} \left(N^e_1(\xi) \, \vth_{\mathcal{K}_e\mathcal{K}_{e+1}}\right) \right]_{ljm} N^e_1(\xi) \left[\pd{\Log_{ \SE(3)}}{\vH} (\vH_{\mathcal{K}_e\mathcal{K}_{e+1}}) \right]_{mno} \\
		\Bigg(\Bigg[\pd{\vH^{-1}_{\mathcal{I}\mathcal{K}_e}}{\vq}\Bigg]_{npk} \left[\vH_{\mathcal{I}\mathcal{K}_{e+1}}\right]_{po} + \left[\vH^{-1}_{\mathcal{I}\mathcal{K}_e}\right]_{np} \left[\pd{\vH_{\mathcal{I}\mathcal{K}_{e+1}}}{\vq}\right]_{pok}\Bigg)
		\Bigg\} \, .
	\end{multlined}
\end{equation}
\correctb{\section{Discrete conservation properties}
	\label{sec:discrete_preservation_properties}
	
	Following Romero and Armero~\cite{Romero2002}, this section discusses the discrete conservation properties of the proposed rod finite element formulation. That is, the conservation of total energy, linear and angular momentum. For the sake of brevity, but without loss of generality, we consider the case for which the centerline points ${}_I \vr_{OP}(\xi, t)$ coincide with the cross-sections' center of mass and for which $\vS_{\rho_0}$ of \eqref{eq:A_S_I_rho} vanishes.
	\subsection{Discrete conservation of total energy}
	With the discrete kinematics from \eqref{eq:relative_interpolation_H_two_node} and \eqref{eq:var_vel_interpolation}, the total energy $E = T + U$ of the discretized rod is given in terms of the kinetic energy
	\begin{equation}
		T(\vq, \vu) = \frac{1}{2} \int_{\mathcal{B}} {}_I \dot{\vr}_{OQ}\T {}_I \dot{\vr}_{OQ} \, \diff[m] = \frac{1}{2} \sum_{e=0}^{n_\mathrm{el} - 1}\int_{\mathcal{J}^e} \left\{
		{}_I \vv_{P}\T A_{\rho_0} {}_I \vv_{P} + {}_K \vom_{IK}\T {}_K \vI_{\rho_0} {}_K \vom_{IK}
		\right\} J \, \diff[\xi]
	\end{equation}
	and the potential energy
	\begin{equation}
		U(\vq) = \sum_{e=0}^{n_\mathrm{el} - 1}\int_{\mathcal{J}^e} W({}_K \vga, {}_K \vka_{IK}; \xi) J \, \diff[\xi] \, .
	\end{equation}
	Their respective time derivatives are
	\begin{equation}
		\dot{T}(\vq, \vu) = \sum_{e=0}^{n_\mathrm{el} - 1}\int_{\mathcal{J}^e} \Big\{
		{}_I \vv_{P}\T A_{\rho_0} {}_I \dot{\vv}_{P}
		+ {}_K \vom_{IK}\T {}_K \vI_{\rho_0} ({}_K \vom_{IK})^{\bigcdot}
		\Big\} J \, \diff[\xi]
	\end{equation}
	and
	\begin{equation}
		\begin{aligned}
			\dot{U}(\vq) &= \sum_{e=0}^{n_\mathrm{el} - 1}\int_{\mathcal{J}^e} \left\{
			{}_K \vn\T ({}_K \bar{\vga})^{\bigcdot} + {}_K \vm\T ({}_K \bar{\vka}_{IK})^{\bigcdot}
			\right\} \diff[\xi] \\
			&= \sum_{e=0}^{n_\mathrm{el} - 1}\int_{\mathcal{J}^e} \left\{ 
			({}_I \vv_{P})_{,\xi}\T \, \vA_{IK} {}_K\vn + ({}_K \vom_{IK})_{,\xi}\T \, {}_K \vm - {}_K \vom_{IK}\T ({}_K \bar{\vga} \times {}_K \vn + {}_K \bar{\vka}_{IK} \times {}_K \vm)
			\right\} \diff[\xi] \, ,
		\end{aligned}
	\end{equation}
	where we have used the identities $({}_K \bar{\vga})^{\bigcdot} = \vA_{IK}\T ({}_I \vv_{P})_{,\xi} - {}_K \vom_{IK} \times {}_K \bar{\vga}$ and $({}_K \bar{\vka}_{IK})^{\bigcdot} = ({}_K \vom_{IK})_{,\xi} - {}_K \vom_{IK} \times {}_K \bar{\vka}_{IK}$, which follow from the derivation given in Appendix \ref{sec:variations_strain_measures} by replacing the variations with the time derivatives.

	Since the principle of virtual work holds for arbitrary virtual displacements, we can choose the nodal virtual displacements $\delta \vs^i = ({}_I \vv_{P_i}, \ {}_{K_i} \vom_{IK_i})$ to be composed of the nodal minimal velocities. Inserting these virtual displacements
	into~\eqref{eq:discretized_inertial_virtual_work} and~\eqref{eq:discretized_internal_virtual_work}, in the absence of any external force contributions, the virtual work principle leads to
	\begin{equation}
		0 = \delta W^\mathrm{dyn}(\delta \vs; \vq, \vu) + \delta W^\mathrm{int}(\delta \vs; \vq) = -\big(T(\vq, \vu) + U(\vq)\big)^{\bigcdot} \, .
	\end{equation}
	Consequently, conservation of the total energy is guaranteed also for the discrete equations.
	\subsection{Discrete conservation of linear momentum}
	The discrete linear momentum of the Cosserat rod is
	\begin{equation}
		{}_I \vL(\vu) = \int_{\mathcal{B}}  {}_I \dot{\vr}_{OQ} \diff m = \sum_{e=0}^{n_\mathrm{el} - 1} \int_{\mathcal{J}^e} A_{\rho_0} {}_I \vv_{P} J \, \diff[\xi] \, .
	\end{equation}
	For an arbitrary ${}_I \vc \in \mR^3$, let the nodal virtual displacements be given by
	\begin{equation}\label{eq:nodal_variation_linear_momentum}
		\delta \vs^i = \begin{pmatrix}
			{}_I \vc \\
			\mathbf{0}_{3 \times 1}
		\end{pmatrix} \, .
	\end{equation}
	This choice represents a virtual translation of the rod. Inserting \eqref{eq:nodal_variation_linear_momentum} into the discrete internal virtual work~\eqref{eq:discretized_internal_virtual_work} and recognizing that $N_{0,\xi}^e(\xi) + N_{1,\xi}^e(\xi) = 0$, one readily sees that $	\delta W^\mathrm{int}(\delta \vs; \vq) = 0$.
	Using \eqref{eq:nodal_variation_linear_momentum} in the discrete inertial virtual work~\eqref{eq:discretized_inertial_virtual_work} together with $N_0^e(\xi) + N_1^e(\xi) = 1$ leads to
	\begin{equation}
		\delta W^\mathrm{dyn}(\delta \vs;\vq,\vu) = - {}_I \vc\T \left(
		\sum_{e=0}^{n_\mathrm{el} - 1} \int_{\mathcal{J}^e} A_{\rho_0} {}_I \vv_{P} J \, \diff[\xi]
		\right)^{\bigcdot} = - {}_I \vc\T {}_I \dot{\vL}(\vu) \stackrel{!}{=} 0 \, ,
	\end{equation}
	which, in agreement with the principle of virtual work, must vanish if no external forces are applied to the rod. Consequently, the rate of change of the linear momentum ${}_I \dot{\vL} = \mathbf{0}_{3 \times 1}$ vanishes and the linear momentum ${}_I {\vL}$ is preserved.
	\subsection{Discrete conservation of angular momentum}
	The present rod finite element formulation uses nodal virtual rotations expressed in the individual cross-section-fixed bases ${}_{K_i} \delta \vph_{IK_i}(t) \in \mR^3$. For an arbitrary ${}_I \vet \in \mR^3$, the nodal virtual displacements $\delta \vs^i = ({}_I \vet \times {}_I \vr_{OP_i}, \ \vA_{IK_i}\T {}_I \vet)$
	lead to a virtual rotation around the origin of all nodal points. With respect to the inertial basis, by construction, the nodal cross-sections are all virtually rotated the same. However, due to the chosen interpolation \eqref{eq:var_vel_interpolation}$_3$, the cross-sections within the elements are not all rotated with the same constant virtual rotation ${}_I \vet$. Hence, for this formulation a virtual rotation of the entire rod is not contained in the set of admissible virtual displacements. Consequently, there is no algorithmic access to the discrete angular momentum and it is not possible to construct a numerical time integration scheme that preserves the discrete total angular momentum.
	
	If such a conservation property is crucial for a specific application, the present rod formulation can easily be modified. Instead of using the nodal virtual rotations ${}_{K_i} \delta \vph_{IK_i}(t) \in \mR^3$ expressed in the cross-section-fixed basis $K_i$, the nodal virtual rotations ${}_{I} \delta \vph_{IK_i}(t) \in \mR^3$ expressed in the inertial basis $I$ are used. Analogously, the nodal angular velocities ${}_{K_i} \vom_{IK_i}(t) \in \mR^3$ are replaced by ${}_{I} \vom_{IK_i}(t) \in \mR^3$. Proceeding similar to the original formulation, new discretized internal and inertial virtual work contributions are obtained. They are
	\begin{equation}\label{eq:discretized_internal_virtual_work_I_frame}
		\begin{aligned}
			\delta W^\mathrm{int}(\delta \vs;\vq) &= \delta \vs\T \vf^{\mathrm{int}}(\vq) \, , \quad \vf^{\mathrm{int}}(\vq) = \sum_{e=0}^{n_\mathrm{el} - 1} \vC_{e}\T \vf^{\mathrm{int}}_e(\vC_{e}\vq) \, , \\
			\vf^{\mathrm{int}}_e(\vq^e) &= -\int_{\mathcal{J}^e} \sum_{i=0}^{1}\Big\{ N^e_{i,\xi} \vC_{\vr, i}\T \vA_{IK} {}_K \vn + N^e_{i,\xi} \vC_{\vps, i}\T \vA_{IK} {}_K \vm 
			-N^e_{i} \vC_{\vps, i}\T \left({}_I \vr_{OP, \xi} \times \vA_{IK} {}_K \vn\right) \Big\} \diff[\xi] \, ,
		\end{aligned}
	\end{equation}
	and
	\begin{equation}\label{eq:discretized_inertial_virtual_work_I_frame}
		\begin{aligned}
			\delta W^\mathrm{dyn}(\delta \vs;\vq,\vu) &= -\delta \vs\T \left\{\vM(\vq) \dot{\vu} + \vf^{\mathrm{gyr}}(\vq, \vu) \right\} \, , \\
			\vM(\vq) &= \sum_{e=0}^{n_\mathrm{el} - 1} \vC_e\T \vM_e(\vC_{e}\vq) \vC_e \, , \quad
			\vf^{\mathrm{gyr}}(\vq, \vu) = \sum_{e=0}^{n_\mathrm{el}-1} \vC_e\T \vf^{\mathrm{gyr}}(\vC_{e} \vq, \vC_{e} \vu) \, , \\
			\vM_e(\vq^e) &= 
			\begin{multlined}[t]
				\int_{\mathcal{J}^e} \sum_{i=0}^{1} \sum_{k=0}^{1} N^e_{i} N^e_{k} \Big\{
				\vC_{\vr, i}\T A_{\rho_0} \mathbf{1}_{3\times3} \vC_{\vr, k} 
				+ \vC_{\vps, i}\T {}_I \vI_{\rho_0} \vC_{\vps, k}
				\Big\} J \diff[\xi] \, , 
			\end{multlined} \\
			\vf^{\mathrm{gyr}}_e(\vq^e, \vu^e) &= \int_{\mathcal{J}^e} \sum_{i=0}^{1} N^e_i \vC_{\vps, i}\T {}_I \widetilde{\vom}_{IK} {}_I \vI_{\rho_0} {}_I \vom_{IK} J \diff[\xi] \, ,
		\end{aligned}
	\end{equation}
	with
	\begin{equation}
		{}_I \vI_{\rho_0}(\vq) = \vA_{IK}(\vq) {}_K \vI_{\rho_0} \vA_{IK}(\vq)\T \, .
	\end{equation}
	Note, with this reformulation the mass matrix $\vM(\vq)$ and the gyroscopic forces $\vf^\mathrm{gyr}(\vq, \vu)$ depend on the generalized coordinates $\vq$, which is clearly a drawback with respect to the original formulation. 
	
	Conservation of energy and linear momentum follows in the same way as just shown and is not repeated here. As the virtual rotations are now formulated with respect to the inertial basis, the nodal virtual displacements
	\begin{equation}\label{eq:nodal_variation_angular_momentum}
		\delta \vs^i = \begin{pmatrix}
			{}_I \vet \times {}_I \vr_{OP_i} \\
			{}_I \vet
		\end{pmatrix} \, ,
	\end{equation} 
	induce a virtual rotation for the entire rod and conservation of total angular momentum can be shown. Indeed, for this discretization, the discrete angular momentum takes the form
	\begin{equation}
		{}_I \vJ(\vq, \vu) = \int_{\mathcal{B}}  {}_I \vr_{OQ} \times {}_I \dot{\vr}_{OQ} \diff m =  \sum_{e=0}^{n_\mathrm{el} - 1} \int_{\mathcal{J}^e} \left\{
		{}_I \vr_{OP} \times A_{\rho_0} {}_I \vv_{P} 
		+ {}_I \vI_{\rho_0} {}_I \vom_{IK}
		\right\} J \, \diff[\xi] \, .
	\end{equation}
	Inserting \eqref{eq:nodal_variation_angular_momentum} into the discrete internal virtual work~\eqref{eq:discretized_internal_virtual_work_I_frame} and recognizing the identities already used for the linear momentum leads to a vanishing internal virtual work functional $\delta W^\mathrm{int}(\delta \vs; \vq) = 0$.
	Using~\eqref{eq:nodal_variation_angular_momentum} in the discrete inertial virtual work~\eqref{eq:discretized_inertial_virtual_work_I_frame} results in
	\begin{equation}
		\delta W^\mathrm{dyn}(\delta \vs;\vq,\vu) = - {}_I \vet\T \sum_{e=0}^{n_\mathrm{el} - 1} \int_{\mathcal{J}^e} \left\{
		{}_I \vr_{OP} \times A_{\rho_0} {}_I \dot{\vv}_{P} 
		+ {}_I \vI_{\rho_0} {}_I \dot{\vom}_{IK}
		+ {}_I \widetilde{\vom}_{IK} {}_I \vI_{\rho_0} {}_I \vom_{IK}
		\right\} J \, \diff[\xi] = - {}_I \vet\T {}_I \dot{\vJ}(\vq, \vu) \stackrel{!}{=} 0 \, ,
	\end{equation}
	which, according to the principle of virtual work, must be zero for vanishing external virtual work functionals. It readily follows that ${}_I \dot{\vJ} = \mathbf{0}_{3 \times 1}$ and that the total angular momentum ${}_I \vJ$ is conserved in these cases.}

\bibliographystyle{ieeetr}


\begin{thebibliography}{10}
	
	\bibitem{Cosserat1909}
	E.~Cosserat and F.~Cosserat, {\em Th\'eorie des Corps d\'eformables}.
	\newblock Librairie Scientifique A. Hermann et Fils, 1909.
	
	\bibitem{Timoshenko1921}
	S.~Timoshenko, ``Lxvi. on the correction for shear of the differential equation
	for transverse vibrations of prismatic bars,'' {\em The London, Edinburgh,
		and Dublin Philosophical Magazine and Journal of Science}, vol.~41, no.~245,
	pp.~744--746, 1921.
	
	\bibitem{Reissner1981}
	E.~Reissner, ``On finite deformations of space-curved beams,'' {\em Zeitschrift
		f{\"u}r angewandte Mathematik und Physik ZAMP}, vol.~32, pp.~734--744, Nov
	1981.
	
	\bibitem{Simo1986}
	J.~C. Simo and L.~Vu-Quoc, ``A three-dimensional finite-strain rod model.
	{P}art {II}: {C}omputational aspects,'' {\em Computer Methods in Applied
		Mechanics and Engineering}, vol.~58, pp.~79--116, 1986.
	
	\bibitem{Antman2005}
	S.~S. Antman, {\em Nonlinear Problems of Elasticity}, vol.~107 of {\em Applied
		Mathematical Sciences}.
	\newblock Springer, 2nd~ed., 2005.
	
	\bibitem{Betsch2002}
	P.~Betsch and P.~Steinmann, ``Frame-indifferent beam finite elements based upon
	the geometrically exact beam theory,'' {\em International Journal for
		Numerical Methods in Engineering}, vol.~54, no.~12, pp.~1775--1788, 2002.
	
	\bibitem{Reddy2020}
	J.~N. Reddy and A.~R. Srinivasa, ``Misattributions and misnomers in mechanics:
	{W}hy they matter in the search for insight and precision of thought,'' {\em
		Vietnam Journal of Mechanics}, vol.~42, no.~3, pp.~283--291, 2020.
	
	\bibitem{Eugster2020b}
	S.~R. Eugster and J.~Harsch, {\em A variational formulation of classical
		nonlinear beam theories}, pp.~95--121.
	\newblock Cham: Springer International Publishing, 2020.
	
	\bibitem{Harsch2021a}
	J.~Harsch, G.~Capobianco, and S.~R. Eugster, ``Finite element formulations for
	constrained spatial nonlinear beam theories,'' {\em Mathematics and Mechanics
		of Solids}, vol.~26, no.~12, pp.~1838--1863, 2021.
	
	\bibitem{Meier2019}
	C.~Meier, A.~Popp, and W.~A. Wall, ``{G}eometrically {E}xact {F}inite {E}lement
	{F}ormulations for {S}lender {B}eams: {K}irchhoff--{L}ove {T}heory {V}ersus
	{S}imo--{R}eissner {T}heory,'' {\em Archives of Computational Methods in
		Engineering}, vol.~26, pp.~163--243, Jan 2019.
	
	\bibitem{Cardona1988}
	A.~Cardona and M.~Geradin, ``A beam finite element non-linear theory with
	finite rotations,'' {\em International Journal for Numerical Methods in
		Engineering}, vol.~26, no.~11, pp.~2403--2438, 1988.
	
	\bibitem{Ibrahimbegovic1995}
	A.~Ibrahimbegović, F.~Frey, and I.~Kožar, ``Computational aspects of
	vector-like parametrization of three-dimensional finite rotations,'' {\em
		International Journal for Numerical Methods in Engineering}, vol.~38, no.~21,
	pp.~3653--3673, 1995.
	
	\bibitem{Crisfield1999}
	M.~A. Crisfield and G.~Jeleni{\'c}, ``Objectivity of strain measures in the
	geometrically exact three-dimensional beam theory and its finite-element
	implementation,'' {\em Proceedings: Mathematical, Physical and Engineering
		Sciences}, vol.~455, no.~1983, pp.~1125--1147, 1999.
	
	\bibitem{Jelenic1999}
	G.~Jelenić and M.~Crisfield, ``{G}eometrically exact {3D} beam theory:
	implementation of a strain-invariant finite element for statics and
	dynamics,'' {\em Computer Methods in Applied Mechanics and Engineering},
	vol.~171, no.~1, pp.~141--171, 1999.
	
	\bibitem{Petrov1940}
	G.~I. Petrov, ``Application of {G}alerkin's method to the problem of stability
	of flow of a viscous fluid,'' {\em J. Appl. Math. Mech., Prikladnaya
		Matematika i Mekhanika, PMM}, vol.~4(3), pp.~3--12, 1994.
	
	\bibitem{Sailer2021}
	S.~Sailer and R.~I. Leine, ``Model reduction of the tippedisk: a path to the
	full analysis,'' {\em Nonlinear Dynamics}, vol.~105, pp.~1955--1975, Aug
	2021.
	
	\bibitem{Borri1994}
	M.~Borri and C.~Bottasso, ``An intrinsic beam model based on a helicoidal
	approximation--{P}art {I}: {F}ormulation,'' {\em International Journal for
		Numerical Methods in Engineering}, vol.~37, no.~13, pp.~2267--2289, 1994.
	
	\bibitem{Cesarek2013}
	P.~Češarek, M.~Saje, and D.~Zupan, ``Dynamics of flexible beams:
	{F}inite-element formulation based on interpolation of strain measures,''
	{\em Finite Elements in Analysis and Design}, vol.~72, pp.~47--63, 2013.
	
	\bibitem{Renda2016}
	F.~Renda, V.~Cacucciolo, J.~Dias, and L.~Seneviratne, ``{D}iscrete {C}osserat
	approach for soft robot dynamics: {A} new piece-wise constant strain model
	with torsion and shears,'' in {\em 2016 IEEE/RSJ International Conference on
		Intelligent Robots and Systems (IROS)}, pp.~5495--5502, 2016.
	
	\bibitem{Sonneville2014}
	V.~Sonneville, A.~Cardona, and O.~Brüls, ``Geometrically exact beam finite
	element formulated on the special {E}uclidean group {SE}(3),'' {\em Computer
		Methods in Applied Mechanics and Engineering}, vol.~268, pp.~451--474, 2014.
	
	\bibitem{Eugster2014c}
	S.~R. Eugster, C.~Hesch, P.~Betsch, and {\relax Ch}.~Glocker, ``Director-based
	beam finite elements relying on the geometrically exact beam theory
	formulated in skew coordinates,'' {\em International Journal for Numerical
		Methods in Engineering}, vol.~97, no.~2, pp.~111--129, 2014.
	
	\bibitem{Bullo1995}
	F.~Bullo and R.~M. Murray, ``{P}roportional {D}erivative ({PD}) {C}ontrol on
	the {E}uclidean {G}roup.'' 1995.
	
	\bibitem{Park2005}
	J.~Park and W.-K. Chung, ``Geometric integration on {E}uclidean group with
	application to articulated multibody systems,'' {\em IEEE Transactions on
		Robotics}, vol.~21, no.~5, pp.~850--863, 2005.
	
	\bibitem{Barfoot2014}
	T.~D. Barfoot and P.~T. Furgale, ``Associating uncertainty with
	three-dimensional poses for use in estimation problems,'' {\em IEEE
		Transactions on Robotics}, vol.~30, no.~3, pp.~679--693, 2014.
	
	\bibitem{Arnold2016}
	M.~Arnold, A.~Cardona, and O.~Br{\"u}ls, {\em A Lie Algebra Approach to Lie
		Group Time Integration of Constrained Systems}, pp.~91--158.
	\newblock Cham: Springer International Publishing, 2016.
	
	\bibitem{Harsch2020a}
	J.~Harsch and S.~R. Eugster, {\em Finite element analysis of planar nonlinear
		classical beam theories}, pp.~123--157.
	\newblock Cham: Springer International Publishing, 2020.
	
	\bibitem{Balobanov2018}
	V.~Balobanov and J.~Niiranen, ``Locking-free variational formulations and
	isogeometric analysis for the timoshenko beam models of strain gradient and
	classical elasticity,'' {\em Computer Methods in Applied Mechanics and
		Engineering}, vol.~339, pp.~137--159, 2018.
	
	\bibitem{Meier2015}
	C.~Meier, A.~Popp, and W.~A. Wall, ``A locking-free finite element formulation
	and reduced models for geometrically exact {K}irchhoff rods,'' {\em Computer
		Methods in Applied Mechanics and Engineering}, vol.~290, pp.~314--341, 2015.
	
	\bibitem{Greco2017}
	L.~Greco, M.~Cuomo, L.~Contrafatto, and S.~Gazzo, ``An efficient blended mixed
	b-spline formulation for removing membrane locking in plane curved kirchhoff
	rods,'' {\em Computer Methods in Applied Mechanics and Engineering},
	vol.~324, pp.~476--511, 2017.
	
	\bibitem{Santos2010}
	H.~Santos, P.~Pimenta, and J.~{Moitinho de Almeida}, ``Hybrid and multi-field
	variational principles for geometrically exact three-dimensional beams,''
	{\em International Journal of Non-Linear Mechanics}, vol.~45, no.~8,
	pp.~809--820, 2010.
	
	\bibitem{Santos2011}
	H.~Santos, P.~Pimenta, and J.~{Moitinho de Almeida}, ``A hybrid-mixed finite
	element formulation for the geometrically exact analysis of three-dimensional
	framed structures,'' {\em Computational Mechanics}, vol.~48, p.~591, Jun
	2011.
	
	\bibitem{Betsch2016}
	P.~Betsch and A.~Janz, ``An energy–momentum consistent method for transient
	simulations with mixed finite elements developed in the framework of
	geometrically exact shells,'' {\em International Journal for Numerical
		Methods in Engineering}, vol.~108, no.~5, pp.~423--455, 2016.
	
	\bibitem{Egeland2002}
	O.~Egeland and J.~T. Gravdahl, {\em Modeling and Simulation for Automatic
		Control}.
	\newblock Marine Cybernetics, 2002.
	
	\bibitem{Bosten2022}
	A.~Bosten, A.~Cosimo, J.~Linn, and O.~Br{\"u}ls, ``A mortar formulation for
	frictionless line-to-line beam contact,'' {\em Multibody System Dynamics},
	vol.~54, pp.~31--52, Jan 2022.
	
	\bibitem{Geradin2001}
	M.~G{\'e}radin and A.~Cardona, {\em Flexible Multibody Dynamics: A Finite
		Element Approach}.
	\newblock Wiley, 2001.
	
	\bibitem{dellIsola2020b}
	F.~dell'Isola and D.~J. Steigmann, {\em Discrete and {C}ontinuum {M}odels for
		{C}omplex {M}etamaterials}.
	\newblock Cambridge: Cambridge University Press, 2020.
	
	\bibitem{Hairer1993}
	E.~Hairer, S.~P. N\o{}rsett, and G.~Wanner, {\em Solving Ordinary Differential
		Equations I}.
	\newblock Springer, 1993.
	
	\bibitem{Jay1995}
	L.~Jay, ``Convergence of runge-kutta methods for differential-algebraic systems
	of index 3,'' {\em Applied Numerical Mathematics}, vol.~17, no.~2,
	pp.~97--118, 1995.
	
	\bibitem{Hairer2002}
	E.~Hairer and G.~Wanner, {\em Solving Ordinary Differential Equations II}.
	\newblock Springer, second~ed., 2002.
	
	\bibitem{Jay1996}
	L.~Jay, ``Symplectic partitioned runge–kutta methods for constrained
	hamiltonian systems,'' {\em SIAM Journal on Numerical Analysis}, vol.~33,
	no.~1, pp.~368--387, 1996.
	
	\bibitem{Hairer2006}
	E.~Hairer, C.~Lubich, and G.~Wanner, {\em Geometric numerical integration:
		structure-preserving algorithms for ordinary differential equations}, vol.~31
	of {\em Springer Series in Computational Mathematics}.
	\newblock Springer Verlag, Berlin Heidelberg, second~ed., 2006.
	
	\bibitem{Bruels2012}
	O.~Brüls, A.~Cardona, and M.~Arnold, ``Lie group generalized-$\alpha$ time
	integration of constrained flexible multibody systems,'' {\em Mechanism and
		Machine Theory}, vol.~48, pp.~121--137, 2012.
	
	\bibitem{Newmark1959}
	N.~M. Newmark, ``A method of computation for structural dynamics,'' {\em
		Journal of the Engineering Mechanics Division}, vol.~85, no.~3, pp.~67--94,
	1959.
	
	\bibitem{Chung1993}
	J.~Chung and G.~Hulbert, ``A time integration algorithm for structural dynamics
	with improved numerical dissipation: the generalized-$\alpha$ method,'' {\em
		Journal of applied mechanics}, vol.~60, no.~2, pp.~371--375, 1993.
	
	\bibitem{Jansen2000}
	K.~E. Jansen, C.~H. Whiting, and G.~M. Hulbert, ``A generalized-$\alpha$ method
	for integrating the filtered navier–stokes equations with a stabilized
	finite element method,'' {\em Computer Methods in Applied Mechanics and
		Engineering}, vol.~190, no.~3, pp.~305--319, 2000.
	
	\bibitem{Ibrahimbegovic1997}
	A.~Ibrahimbegovic, ``On the choice of finite rotation parameters,'' {\em
		Computer Methods in Applied Mechanics and Engineering}, vol.~149, no.~1,
	pp.~49--71, 1997.
	\newblock Containing papers presented at the Symposium on Advances in
	Computational Mechanics.
	
	\bibitem{Huynh2009}
	D.~Q. Huynh, ``Metrics for {3D} {R}otations: {C}omparison and {A}nalysis,''
	{\em Journal of Mathematical Imaging and Vision}, vol.~35, pp.~155--164, Oct
	2009.
	
	\bibitem{Ogden1997}
	R.~W. Ogden, {\em Non-linear Elastic Deformations}.
	\newblock Dover Publications, 1997.
	
	\bibitem{Meier2014}
	C.~Meier, A.~Popp, and W.~A. Wall, ``An objective {3D} large deformation finite
	element formulation for geometrically exact curved {K}irchhoff rods,'' {\em
		Computer Methods in Applied Mechanics and Engineering}, vol.~278,
	pp.~445--478, 2014.
	
	\bibitem{Maekinen2007}
	J.~Mäkinen, ``{T}otal {L}agrangian {R}eissner's geometrically exact beam
	element without singularities,'' {\em International Journal for Numerical
		Methods in Engineering}, vol.~70, no.~9, pp.~1009--1048, 2007.
	
	\bibitem{Magnus1971}
	K.~Magnus, {\em Kreisel: Theorie und Anwendungen}.
	\newblock Berlin, Heidelberg: Springer Berlin Heidelberg, 1~ed., 1971.
	
	\bibitem{Hall2015}
	B.~Hall, {\em Lie Groups, Lie Algebras, and Representations: An Elementary
		Introduction}.
	\newblock Cham: Springer International Publishing, 2015.
	
	\bibitem{Baker2002}
	A.~Baker, {\em Matrix Groups: An Introduction to Lie Group Theory}.
	\newblock Springer London, 2002.
	
	\bibitem{Iserles2000}
	A.~Iserles, H.~Z. Munthe-Kaas, S.~P. Nørsett, and A.~Zanna, ``Lie-group
	methods,'' {\em Acta Numerica}, vol.~9, p.~215–365, 2000.
	
	\bibitem{Murray1994}
	R.~M. Murray, Z.~Li, and S.~S. Sastry, {\em A {M}athematical {I}ntroduction to
		{R}obotic {M}anipulation}.
	\newblock Taylor \& Francis, 1994.
	
	\bibitem{Gallego2015}
	G.~Gallego and A.~Yezzi, ``{A} {C}ompact {F}ormula for the {D}erivative of a
	3-{D} {R}otation in {E}xponential {C}oordinates,'' {\em Journal of
		Mathematical Imaging and Vision}, vol.~51, pp.~378--384, Mar 2015.
	
	\bibitem{Romero2002}
	I.~Romero and F.~Armero, ``An objective finite element approximation of the
	kinematics of geometrically exact rods and its use in the formulation of an
	energy-momentum conserving scheme in dynamics,'' {\em International Journal
		for Numerical Methods in Engineering}, vol.~54, pp.~1683--1716, 2002.
	
\end{thebibliography}

\end{document}